\documentclass[12pt]{article}
\usepackage{graphicx}
\usepackage{epsfig}
\usepackage{amsmath}
\usepackage{amssymb}
\usepackage{amsthm}
\usepackage{latexsym}
\usepackage{cite}
\usepackage{amsbsy,amsfonts,euscript,eepic,rotating}

\newcommand{\comment}[1]{}
\newtheorem{theorem}{Theorem}
\newtheorem{lemma}{Lemma}[section]
\newtheorem{remark}{Remark}[section]
\newtheorem{corollary}{Corollary}[section]

\newtheorem{definition}{Definition}[section]
\begin{document}

\title{\LARGE
{\bf Critical Percolation Exploration Path \\ and $SLE_6$: a Proof of Convergence}
%On the Convergence of the Percolation Exploration Process to Chordal $SLE_6$}
}

\author{
{\bf Federico Camia}
\thanks{Research partially supported by a Marie Curie Intra-European Fellowship
under contract MEIF-CT-2003-500740 and by a Veni grant of the Dutch Organization
for Scientific Research (NWO).}\,
\thanks{E-mail: fede@few.vu.nl}\\
{\small \sl Department of Mathematics, Vrije Universiteit Amsterdam}\\
\and
{\bf Charles M.~Newman}
\thanks{Research partially supported by the
U.S. NSF under grants DMS-01-04278 and DMS-06-06696.}\,
\thanks{E-mail: newman@courant.nyu.edu}\\
{\small \sl Courant Inst.~of Mathematical Sciences, New York University}
}

\date{}

\maketitle

\begin{abstract}
It was argued by Schramm and Smirnov that the critical site percolation
exploration path on the triangular lattice converges in distribution
to the trace of chordal $SLE_6$.
We provide here a detailed proof, which relies on Smirnov's theorem
that crossing probabilities have a conformally invariant scaling limit
(given by Cardy's formula).
The version of convergence to $SLE_6$ that we prove suffices for the
Smirnov-Werner derivation of certain critical percolation crossing
exponents and for our analysis of the critical percolation full scaling
limit as a process of continuum nonsimple loops.
\end{abstract}

\noindent {\bf Keywords:} continuum scaling limit, percolation, SLE,
critical behavior, triangular lattice, conformal invariance.

\noindent {\bf AMS 2000 Subject Classification:} 82B27, 60K35, 82B43,
60D05, 30C35.

\section{Synopsis}

The purpose of this paper is to present a detailed, self-contained proof
of the convergence of the critical site percolation exploration path
(on the triangular lattice) to the trace of chordal $SLE_6$ for Jordan domains.
We will prove convergence in a strong sense: in the topology induced by
the uniform metric on continuous curves (modulo monotonic reparametrization),
and ``locally uniformly" in the ``shape" of the domain and the starting and
ending points of the curve.

The main technical difficulty (in the approach followed here) appears
in a rather surprising way --- to obtain a Markov property for any
scaling limit of the percolation exploration path.
The surprise is that an even stronger Markov property trivially holds
for the exploration path itself.
To show that an analogous property holds in the scaling limit
(see Theorem~\ref{spatial-markov} and Remark~\ref{explanation}) is
largely responsible for the length of the paper.
Roughly, the difficulty is that in the scaling limit the exploration path
touches itself and the boundary of the domain (infinitely many times).
The touching of the domain boundary in particular requires a lengthy analysis
(see Lemmas~\ref{double-crossing}-\ref{mushroom}) since the standard percolation
bound on multiple crossings of a ``semi-annulus" only applies to the case of
a ``flat" boundary (see~\cite{smirnov-long} and Appendix A of~\cite{lsw5}).
This issue is resolved here by using the continuity of Cardy's formula with
respect to changes in the domain.
We remark that the results and methods developed here about touching of domain
boundaries have other applications --- e.g., to the existence and conformal
invariance of the full scaling limit in general (non-flat) domains; these extensions
of the results of~\cite{cn1,cn2} will be discussed elsewhere~\cite{cn4}.

%Another important technical point is the need to strengthen Smirnov's
%theorem~\cite{smirnov} about the convergence of crossing probabilities to Cardy's
%formula (Theorem~\ref{cardy-smirnov} here) to obtain convergence ``locally uniformly"
%in the ``shape" of the domain and the positions of the two boundary arcs
%(see Theorem~\ref{strong-cardy}).

The proof has two parts: a characterization for $SLE_6$ curves
(Theorem~\ref{characterization}), which is similar to Schramm's argument identifying
$SLE_{\kappa}$ but only uses conformal invariance of hulls at special stopping times,
and a series of results showing that any subsequential scaling limit of the exploration
path satisfies the hypotheses of Theorem~\ref{characterization}.
%One of these results is the already mentioned theorem of Smirnov (Theorem~\ref{cardy-smirnov}).

The theorem of Smirnov (Theorem~\ref{cardy-smirnov} here) about convergence
to Cardy's formula~\cite{smirnov} is a key tool throughout.
It is used in the proof of Theorem~\ref{characterization}, which follows roughly
Smirnov's sketch in~\cite{smirnov,smirnov-long} (but with one significant
difference --- see Remarks~\ref{remark-differ} and~\ref{remark-difference}), and
is also crucial, in its strengthened version, Theorem~\ref{strong-cardy}, in proving
that the ``filling" of the exploration path converges to a hull process having the
Markov property necessary to apply Theorem~\ref{characterization}.
As mentioned, this step, implicitly assumed in Smirnov's sketch~\cite{smirnov,smirnov-long},
turns out to be the most technically difficult one.
Despite its length, we believe that a detailed proof is needed, since the result,
beside its own interest, has important applications --- notably the rigorous
derivation of certain critical exponents~\cite{sw}, of Watt's crossing
formula~\cite{dubedat} and Schramm's percolation formula~\cite{schramm1}, and
the derivation of the full scaling limit~\cite{cn1,cn2}.
We note that Smirnov has recently sketched in~\cite{smirnov06}
a proof different from that of~\cite{smirnov,smirnov-long}.

\section{Introduction} \label{intro}

The percolation exploration path was introduced by Schramm in 1999 in a
seminal paper~\cite{schramm} where
%about the continuum scaling limit of discrete
%models in two dimensions and their purposed conformal invariance.
it was used to give a precise formulation (beyond crossing probabilities)
to the conjecture that the scaling limit of two-dimensional critical
percolation is conformally invariant.
Schramm's formulation of the conjecture involves his Stochastic Loewner
Evolution or Schramm-Loewner Evolution
($SLE_{\kappa}$), and can be expressed, roughly speaking, by
saying that the percolation exploration path converges in distribution
to $SLE_6$.

A simple and elegant argument, due again to Schramm, shows that if the scaling
limit of the percolation exploration path exists and is conformally invariant,
then it must necessarily be an $SLE_{\kappa}$ curve; the value $\kappa=6$ can
be determined by looking at crossing probabilities, since $SLE_6$ is the only
$SLE_{\kappa}$ that satisfies Cardy's formula~\cite{cardy}.

Shortly after Schramm's paper appeared, Smirnov published a
proof~\cite{smirnov}, for site percolation on the triangular lattice, of the
conformal invariance of the scaling limit of crossing probabilities, opening
the way to a complete proof of Schramm's conjecture.
In~\cite{smirnov} (see also~\cite{smirnov-long}) Smirnov also outlined a
possible strategy for using the conformal invariance of crossing probabilities
to prove Schramm's conjecture.
Roughly at the same time, convergence of the exploration path to $SLE_6$
was used by Smirnov and Werner~\cite{sw} and by Lawler, Schramm and Werner~\cite{lsw5}
as a key step in a derivation of the values of various percolation
critical exponents, most of which had been previously predicted in the
physics literature (see the references in~\cite{sw}).
Later, it was used by the authors of this paper to obtain the full
scaling limit of two-dimensional critical percolation (see~\cite{cn,cn1}).

However, a detailed proof of the convergence of the exploration path to
$SLE_6$ did not appear until 2005, in an appendix of~\cite{cn1}, where
we followed a modified version of Smirnov's strategy.
The purpose of the present paper is to present essentially that
proof in a self-contained form.
(We note that Lemma A.3 in~\cite{cn1}, whose proof had an error,
has been replaced by Lemma~\ref{close-encounters} here.)
Our proof roughly follows Smirnov's outline of~\cite{smirnov,smirnov-long},
based on the convergence to Cardy's formula~\cite{smirnov,smirnov-long} and
on Markov properties (see Theorem~\ref{characterization} below and the
discussion preceding it).
But there are two significant modifications, which we found necessary for a proof.
The first is to use a different sequence of stopping times, which results
in a different geometry for the Markov chain approximation to $SLE_6$ (see
Remark~\ref{remark-difference}).
The second is that ``close encounters" by the exploration
path to the domain boundary are not handled by general results for
``three-arm" events at the boundary of a half-plane, but rather by
a more complex argument based partly on continuity of crossing probabilities
with respect to domain boundaries (see Lemmas~\ref{double-crossing}-\ref{mushroom}).
%Moreover, we cannot use directly Smirnov's result on convergence of
%crossing probabilities (see Theorem~\ref{cardy-smirnov}), but need
%an extended version which is given in Theorem~\ref{strong-cardy} of
%Section~\ref{convergence-to-sle}.

We note (see Remark~\ref{remark-differ} below for more discussion)
that our choice of stopping times is closer in spirit than is the
choice in~\cite{smirnov,smirnov-long} to the proofs of convergence
of the loop erased random walk to $SLE_2$~\cite{lsw6} and of the
harmonic explorer to $SLE_4$~\cite{ss}.
It may also be applicable to other systems in which an $SLE_{\kappa}$ limit
is expected, provided that sufficient information can be obtained about
conformal invariance of the scaling limit of the analogues of exploration
hitting distributions in those systems.

Schramm's conjecture, as stated in Smirnov's paper~\cite{smirnov}, concerns
the convergence in distribution of the percolation exploration path to the
trace of chordal $SLE_6$ in a fairly arbitrary fixed domain.
Here (see Theorem~\ref{conv-to-sle-1}) we will prove a version of convergence
which is slightly stronger but somewhat less general: we will show that the
distribution of the percolation exploration path converges to that of the
trace of chordal $SLE_6$ ``\emph{locally uniformly}" in the ``shape" of the
domain and in the positions of the starting and ending points of the path,
but we will restrict attention to Jordan domains (i.e., domains whose boundary
is a simple closed curve).
Our main motivation in this specific formulation is to provide the key tool
needed in~\cite{cn2} to prove that the set of \emph{all} critical percolation
interfaces converges (in distribution) in the scaling limit to a certain
countable collection of continuous, nonsimple, noncrossing, fractal loops in
the plane.
Our formulation of convergence to $SLE_6$ is also sufficient for a key step
in the proof of certain critical exponents~\cite{sw} -- namely for
$j (\geq 1)$ crossings of a semi-annulus and for $j (\geq 2)$ crossings, not
all of the same color, of an annulus.
It does not appear to be sufficient, without at least also using some of~\cite{cn2},
for the derivation in~\cite{lsw5} of the ``one-arm" exponent (i.e., for one crossing
of an annulus) and thus not sufficient for proofs of other exponents based
on the one-arm exponent (see~\cite{kesten87} and Sect.~1 of~\cite{sw}).

%The rest of the paper is organized as follows.
In the next section, we give some preliminary definitions.
%including that of chordal $SLE_{\kappa}$.
In Sect.~\ref{lapa}, we define the percolation exploration path.
In Sect.~\ref{convergence-to-sle}, we introduce Cardy's formula, give
a characterization result for $SLE_6$, and state an extended version of
Smirnov's result on the scaling limit of crossing probabilities.
Sect.~\ref{boundary} contains results concerning the ``envelope" of
the hull of exploration paths and $SLE_6$ paths.
Those results are needed in Sect.~\ref{convergence}, which is devoted
to the proof of the main convergence result (see Theorem~\ref{conv-to-sle-1}).
The paper ends with an appendix about sequences of conformal maps.

\section{Preliminary Definitions} \label{defs}

We identify the real plane ${\mathbb R}^2$ and the complex plane $\mathbb C$
and use
%(02-24-06) (i.e., the one-point compactification of
%${\mathbb C}$ FEDERICO: MAYBE WE SHOULD GET RID OF THE RIEMANN
%SPHERE SINCE WE NEVER MENTION IT AGAIN)
the open half-plane
$\mathbb H = \{ x+iy : y>0 \}$ (and its closure $\overline{\mathbb H}$).
%where chordal $SLE$ will be defined (see Section~\ref{sle1}).
$\mathbb D$ denotes the open unit disc
${\mathbb D} = \{ z \in {\mathbb C} : |z|<1 \}$.
A domain $D$ of the complex plane $\mathbb C$ is a
nonempty, connected, open subset of $\mathbb C$; a
simply connected domain $D$ is said to be Jordan
if its (topological) boundary $\partial D$
is a Jordan curve (i.e., a simple continuous loop).

We often use Riemann's mapping theorem --- that if $D$ is
any simply connected domain other than the entire plane
$\mathbb C$ and $z_0 \in D$, then there is a unique conformal
map $\phi$ of $\mathbb D$ onto $D$ such that $\phi(0)=z_0$
and $\phi'(0)>0$.

%\subsection{Compactification of ${\mathbb R}^2$}
When taking the scaling limit $\delta \to 0$ one can focus
on fixed bounded
regions, $\Lambda \subset {\mathbb R}^2$,
or consider the whole ${\mathbb R}^2$ at once.
The second option avoids dealing with boundary conditions,
but requires an appropriate choice of metric.
A convenient way of dealing with the whole ${\mathbb R}^2$
is to replace the Euclidean metric with a distance function
%$\Delta(\cdot,\cdot)$ defined on ${\mathbb R}^2 \times {\mathbb R}^2$ by
\begin{equation}
\Delta(u,v) = \inf_h \int (1 + | {h} |^2)^{-1} \, ds,
\end{equation}
where the infimum is over all smooth curves $h(s)$
joining $u$ with $v$, parametrized by arclength $s$, and
where $|\cdot|$ denotes the Euclidean norm.
This metric is equivalent to the Euclidean metric in bounded
regions, but it has the advantage of making ${\mathbb R}^2$
precompact.
Adding a single point at infinity yields the compact space
$\dot{\mathbb R}^2$ which is isometric, via stereographic
projection, to the two-dimensional sphere.

\subsection{The Space of Curves} \label{space}

In dealing with the scaling limit we use the approach of
Aizenman-Burchard~\cite{ab}.
Denote by ${\cal S}_{\Lambda}$ the complete separable metric
space of continuous curves in a closed (bounded) subset
$\Lambda \subset {\mathbb R}^2$
with the metric~(\ref{distance}) defined below.
Curves are regarded as equivalence classes of continuous
functions from the unit interval to ${\mathbb R}^2$,
modulo monotonic reparametrizations.
$\gamma$ will represent a particular curve and $\gamma(t)$ a
parametrization of $\gamma$; ${\cal F}$ will represent a set
of curves (more precisely, a closed subset of ${\cal S}_{\Lambda}$).
%$\text{d}(\cdot,\cdot)$ will denote the uniform metric
%on curves, defined by
We define a metric on curves by
\begin{equation} \label{distance}
\text{d} (\gamma_1,\gamma_2) \equiv \inf
\sup_{t \in [0,1]} |\gamma_1(t) - \gamma_2(t)|,
\end{equation}
where the infimum is over %all choices of
parametrizations of $\gamma_1$ and $\gamma_2$.
%from the interval $[0,1]$.
The distance between two closed sets of curves is defined
by the induced Hausdorff metric:
\begin{equation} \label{hausdorff}
\text{dist}({\cal F},{\cal F}') \leq \varepsilon
\Leftrightarrow \forall \, \gamma \in {\cal F}, \, \exists \,
\gamma' \in {\cal F}' \text{ with }
\text{d} (\gamma,\gamma') \leq \varepsilon,
\text{ and vice versa.}
\end{equation}
The space $\Omega_{\Lambda}$ of closed subsets of ${\cal S}_{\Lambda}$
(i.e., collections of curves in $\Lambda$) with the metric~(\ref{hausdorff})
is also a complete separable metric space.
We denote by ${\cal B}_{\Lambda}$ its Borel $\sigma$-algebra.
For each $\delta>0$, the random curves we consider are polygonal paths
% corresponding to exploration paths $\gamma^{\delta}$
on the edges of the hexagonal lattice $\delta {\cal H}$,
dual to the triangular lattice $\delta {\cal T}$.
A superscript $\delta$ indicates that the curves correspond to
a model with a ``short distance cutoff'' of magnitude $\delta$.
We also consider the complete separable metric space ${\cal S}$
of continuous curves in $\dot{\mathbb R}^2$ with distance
\begin{equation} \label{Distance}
\text{D} (\gamma_1,\gamma_2) \equiv \inf
\sup_{t \in [0,1]} \Delta(\gamma_1(t),\gamma_2(t)),
\end{equation}
where the infimum is again over %all choices of
parametrizations of $\gamma_1$ and $\gamma_2$. %from the interval $[0,1]$.
The distance between two closed sets of curves is again
defined by the induced Hausdorff metric:
\begin{equation} \label{hausdorff-D}
\text{Dist}({\cal F},{\cal F}') \leq \varepsilon
\Leftrightarrow \forall \, \gamma \in {\cal F}, \, \exists \,
\gamma' \in {\cal F}' \text{ with }
\text{D} (\gamma,\gamma') \leq \varepsilon,
\text{ and vice versa.}
\end{equation}
The space $\Omega$ of closed sets of $\cal S$
(i.e., collections of curves in $\dot{\mathbb R}^2$)
with the metric~(\ref{hausdorff-D}) is also a complete
separable metric space.
We denote by ${\cal B}$ its Borel $\sigma$-algebra.
When we talk about convergence in distribution of random curves,
we %always mean with respect to
refer to the uniform metric~(\ref{distance}),
while %when we deal with
for closed collections of curves, we refer to the
metric~(\ref{hausdorff}) or~(\ref{hausdorff-D}).
\begin{remark} \label{ab}
In~\cite{cn,cn2}, the space $\Omega$ of closed sets of $\cal S$
was used for collections of exploration paths and cluster
boundary loops and their scaling limits, $SLE_6$ paths and
continuum nonsimple loops.
Here, in the statements and proofs of
Lemmas~\ref{double-crossing}, \ref{equal} and~\ref{mushroom},
we apply $\Omega$ in essentially the original setting of
Aizenman and Burchard~\cite{aizenman,ab}, i.e., for collections
of blue and yellow simple $\cal T$-paths (see Sect.~\ref{lapa})
and their scaling limits.
The slight modification needed to keep track of the colors is
easily managed.
\end{remark}

\subsection{Chordal $SLE_{\kappa}$} \label{sle1}

The Stochastic Loewner Evolution ($SLE_{\kappa}$) was introduced
by Schramm~\cite{schramm} to study two-dimensional probabilistic
lattice models whose scaling limits are expected to be conformally
invariant.
Here we describe the chordal version of $SLE_{\kappa}$;
for more, see~\cite{schramm} as well as the fine reviews by
Lawler~\cite{lawler1}, Kager and Nienhuis~\cite{kn}, and
Werner~\cite{werner4}, and Lawler's book~\cite{lawler2}.

Let $\mathbb H$ denote the upper half-plane.
%For a given continuous real function $U_t$ with $U_0 = 0$,
%define, for each $z \in \overline{\mathbb H}$, the function
%$g_t(z)$ as the solution to the ODE
%\begin{equation}
%\partial_t g_t(z) = \frac{2}{g_t(z) - U_t},
%\end{equation}
%with $g_0(z) = z$.
%This is well defined as long as $g_t(z) - U_t \neq 0$,
%i.e., for all $t < T(z)$, where
%\begin{equation}
%T(z) \equiv \sup \{ t \geq 0 : \min_{s \in [0,t]} | g_s(z) - U_s| > 0 \}.
%\end{equation}
%Let $K_t \equiv \{ z \in \overline{\mathbb H} : T(z) \leq t \}$
%and let ${\mathbb H}_t$ be the unbounded component of
%${\mathbb H} \setminus K_t$; it can be shown that $K_t$ is bounded
%and that $g_t$ is a conformal map from ${\mathbb H}_t$ onto $\mathbb H$.
%For each $t$, it is possible to write $g_t(z)$ as
%\begin{equation}
%g_t(z) = z + \frac{2t}{z} + o(\frac{1}{z}),
%\end{equation}
%when $z \to \infty$.
%The family $(K_t, t \geq 0)$ is called the {\bf Loewner chain}
%associated to the driving function $(U_t, t \geq 0)$.
%
%\begin{definition} \label{def-sle}
%{\bf Chordal $SLE_{\kappa}$} is the Loewner chain $(K_t, t \geq 0)$
%that is obtained when the driving function
%$U_t = \sqrt{\kappa} B_t$ is $\sqrt{\kappa}$ times a standard
%real-valued Brownian motion $(B_t, t \geq 0)$ with $B_0 = 0$.
%\end{definition}
%
For all $\kappa \geq 0$, chordal $SLE_{\kappa}$ in $\mathbb H$ is a certain
random family $(K_t,t \geq 0)$ of bounded subsets of $\overline{\mathbb H}$
that is generated by a continuous random curve $\gamma$ (with $\gamma(0)=0$)
in the sense that, for all $t \geq 0$,
${\mathbb H}_t \equiv {\mathbb H} \setminus K_t$ is the unbounded connected
component of ${\mathbb H} \setminus \gamma[0,t]$; $\gamma$ is called the
{\bf trace} of chordal $SLE_{\kappa}$.

%\subsection{Chordal $SLE_{\kappa}$ in an Arbitrary Simply Connected Domain} \label{sle2}

Let $D \subset {\mathbb C}$ ($D \neq {\mathbb C}$) be a simply
connected domain whose boundary is a continuous curve.
%By Riemann's mapping theorem, there are (many) conformal maps
%from the upper half-plane $\mathbb H$ onto $D$.
Given two distinct points $a,b \in \partial D$ (or more accurately,
two distinct prime ends), there exists a conformal map $f$ from
$\mathbb H$ onto $D$ such that $f(0)=a$ and
$f(\infty) \equiv \lim_{|z| \to \infty} f(z) = b$.
The choice of $a$ and $b$ only characterizes $f(\cdot)$ up to a
scaling factor $\lambda>0$, since $f(\lambda \cdot)$ would also do.

Suppose that $(K_t, t \geq 0)$ is a chordal $SLE_{\kappa}$ in
$\mathbb H$; chordal $SLE_{\kappa}$ $(\tilde K_t, t \geq 0)$ in $D$
from $a$ to $b$ as the image of $(K_t, t \geq 0)$ under $f$.
%It is possible to show, using scaling properties of
%$SLE_{\kappa}$, that
The law of $(\tilde K_t, t \geq 0)$ is unchanged, up to a linear
time-change, if we replace $f(\cdot)$ by $f(\lambda \cdot)$.
One considers $(\tilde K_t, t \geq 0)$ as a process from $a$
to $b$ in $D$, ignoring the role of $f$.

In the case $\kappa = 6$, $(K_t, t \geq 0)$ is generated by a
continuous nonsimple curve $\gamma$ with Hausdorff dimension $7/4$.
We will denote by $\gamma_{D,a,b}$ the image of $\gamma$ under $f$
and call it the trace of chordal $SLE_6$ in $D$ from $a$ to $b$;
$\gamma_{D,a,b}$ is a continuous nonsimple curve inside $\overline D$
from $a$ to $b$, and it can be given a parametrization $\gamma_{D,a,b}(t)$
such that $\gamma_{D,a,b}(0)=a$ and $\gamma_{D,a,b}(1)=b$, so that
we are in the metric framework described in Section~\ref{space}.
%It will be convenient to think of $\gamma_{D,a,b}$ as an
%oriented path, with orientation from $a$ to $b$.

\section{Lattices and Paths} \label{lapa}

We will denote by $\cal T$ the two-dimensional triangular lattice,
whose sites we think of as the elementary cells of a regular hexagonal
lattice $\cal H$ embedded in the plane as in Figure~\ref{fig2-sec4}.
We say that two hexagons are neighbors (or that they are adjacent) if
they have a common edge.
A sequence $(\xi_0, \ldots, \xi_n)$ of %sites of $\cal T$
hexagons of $\cal H$ such that $\xi_{i-1}$ and $\xi_i$
are neighbors for all $i= 1, \ldots, n$ and $\xi_i \neq \xi_j$
whenever $i \neq j$ will be called a {\bf $\cal T$-path} and denoted
by $\pi$.
If the first and last hexagons of the path are neighbors,
the path will be called a {\bf $\cal T$-loop}.
%The $\cal T$-paths we will consider are self-avoiding
%(i.e., $\xi_i \neq \xi_j$ whenever $i \neq j$).

A set $D$ of hexagons is {\bf connected} if any two hexagons in $D$
can be joined by a $\cal T$-path contained in $D$.
We say that a finite set $D$ of hexagons is {\bf simply connected}
if both $D$ and ${\cal T} \setminus D$ are connected.
%(by the edges of $\cal T$).
For a simply connected set $D$ of hexagons, we denote by
$\Delta D$ its {\bf external site boundary}, or {\bf s-boundary}
(i.e., the set of hexagons that do not belong to $D$
but are adjacent to hexagons in $D$), and by $\partial D$ the
topological boundary of $D$ when $D$ is considered as a domain of
$\mathbb C$.
We will call a bounded, simply connected subset $D$ of $\cal T$
a {\bf Jordan set} if its s-boundary $\Delta D$ is a $\cal T$-loop.

For a Jordan set $D \subset {\cal T}$, a vertex $x \in {\cal H}$
that belongs to $\partial D$ can be either of two types, according to whether
the edge incident on $x$ that is not in $\partial D$ belongs to a hexagon
in $D$ or not.
We call a vertex of the second type an {\bf e-vertex} (e for ``external''
or ``exposed'').

Given a Jordan set $D$ and two e-vertices $x,y$ in $\partial D$,
we denote by $\partial_{x,y} D$ the portion of $\partial D$
traversed counterclockwise from $x$ to $y$, and call it the
{\bf right boundary}; the remaining part of the boundary is
denote by $\partial_{y,x} D$ and is called the {\bf left boundary}.
Analogously, the portion of $\Delta_{x,y} D$ of $\Delta D$ whose
hexagons are adjacent to $\partial_{x,y} D$ is called the
{\bf right s-boundary} and the remaining part the {\bf left s-boundary}.

%\begin{figure}[!ht]
%\begin{center}
%\includegraphics[width=6cm]{fig1-sec4.ps}
%\caption{Portion of the hexagonal lattice.}
%\label{fig1-sec4}
%\end{center}
%\end{figure}

A {\bf percolation configuration}
$\sigma = \{ \sigma(\xi) \}_{\xi \in \cal T} \in \{ -1, +1 \}^{\cal T}$
on $\cal T$ is an assignment of $-1$ (equivalently, yellow) or $+1$
(blue) to each site of $\cal T$.
For a domain $D$ of the plane, the restriction to $D \cap \cal T$ of
%$\cal T$ of the percolation configuration
$\sigma$ is denoted by $\sigma_D$.
On the space of configurations $\Sigma = \{ -1,+1 \}^{\cal T}$,
we consider the usual product topology and denote by $\mathbb P$
the uniform measure, corresponding to Bernoulli percolation with
equal density of yellow (minus) and blue (plus) hexagons, which is
critical percolation in the case of the triangular lattice.

A (percolation) {\bf cluster} is a maximal, connected, monochromatic
subset of $\cal T$; we will distinguish between blue (plus) and
yellow (minus) clusters.
The {\bf boundary} of a cluster $D$ is the set of edges of
$\cal H$ that surround the cluster (i.e., its Peierls contour);
it coincides with the topological boundary of $D$ considered as a
domain of $\mathbb C$.
The set of all boundaries is a collection of ``nested'' simple loops
along the edges of $\cal H$.

Given a percolation configuration $\sigma$, we associate an
arrow to each edge of $\cal H$ belonging to the boundary of
a cluster in such a way that the hexagon to the right of the edge
with respect to the direction of the arrow is blue (plus).
The set of all boundaries then becomes a collection of nested,
oriented, simple loops.
A {\bf boundary path} (or {\bf b-path}) $\gamma$ is a sequence
$(e_0, \ldots, e_n)$ of distinct edges of $\cal H$ belonging
to the boundary of a cluster and such that $e_{i-1}$ and $e_i$
meet at a vertex of $\cal H$ for all $i= 1, \ldots, n$.
To each b-path, we can associate a direction according to the
direction of the edges in the path.

Given a b-path $\gamma$, we denote by $\Gamma_B(\gamma)$
(respectively, $\Gamma_Y(\gamma)$) the set of blue (resp.,
yellow) hexagons (i.e., sites of $\cal T$) adjacent to
$\gamma$; we also let
$\Gamma(\gamma) \equiv \Gamma_B(\gamma) \cup \Gamma_Y(\gamma)$.

\subsection{The Percolation Exploration Process and Path} \label{explo}

\begin{figure}[!ht]
\begin{center}
\includegraphics[width=6cm]{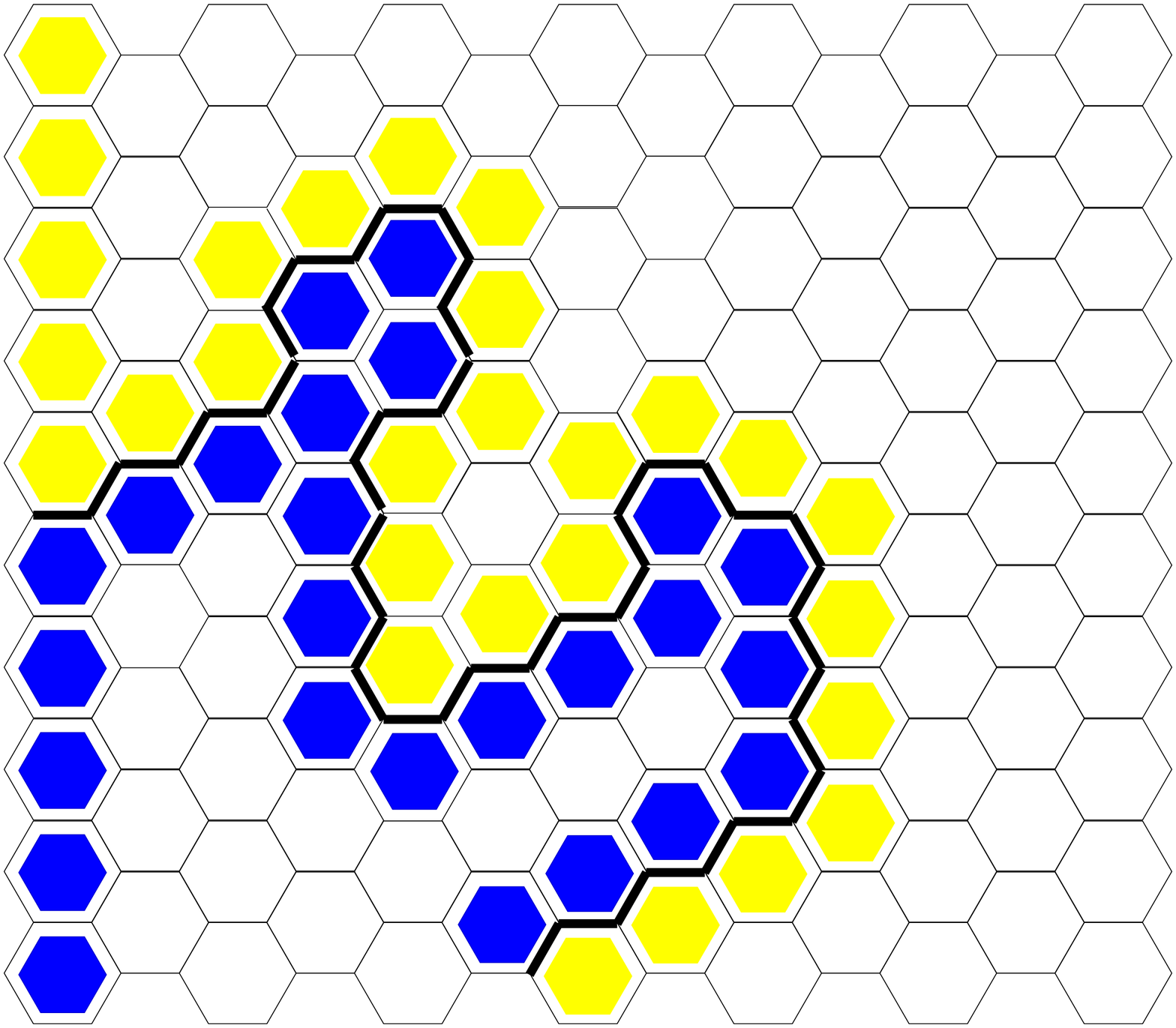}
\caption{Percolation exploration process in a portion of the
hexagonal lattice with blue/yellow boundary conditions on the
first column, corresponding to the boundary of the region where
the exploration is carried out.
The colored hexagons that do not belong to the first column
have been ``explored'' during the exploration process.
The heavy line between yellow (light) and blue (dark) hexagons
is the exploration path produced by the exploration process.}
%All the hexagons in the first line from below belong to the
%site boundary of the domain and are yellow (light shade in
%the figure) if they are to the left of the origin $0$ and
%blue (heavy shade in the figure) if they are to the right $0$.
%The other shaded hexagons have been ``explored'' during the
%process, which produces an exploration path (heavy line)
%along the edges of $\cal H$.}
\label{fig2-sec4}
\end{center}
\end{figure}

For a Jordan set $D \subset {\cal T}$ and two e-vertices $x,y$
in $\partial D$, imagine coloring blue all the hexagons in
$\Delta_{x,y} D$ and yellow all those in $\Delta_{y,x} D$.
Then, for any percolation configuration $\sigma_D$ inside
$D$, there is a unique b-path $\gamma$ from $x$ to $y$ which
separates the blue cluster adjacent to $\Delta_{x,y} D$ from
the yellow cluster adjacent to $\Delta_{y,x} D$.
We call $\gamma = \gamma^1_{D,x,y}(\sigma_D)$ a
{\bf percolation exploration path} (see Figure~\ref{fig2-sec4})
in $D$ from $x$ to $y$ with mesh size $\delta=1$.

%An exploration path $\gamma$ can be decomposed into
%{\bf left excursions} $\cal E$, i.e., maximal b-subpaths
%of $\gamma$ that do not use edges of the left boundary $\partial_{y,x} D$.
%Successive left excursions are separated by portions of
%$\gamma$ that contain only edges of the left boundary $\partial_{y,x} D$.
%Analogously, $\gamma$ can be decomposed into {\bf right excursions},
%i.e., maximal b-subpaths of $\gamma$ that do not use edges of the right
%boundary $\partial_{x,y} D$.
%Successive right excursions are separated by portions of
%$\gamma$ that contain only edges of the right boundary $\partial_{x,y} D$.

Notice that the exploration path $\gamma$
%$\gamma=\gamma_{D,x,y}(\sigma_D)$
only depends on
the percolation configuration $\sigma_D$ inside $D$
and the positions of the e-vertices $x$ and $y$;
in particular, it does not depend on the color of
the hexagons in $\Delta D$, since it is defined by
imposing fictitious boundary conditions on $\Delta D$.
To see this more clearly, we next show how to construct
the percolation exploration path dynamically, via
the {\bf percolation exploration process} defined
below.

Given a Jordan set $D \subset {\cal T}$ and two
e-vertices $x,y$ in $\partial D$, assign to $\partial_{x,y} D$
a counterclockwise orientation (i.e., from $x$ to $y$) and
to $\partial_{y,x} D$ a clockwise orientation.
%Call $e_x$ the edge incident on $x$ that does not belong to
%$\partial D$ and orient it in the direction of $x$;
%this is the ``starting edge'' of an exploration procedure
%that will produce an oriented path inside $D$ along
%the edges of $\cal H$, together with two \emph{nonsimple}
%monochromatic paths on $\cal T$.
%From $e_x$, the process moves along the edges of
%hexagons in $D$ according to the rules below.
The edge $e_x$ incident on $x$ that does not belong to $\partial D$
is oriented in the direction of $x$.
From there one starts an exploration procedure that produces
an oriented path inside $D$ along the edges of $\cal H$, together
with two \emph{nonsimple} monochromatic paths on $\cal T$, as
follows.
At each step there are two possible edges (left or right edge
with respect to the current direction of exploration) to choose
from, both belonging to the same hexagon $\xi$ contained in $D$
or $\Delta D$.
If $\xi$ belongs to $D$ and has not been previously ``explored,"
its color is determined by flipping a fair coin and then
the edge to the left (with respect to the direction in which
the exploration is moving) is chosen if $\xi$ is blue (plus),
or the edge to the right is chosen if $\xi$ is yellow (minus).
If $\xi$ belongs to $D$ and has been previously explored, the
color already assigned to it is used to choose an edge according to
the rule above.
If $\xi$ belongs to the right external boundary $\Delta_{x,y} D$,
the left edge is chosen.
If $\xi$ belongs to the left external boundary $\Delta_{y,x} D$,
the right edge is chosen.
The exploration stops when it reaches $y$.

Next, we introduce a class of domains of the plane which will appear
later in Theorems~\ref{strong-cardy} and~\ref{spatial-markov} and
various lemmas.
Let $D$ be a bounded simply connected domain whose boundary
$\partial D$ is a continuous curve.
Let $\phi:{\mathbb D} \to D$ be the (unique) conformal map from the
unit disc $\mathbb D$ to $D$ with
$\phi(0)=z_0 \in D$ and $\phi'(0)>0$;
note that by Theorem~\ref{cont-thm} of Appendix~\ref{rado}, $\phi$
has a continuous extension to $\overline{\mathbb D}$.
Let $a,c,d$ be three points of $\partial D$ (or more accurately,
three prime ends) in counterclockwise order -- i.e., such that
$a=\phi(a^*)$, $c=\phi(c^*)$ and $d=\phi(d^*)$, with $a^*$, $c^*$
and $d^*$ three distinct points of $\partial {\mathbb D}$ in
counterclockwise order.
We will call $D$ {\bf admissible} with respect to $(a,c,d)$ if
the counterclockwise arcs $J_1 \equiv \overline{da}$, $J_2 \equiv \overline{ac}$
and $J_3 \equiv \overline{cd}$ are simple curves, $J_3$ does not
touch the interior of either $J_1$ or $J_2$, and from each point
in $J_3$ there is a path to infinity that does not cross $\partial D$.
(Note that a Jordan $D$ is admissible for any counterclockwise
$a,c,d$ on $\partial D$.)
%For each $i=1,2,3$, let ${\cal J}_i: [0,1] \to J_i$, be a
%parametrization of $J_i$.

Notice that, according to our definition, the interiors of the
arcs $J_1$ and $J_2$ can touch.
If that happens, the double-points of the boundary (belonging to
both arcs) are counted twice and considered as two distinct points
(and are two different prime ends).
The significance of the notion of admissible is that certain domains
arising naturally in the proof of Theorem~\ref{spatial-markov} are
not Jordan but are admissible; this is because the
hulls $K_t$ generated by
%sets $K_t$ of
chordal $SLE_6$ paths
have cut-points~\cite{beffara1} -- see Figure~\ref{cut-point-fig}.

With $D,J_1,J_2,J_3$ as just discussed, let now $\{D^{\delta}\}$ be
a sequence of Jordan sets in $\delta {\cal T}$ (i.e., composed of the
hexagons of the scaled hexagonal lattices $\delta {\cal H}$).
If we can split $\partial D^{\delta}$ into three Jordan arcs,
$J^{\delta}_1$, $J^{\delta}_2$, $J^{\delta}_3$, such that
$\text{d}(J^{\delta}_i,J_i) \to 0$ for each $i=1,2,3$ as $\delta \to 0$,
we say that $\partial D^{\delta}$ converges to $\partial D$
%[FEDERICO: SHOULD WE KEEP THIS? in the uniform metric~(\ref{distance}) on continuous curves,]
as $\delta \to 0$ and write $\partial D^{\delta} \to \partial D$
or, equivalently, $D^{\delta} \to D$.

Let $a^{\delta}$ and $b^{\delta}$ be distinct e-vertices of $\partial D^{\delta}$
and let
%({\tt Federico (02-06-06): I want to avoid having
%$\gamma^{\delta}_{D,a,b} = \gamma^{\delta}_{D^{\delta},a^{\delta},b^{\delta}}$,
%which I don't like much. But I'm not sure what is worse. Here is a possible way out.})
%$\gamma = \delta \gamma^1_{\frac{1}{\delta}D^{\delta},\frac{1}{\delta}a^{\delta},\frac{1}{\delta}b^{\delta}}$
%$\gamma = \gamma^{\delta}_{D^{\delta},a^{\delta},b^{\delta}}$
$\gamma$ be the exploration path in $D^{\delta}$ from $a^{\delta}$ to $b^{\delta}$.
If, as $\delta \to 0$, $\partial D^{\delta} \to \partial D$, $a^{\delta} \to a$
and $b^{\delta} \to b$, where $D$ is a domain admissible with respect to $(a,c,d)$
and $b \in J_3=\overline{cd}$, we say that $(D^{\delta},a^{\delta},b^{\delta})$ is
a {\bf $\delta$-approximation} to $(D,a,b)$, write
$(D^{\delta},a^{\delta},b^{\delta}) \to (D,a,b)$, and denote the exploration
path $\gamma$ by $\gamma^{\delta}_{D,a,b}$.
Note that $\gamma^{\delta}_{D,a,b}$ depends not only on $(D,a,b)$, but
also on the $\delta$-approximation to $(D,a,b)$.
For simplicity of notation, we do not make explicit this dependence.

For a fixed $\delta>0$, the probability measure $\mathbb P$ on percolation
configurations induces a probability measure $\mu^{\delta}_{D,a,b}$
on exploration paths $\gamma^{\delta}_{D,a,b}$.
In the continuum scaling limit, $\delta \to 0$, one is interested in the weak
convergence with respect to the uniform metric~(\ref{distance})
of $\mu^{\delta}_{D,a,b}$ to a probability measure $\mu_{D,a,b}$
supported on continuous curves.

Before concluding this section, we give some more definitions.
Consider the exploration path $\gamma = \gamma^{\delta}_{D,a,b}$
and the set $\Gamma(\gamma) = \Gamma_Y(\gamma) \cup \Gamma_B(\gamma)$.
The set $D \setminus \Gamma(\gamma)$ is the
union of its connected components (in the lattice sense),
which are simply connected.
For $\delta$ small and $a,b \in \partial D$ not too close to each other,
with high probability the exploration process inside $D^{\delta}$ will make
large excursions into $D^{\delta}$, so that $D^{\delta} \setminus \Gamma(\gamma)$
will have more than one component.
Given a point $z \in {\mathbb C}$ contained in $D^{\delta} \setminus \Gamma(\gamma)$,
we will denote by
$D^{\delta}_{a,b}(z)$ the domain corresponding
to the unique element of $D^{\delta} \setminus \Gamma(\gamma)$ that contains $z$
(notice that for a deterministic $z \in D$,
$D^{\delta}_{a,b}(z)$
is well defined with high probability for $\delta$ small, i.e., when $z \in D^{\delta}$
and $z \notin \Gamma(\gamma)$).

There are four types of domains which may be usefully
thought of in terms of their external site boundaries:
(1) those components  whose  site boundary contains both
sites in $\Gamma_Y(\gamma)$ and $\Delta_{b^{\delta},a^{\delta}} D^{\delta}$,
(2) the analogous components  with $\Delta_{b^{\delta},a^{\delta}} D^{\delta}$
replaced by $\Delta_{a^{\delta},b^{\delta}} D^{\delta}$ and $\Gamma_Y(\gamma)$
by $\Gamma_B(\gamma)$, (3) those components whose site boundary only contains
sites in $\Gamma_Y(\gamma)$, and finally (4) the analogous components  with
$\Gamma_Y(\gamma)$  replaced by $\Gamma_B(\gamma)$.
These different types will appear in the proof of Lemma~\ref{boundaries}.

\section{Cardy's Formula and a Characterization of $SLE_6$} \label{convergence-to-sle}

The existence of subsequential limits for the percolation exploration
path, which follows from the work of Aizenman and Burchard~\cite{ab},
means that the proof of convergence to $SLE_6$ can be divided into two
parts: first we will give a characterization of chordal $SLE_6$ in terms
of two properties that determine it uniquely; then we will show that any
subsequential scaling limit of the percolation exploration path satisfies
these two properties.

The characterization part will follow from known properties of hulls
and of $SLE_6$ (see~\cite{lsw7} and~\cite{werner}).
The second part will follow from an extension of Smirnov's result about the
convergence of crossing probabilities to Cardy's formula~\cite{cardy}
(see Theorem~\ref{strong-cardy} below) for sequences of Jordan domains $D_k$,
with the domain $D_k$ changing together with the mesh $\delta_k$ of the lattice,
combined with the proof of a certain spatial Markov property for subsequential
limits of percolation exploration hulls (Theorem~\ref{spatial-markov}).
We note that although Theorem~\ref{strong-cardy} represents only a slight
extension to Smirnov's result on convergence of crossing probabilities, this
extension and its proof play a major role in the technically important
Lemmas~\ref{double-crossing}, \ref{equal} and~\ref{mushroom}, which control
the ``close encounters" of exploration paths to domain boundaries.
The proof of Theorem~\ref{strong-cardy} is modelled after a simpler
geometric argument involving only rectangles used in~\cite{cns1}.

Let $D$ be a bounded simply connected domain containing the origin
whose boundary $\partial D$ is a continuous curve.
Let $\phi:{\mathbb D} \to D$ be the (unique) conformal map from the
unit disc to $D$ with $\phi(0)=0$ and $\phi'(0)>0$; note that by
Theorem~\ref{cont-thm} of Appendix~\ref{rado}, $\phi$ has a continuous
extension to $\overline{\mathbb D}$.
Let $z_1,z_2,z_3,z_4$ be four points of $\partial D$ in counterclockwise
order -- i.e., such that $z_j=\phi(w_j), \,\,\, j=1,2,3,4$, with
$w_1,\ldots,w_4$ in counterclockwise order.
Also, let $\eta = \frac{(w_1-w_2)(w_3-w_4)}{(w_1-w_3)(w_2-w_4)}$.
Cardy's formula~\cite{cardy} for the probability $\Phi_{D}(z_1,z_2;z_3,z_4)$
of a ``crossing" in $D$ from the counterclockwise arc $\overline{z_1 z_2}$
to the counterclockwise arc $\overline{z_3 z_4}$ is
\begin{equation} \label{cardy-formula}
\Phi_{D}(z_1,z_2;z_3,z_4) =
\frac{\Gamma(2/3)}{\Gamma(4/3) \Gamma(1/3)} \eta^{1/3} {}_2F_1(1/3,2/3;4/3;\eta),
\end{equation}
where ${}_2F_1$ is a hypergeometric function.

For a given mesh $\delta>0$, the probability of a blue crossing inside $D$
from the counterclockwise arc $\overline{z_1 z_2}$ to the counterclockwise
arc $\overline{z_3 z_4}$ is the probability of the existence of a blue
$\cal T$-path $(\xi_0,\ldots,\xi_n)$ such that $\xi_0$ intersects the
counterclockwise arc $\overline{z_1 z_2}$, $\xi_n$ intersects the
counterclockwise arc $\overline{z_3 z_4}$ and $\xi_1,\ldots,\xi_{n-1}$
are all contained in $D$.
Smirnov proved the following major theorem, concerning the conjectured
behavior~\cite{cardy} of crossing probabilities in the scaling limit
(see also~\cite{beffara2}).

%For a given mesh $\delta>0$, we denote by ${\cal G}(D,\delta)$ the largest
%Jordan set of hexagons of the scaled hexagonal lattice $\delta {\cal H}$ that
%is contained in $D$, and by $z_1^{\delta},z_2^{\delta},z_3^{\delta},z_4^{\delta}$
%the e-vertices of $\partial {\cal G}(D,\delta)$ closest to $z_1,z_2,z_3,z_4$
%respectively (if there are two such vertices closest to $v_1$, we choose, say,
%the first one encountered going clockwise along $\partial {\cal G}(D,\delta)$,
%and analogously for $v_2,v_3,v_4$).
%The probability of a crossing inside ${\cal G}(D,\delta)$ from the
%counterclockwise arc $\overline{z_1^{\delta} z_2^{\delta}}$ to the
%counterclockwise arc $\overline{z_3^{\delta} z_4^{\delta}}$ is the probability
%of the existence of a blue $\cal T$-path contained in ${\cal G}(D,\delta)$ that
%starts at a hexagon touching $\overline{z_1^{\delta} z_2^{\delta}}$ and ends at
%a hexagon touching $\overline{z_3^{\delta} z_4^{\delta}}$.
%Smirnov proved the following major theorem (see also~\cite{beffara}), concerning
%the conjectured behavior~\cite{cardy} of crossing probabilities in the scaling limit.

\begin{theorem} \emph{(Smirnov~\cite{smirnov})} \label{cardy-smirnov}
Let $D$ be a Jordan domain whose boundary $\partial D$ is a finite
union of smooth (e.g., $C^2$) curves.
%piecewise smooth.
As $\delta \to 0$, the limit of the probability of a blue crossing inside
$D$ from the counterclockwise arc $\overline{z_1^{\delta} z_2^{\delta}}$
to the counterclockwise arc $\overline{z_3^{\delta} z_4^{\delta}}$ exists,
is a conformal invariant of $(D,z_1,z_2,z_3,z_4)$ and is given by Cardy's
formula~(\ref{cardy-formula}).
\end{theorem}

%\begin{theorem} \emph{(Smirnov~\cite{smirnov})} \label{cardy-smirnov}
%For fixed $D'$, $z_1$, $z_2$, $z_3$, $z_4$, the limit as $\delta \to 0$
%of the above crossing probability exists, is a conformal invariant of
%$D'$, $z_1$, $z_2$, $z_3$, $z_4$, and is
%given by Cardy's formula~(\ref{cardy-formula}).
%\end{theorem}
%
%\begin{theorem} \emph{(Smirnov~\cite{smirnov})} \label{cardy-smirnov}
%For fixed $(G,z_1,z_2,z_3,z_4)$, there is {\bf some}
%$\delta$-approximation $(D'^{\delta},z_1^{\delta},z_2^{\delta},
%z_3^{\delta},z_4^{\delta})$ converging to $(D',z_1,z_2,z_3,z_4)$
%such that the limit as $\delta \to 0$
%of the above crossing probability exists, is a conformal invariant of
%$(D',z_1,z_2,z_3,z_4)$ and is given by Cardy's formula~(\ref{cardy-formula}).
%\end{theorem}

\begin{remark}
We have stated Smirnov's result in the form that will be used in this paper,
but note that Smirnov does not restrict attention to Jordan domains with
piecewise smooth boundary but rather allows for more general bounded
domains (see~\cite{smirnov,smirnov-long}).
We also remark that Theorem~\ref{strong-cardy} below extends Theorem~\ref{cardy-smirnov}
to a larger class of domains, including in particular all Jordan domains.
\end{remark}

%\begin{remark} \label{remark-hausdorff}
%NEW REMARK
%We note that, although ${\cal G}(D,\delta)$ is not a priori a $\delta$-ap\-proxi\-ma\-tion
%of $D$ in the sense of Definition~\ref{approx}, it is easy to see that $\partial{\cal G}(D,\delta)$
%converges to $\partial D$ at least in the Hausdorff metric as long as $G$ is a Jordan domain.
%\end{remark}

Let us now specify the objects that we are interested in.
Suppose $D$ is a simply connected domain whose boundary $\partial D$
is a continuous curve, and $a,b$ are two distinct points in $\partial D$
(or more accurately, two distinct prime ends), and let $\tilde\mu_{D,a,b}$
be a probability measure on continuous curves
$\tilde\gamma = \tilde\gamma_{D,a,b}:[0,\infty) \to \overline D$
with $\tilde\gamma(0)=a$,
$\tilde\gamma(\infty) \equiv \lim_{t \to \infty} \tilde\gamma(t) = b$,
and $\tilde\gamma(t) \neq b$ for $t$ finite
(we remark that the use of $[0,\infty)$ instead of $[0,1]$ for the time
parametrization is purely for convenience).
Let $D_t \equiv D \setminus \tilde K_t$ denote the (unique) connected
component of $D \setminus \tilde\gamma[0,t]$ whose closure contains $b$,
where $\tilde K_t$, the {\bf filling} of $\tilde\gamma[0,t]$, is a closed
connected subset of $\overline D$.
$\tilde K_t$ is called a {\bf hull} if it satisfies the condition
\begin{equation} \label{hulls}
\overline{\tilde K_t \cap D} = \tilde K_t.
\end{equation}
We will consider curves $\tilde\gamma$ such that
%(i) $\tilde K_t$ is a hull for each $t$, although we normally only consider
%$\tilde K_T$ at certain stopping times $T$.
%(An example of such a curve that we are particularly interested in is the
%trace of chordal $SLE_6$.)
%We also want that
(i) $\tilde \gamma$ is the limit (in distribution using the metric~(\ref{distance}))
of random simple curves and such that, for $0<t_1<t_2$ with
$\tilde \gamma (t_1) \neq \tilde \gamma (t_2)$,
(ii) $\tilde \gamma (t_2) \notin \tilde K_{t_1} \setminus \partial \tilde K_{t_1}$
and (iii) $\exists t' \in (t_1,t_2)$ with $\tilde \gamma(t') \in D_{t_1} \equiv
D \setminus \tilde K_{t_1}$.
We note that an example of a curve satisfying these properties is the trace
of chordal $SLE_6$.

%Let $\tilde K_t$ denote the hull of $\tilde\gamma$ at time $t$ and
%let $D_t \equiv D \setminus \tilde K_t$; the {\bf hull} $\tilde K_t$ is
%defined to be the closed set such that $D_t$ is the (unique) connected
%component of $D \setminus \tilde\gamma[0,t]$ whose closure contains $b$.

\begin{figure}[!ht]
\begin{center}
\includegraphics[width=7cm]{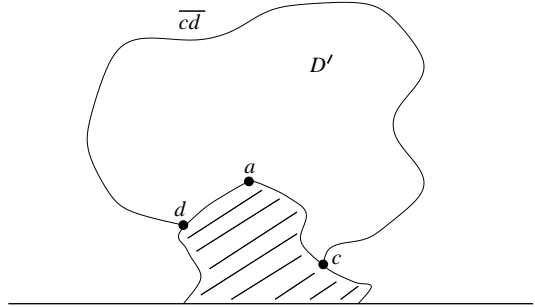}
\caption{$D$ is the upper half-plane $\mathbb H$ with the shaded portion removed,
$b=\infty$, $C'$ is an unbounded subdomain, and $D' = D \setminus C'$ is indicated
in the figure. The counterclockwise arc $\overline{cd}$ indicated in the figure
belongs to $\partial D'$.}
\label{fig1-sec7}
\end{center}
\end{figure}

Let $C' \subset D$ be a closed subset of $\overline D$ such that $a \notin C'$,
$b \in C'$, and $D' = D \setminus C'$ is a bounded simply connected domain
whose boundary
is a continuous curve containing the counterclockwise arc $\overline {cd}$
that does not belong to $\partial D$ (except for its endpoints $c$ and $d$ --
see Figure~\ref{fig1-sec7}).
Let $T'=\inf \{ t:\tilde K_t \cap C' \neq \emptyset \}$ be the first time
that $\tilde\gamma(t)$ hits $C'$.
We say that the hitting distribution of $\tilde\gamma(t)$ at the stopping
time $T'$ is determined by Cardy's formula if, for any $C'$ and any
counterclockwise arc $\overline{xy}$ of $\overline{cd}$, the probability
that $\tilde\gamma$ hits $C'$ at time $T'$ on $\overline{xy}$ is given by
\begin{equation}
{\mathbb P}^*(\tilde\gamma(T') \in \overline{xy}) = \Phi_{D'}(a,c;x,d) - \Phi_{D'}(a,c;y,d).
\end{equation}

\begin{figure}[!ht]
\begin{center}
\includegraphics[width=7cm]{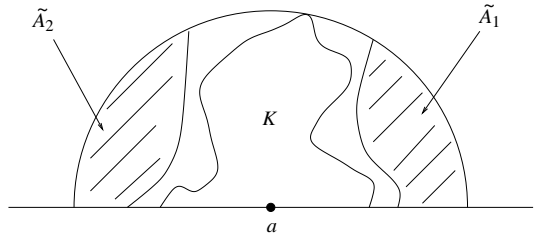}
\caption{Example of a hull $K$ and a set $\tilde A_1 \cup \tilde A_2$ in $\cal A$.
Here, $D = {\mathbb H}$ and $D'$ is the semi-disc centered at $a$.}
\label{fig2-sec7}
\end{center}
\end{figure}

Assume that the filling $\tilde K_{T'}$ of $\tilde\gamma[0,T']$ is a hull;
we denote by $\tilde\nu_{D',a,c,d}$ the distribution of $\tilde K_{T'}$.
To explain what we mean by the distribution of a hull, consider the set
$\tilde{\cal A}$ of closed subsets $\tilde A$ of $\overline{D'}$ that do not
contain $a$ and such that $\partial \tilde A \setminus \partial D'$ is a simple
(continuous) curve contained in $D'$ except for its endpoints, one of which is
on $\partial D' \cap D$ and the other is on $\partial D$ (see Figure~\ref{fig2-sec7}).
Let $\cal A$ be the set of closed subsets of $\overline{D'}$ of the form
$\tilde A_1 \cup \tilde A_2$, where $\tilde A_1, \tilde A_2 \in \tilde{\cal A}$
and $\tilde A_1 \cap \tilde A_2 = \emptyset$.
%Let $\cal A$ be the set of closed subsets of $\overline{D'}$ of the form
%$\tilde A_1 \cup \tilde A_2$, where $\tilde A_1, \tilde A_2 \in \tilde{\cal A}$
%and $\tilde A_1 \cap \tilde A_2 = \emptyset$.

For a given $C'$ and corresponding $T'$, let $\cal K$ be the set whose
elements are possible hulls at time $T'$; we claim that the events
$\{ K \in {\cal K}: K \cap A = \emptyset \}$, for $A \in {\cal A}$,
form a $\pi$-system $\Pi$ (i.e., they are closed under finite intersections;
we also include the empty set in $\Pi$), and we consider the $\sigma$-algebra
$\Sigma=\sigma(\Pi)$ generated by these events.
To see that $\Pi$ is closed under pairwise intersections, notice that, if
$A_1, A_2 \in {\cal A}$, then $\{ K \in {\cal K}: K \cap A_1 = \emptyset \}
\cap \{ K \in {\cal K}: K \cap A_2 = \emptyset \} =
\{ K \in {\cal K}: K \cap \{ A_1 \cup A_2 \} = \emptyset \}$
and $A_1 \cup A_2 \in {\cal A}$ (or else
$\{ K \in {\cal K} : K \cap \{ A_1 \cup A_2 \} = \emptyset \}$ is empty).
We are interested in probability spaces of the form $({\cal K},\Sigma,{\mathbb P}^*)$.

%Consider now a closed subset $C''$ of $D \setminus \tilde K_{T'}$ such that
%$\tilde\gamma(T') \notin C''$, $b \in C''$, and $D''  = D_{T'} \setminus C''$
%is a bounded Jordan domain whose boundary contains the counterclockwise arc
%$\overline {ef}$ that does not belong to $\partial D_{T'}$ (except for its endpoints
%$e$ and $f$).
%Let $T''$ be the first time $\tilde\gamma(t)$ exits $D''$; we say that $\tilde K_t$
%satisfies the {\bf spatial Markov property} if, conditioned on $\tilde K_{T'}$, the
%distribution of $\tilde K_{T''} \setminus \tilde K_{T'}$ is the same as
%$\tilde\nu_{D'',\tilde\gamma(T'),e,f}$.

It is easy to see that if the hitting distribution of $\tilde\gamma(t)$ is
determined by Cardy's formula, then the probabilities of events in $\Pi$
are also determined by Cardy's formula in the following way.
Let $A \in {\cal A}$ be the union of $\tilde A_1, \tilde A_2 \in \tilde{\cal A}$,
with $\partial \tilde A_1 \setminus \partial D'$ given by a curve from
$u_1 \in \partial D' \cap D$ to $v_1 \in \partial D$ and
$\partial \tilde A_2 \setminus \partial D'$ given by a curve from
$u_2 \in \partial D' \cap D$ to $v_2 \in \partial D$; then, assuming that
$a$, $v_1$, $u_1$, $u_2$, $v_2$ are ordered counterclockwise around $\partial D'$,
\begin{equation}
{\mathbb P}^*(\tilde K_{T'} \cap A = \emptyset) =
\Phi_{D' \setminus A}(a,v_1;u_1,v_2,) - \Phi_{D' \setminus A}(a,v_1;u_2,v_2).
\end{equation}
Since $\Pi$ is a $\pi$-system, the probabilities of the events in $\Pi$ determine
uniquely the distribution of the hull in the sense described above.
Therefore, if we let $\gamma_{D,a,b}$ denote the trace of chordal $SLE_6$
inside $D$ from $a$ to $b$, $K_t$ its hull up to time $t$, and
$\tau = \inf \{ t : K_t \cap C' \neq \emptyset \}$ the first time that
$\gamma_{D,a,b}$ hits $C'$, we have the following simple but useful lemma.

\begin{lemma} \label{hull}
With the notation introduced above,
if $\tilde K_{T'}$ is a hull and the hitting
distribution of $\tilde\gamma_{D,a,b}$
at the stopping time $T'$ is determined by Cardy's formula,
then $\tilde K_{T'}$ is distributed like the hull $K_{\tau}$ of chordal $SLE_6$.
\end{lemma}

\noindent {\bf Proof.} It is enough to note that the hitting distribution for
%the hull $K_t$ of
chordal $SLE_6$ is determined by Cardy's formula~\cite{lsw1}. \fbox{} \\

Now let $\tilde f_0$ be a conformal map from the upper half-plane $\mathbb H$
to $D$ such that ${\tilde f_0}^{-1}(a)=0$ and ${\tilde f_0}^{-1}(b)=\infty$.
(Since $\partial D$ is a continuous curve, the map ${\tilde f_0}^{-1}$ has a
continuous extension from $D$ to $D \cup \partial D$ -- see Theorem~\ref{cont-thm}
of Appendix~\ref{rado} -- and, by a slight abuse of notation, we do not
distinguish between ${\tilde f_0}^{-1}$ and its extension; the same applies
to $\tilde f_0$.)
These two conditions determine $\tilde f_0$ only up to a scaling factor.
For $\varepsilon>0$ fixed, let
$C(u,\varepsilon) = \{ z:|u-z|<\varepsilon \} \cap {\mathbb H}$
denote the semi-ball of radius $\varepsilon$ centered at $u$ on
the real line and let $\tilde T_1=\tilde T_1(\varepsilon)$ denote
the first time $\tilde\gamma(t)$ hits $D \setminus \tilde G_1$,
where $\tilde G_1 \equiv \tilde f_0(C(0,\varepsilon))$.
Define recursively $\tilde T_{j+1}$ as the first time
$\tilde\gamma[\tilde T_j,\infty)$ hits $\tilde D_{\tilde T_j} \setminus \tilde G_{j+1}$,
where $\tilde D_{\tilde T_j} \equiv D \setminus \tilde K_{\tilde T_j}$,
$\tilde G_{j+1} \equiv \tilde f_{\tilde T_j}(C(0,\varepsilon))$, and
$\tilde f_{\tilde T_j}$ is a conformal map from $\mathbb H$ to $\tilde D_{\tilde T_j}$
whose inverse maps $\tilde\gamma(\tilde T_j)$ to $0$ and $b$ to $\infty$.
We also define $\tilde\tau_{j+1} \equiv \tilde T_{j+1} - \tilde T_j$,
so that $\tilde T_j = \tilde\tau_1+\ldots+\tilde\tau_j$.
%We note that, like $\tilde f_0$, the conformal maps $\tilde f_{\tilde T_j}$
%are only defined up to a scaling factor.
We choose $\tilde f_{\tilde T_j}$ so that its inverse is the composition
of the restriction of $\tilde{f_0}^{-1}$ to $\tilde D_{\tilde T_j}$ with
$\tilde \varphi_{\tilde T_j}$, where $\tilde \varphi_{\tilde T_j}$ is the
unique conformal transformation from
${\mathbb H} \setminus \tilde{f_0}^{-1}(\tilde K_{\tilde T_j})$ to $\mathbb H$
that maps $\infty$ to $\infty$ and $\tilde{f_0}^{-1} (\tilde\gamma(\tilde T_j))$
to the origin of the real axis, and has derivative at $\infty$ equal to $1$.
%with expansion at infinity:
%\begin{equation}
%\tilde g_t (z) = z + \frac{\tilde a(t)}{z} + o(\frac{1}{z}).
%\end{equation}
%For concreteness, we may specify the time parametrization of $\tilde\gamma$
%so that $a(t)=2t$.

Notice that $\tilde G_{j+1}$ is a bounded simply connected domain chosen
so that the conformal transformation which maps $\tilde D_{\tilde T_j}$
to $\mathbb H$ maps $\tilde G_{j+1}$ to the semi-ball $C(0,\varepsilon)$
centered at the origin on the real line.
%(We also note that, because of properties of conformal transformations, the image
%of the semicircle $\{ z:|u-z|=\varepsilon \} \cap {\mathbb H}$ under $f_{\tilde T_j}$
%is a continuous curve whose second derivative, away from its endpoints, is uniformly
%bounded.)
With these definitions, consider the (discrete-time) stochastic process
$\tilde X_j \equiv (\tilde K_{\tilde T_j}, \tilde\gamma(\tilde T_j))$ for $j=1,2,\ldots$;
we say that $\tilde K_t$ satisfies the {\bf spatial Markov property} if each
$\tilde K_{\tilde T_j}$ is a hull and $\tilde X_j$ for $j=1,2,\ldots$ is a
Markov chain (for any choice of the map $\tilde f_0$).
%scaling factors for $\tilde f_0, \tilde f_{\tilde T_1}, \tilde f_{\tilde T_2}, \ldots$).
Notice that the hull of chordal $SLE_6$ satisfies the spatial Markov property,
due to the conformal invariance and Markovian properties~\cite{schramm} of $SLE_6$.

\begin{remark} \label{remark-differ}
The next theorem, our main characterization result for $SLE_6$, uses the
choice of stopping times we have just discussed rather than that proposed
by Smirnov~\cite{smirnov,smirnov-long}.
A technical reason for this revision of Smirnov's strategy is discussed
in Remark~\ref{remark-difference} below.
But a conceptually more important reason is that it naturally gives rise
(in the scaling limit) to a certain random walk on the real line (the
sequence of points to which the tips of the hulls are conformally mapped
at the stopping times) whose increments are i.i.d. random variables.
As the stopping time parameter $\varepsilon \to 0$, this random walk
converges to the driving Brownian motion of the $SLE_6$.
\end{remark}

\begin{theorem} \label{characterization}
If the filling process $\tilde K_t$ of $\tilde\gamma_{D,a,b}$ satisfies the spatial
Markov property and its hitting distribution is determined by Cardy's formula, then
$\tilde\gamma_{D,a,b}$ is distributed like the trace $\gamma_{D,a,b}$ of chordal
$SLE_6$ inside $D$ started at $a$ and aimed at $b$.
\end{theorem}

\noindent {\bf Proof.} Since the trace $\gamma_{D,a,b}$ of chordal $SLE_6$ in a
Jordan domain $D$ is defined (up to a linear time change) as $f(\gamma)$,
where $\gamma = \gamma_{{\mathbb H},0,\infty}$ is the trace of chordal $SLE_6$
in the upper half-plane started at $0$ and $f$ is any conformal map from the upper
half-plane $\mathbb H$ to $D$ such that $f^{-1}(a)=0$ and $f^{-1}(b)=\infty$, it is
enough to show that $\hat\gamma = f^{-1}(\tilde\gamma_{D,a,b})$ is distributed like
the trace of chordal $SLE_6$ in the upper half-plane.
Let $\hat K_t$ denote the filling of $\hat\gamma(t)$ at time $t$ and let $\hat g_t(z)$
be the unique conformal transformation that maps ${\mathbb H} \setminus \hat K_t$
onto $\mathbb H$ with the following expansion at infinity:
\begin{equation} \label{conf-map1}
\hat g_t(z) = z + \frac{\hat a(t)}{z} + o(\frac{1}{z}).
\end{equation}
We choose to parametrize $\hat\gamma(t)$ so that $t=\hat a(t)/2$
(this is often called parametrization by capacity, $\hat a(t)$ being
the half-plane capacity of the filling up to time $t$).

We want to compare $\hat\gamma(t)$ with the trace $\gamma(t)$ of chordal $SLE_6$
in the upper half-plane parameterized in the same way (i.e., with $a(t) = 2 \, t$),
so that, if $K_t$ is the filling of $\gamma$ at time $t$,
${\mathbb H} \setminus K_t$ is mapped onto $\mathbb H$ by a conformal
$g_t$ with the following expansion at infinity:
\begin{equation} \label{conf-map2}
g_t(z) = z + \frac{2t}{z} + o(\frac{1}{z}).
\end{equation}

Our strategy, following~\cite{smirnov,smirnov-long} but with modifications
(see Remark~\ref{remark-difference}), will be to construct suitable polygonal
approximations $\hat\gamma_{\varepsilon}$ and $\gamma_{\varepsilon}$ of $\hat\gamma$
and $\gamma$ which converge, as $\varepsilon \to 0$, to the original curves
(in the uniform metric on continuous curves~(\ref{distance})), and show that
$\hat\gamma_{\varepsilon}$ and $\gamma_{\varepsilon}$ have the same distribution.
This implies the equidistribution of $\hat\gamma$ and $\gamma$.

Let us describe first the construction for $\gamma_{\varepsilon}(t)$;
we do the same for $\hat\gamma_{\varepsilon}(t)$.
The important features in the construction of the polygonal approximations
are the spatial Markov property of the fillings and Cardy's formula, which
are valid for both $\gamma$ and $\hat\gamma$.

For $\varepsilon>0$ fixed, as above let
$C(u,\varepsilon) = \{ z:|u-z|<\varepsilon \} \cap {\mathbb H}$
denote the semi-ball of radius $\varepsilon$ centered at $u$ on the real line.
Let $T_1=T_1(\varepsilon)$ denote the first time $\gamma(t)$ hits
${\mathbb H} \setminus G_1$, where $G_1 \equiv C(0,\varepsilon)$, and
define recursively $T_{j+1}$ as the first time $\gamma[T_j,\infty)$
hits ${\mathbb H}_{T_j} \setminus G_{j+1}$, where
${\mathbb H}_{T_j} = {\mathbb H} \setminus K_{T_j}$ and
$G_{j+1} \equiv g^{-1}_{T_j}(C(g_{T_j}(\gamma(T_j)),\varepsilon))$.
Notice that $G_{j+1}$ is a bounded simply connected domain chosen so that the
conformal transformation which maps ${\mathbb H}_{T_j}$ to $\mathbb H$ maps
$G_{j+1}$ to the semi-ball $C(g_{T_j}(\gamma(T_j)),\varepsilon)$ centered at
the point of the real line where the ``tip" $\gamma(T_j)$ of the hull $K_{T_j}$
is mapped.
The spatial Markov property and the conformal invariance of the hull of $SLE_6$
imply that if we write $T_j = \tau_1+\ldots+\tau_j$, with
$\tau_{j+1} \equiv T_{j+1} - T_j$, the $\tau_j$'s are i.i.d. random variables,
and also that the distribution of $K_{T_{j+1}}$ is the same as that
of $K_{T_j} \cup g^{-1}_{T_j}(K'_{T_1} + g_{T_j}(\gamma(T_j)))$, where $K'_{T_1}$
is a hull equidistributed with $K_{T_1}$, but also is independent of $K_{T_1}$,
and ``$+ g_{T_j}(\gamma(T_j))$" indicates that it is translated by
$g_{T_j}(\gamma(T_j))$ along the real line.
The polygonal approximation $\gamma_{\varepsilon}$ is obtained by joining, for
all $j$, $\gamma(T_j)$ to $\gamma(T_{j+1})$ with a straight segment, where the
speed $\gamma_{\varepsilon}'(t)$ is constant.

Now let $\hat T_1 = \hat T_1(\varepsilon)$ denote the first time
$\hat\gamma(t)$ hits ${\mathbb H} \setminus \hat G_1$, where
$\hat G_1 \equiv C(0,\varepsilon)$, and define recursively
$\hat T_{j+1}$ as the first time $\hat\gamma[\hat T_j,\infty)$ hits
$\hat{\mathbb H}_{\hat T_j} \setminus \hat G_{j+1}$, where
$\hat{\mathbb H}_{\hat T_j} \equiv {\mathbb H} \setminus \hat K_{\hat T_j}$ and
$\hat G_{j+1} \equiv \hat g^{-1}_{T_j}(C(\hat g_{\hat T_j}(\hat\gamma(\hat T_j)),\varepsilon))$.
We also define $\hat\tau_{j+1} \equiv \hat T_{j+1} - \hat T_j$, so that
$\hat T_j = \hat\tau_1+\ldots+\hat\tau_j$.
Once again, $\hat G_{j+1}$ is a bounded simply connected domain chosen so that the
conformal transformation which maps $\hat{\mathbb H}_{\hat T_j}$ to $\mathbb H$ maps
$\hat G_{j+1}$ to the semi-ball $C(\hat g_{\hat T_j}(\hat\gamma(\hat T_j)),\varepsilon)$
centered at the point on the real line where the ``tip" $\hat\gamma(\hat T_j)$ of the
hull $\hat K_{\hat T_j}$ is mapped.
The polygonal approximation $\hat\gamma_{\varepsilon}$ is obtained by joining,
for all $j$, $\hat\gamma(\hat T_j)$ to $\hat\gamma(\hat T_{j+1})$ with a straight
segment, where the speed $\hat\gamma_{\varepsilon}'(t)$ is constant.

Consider the sequence of times
$\tilde T_j$ defined in the natural way so that
$\tilde\gamma(\tilde T_j) = f(\hat\gamma(\hat T_j))$
and the (discrete-time) stochastic processes
$\hat X_j \equiv (\hat K_{\hat T_j}, \hat\gamma(\hat T_j))$ and
$\tilde X_j \equiv (\tilde K_{\tilde T_j}, \tilde\gamma(\tilde T_j))$
related by $\hat X_j = f^{-1}(\tilde X_j)$.
If for $x \in {\mathbb R}$ we let $\theta[x]$ denote the translation
that maps $x$ to $0$ and define the family of conformal maps
$(\tilde f_{\tilde T_j})^{-1} = \theta[g_{\hat T_j}(\hat\gamma(\hat T_j))] \circ g_{\hat T_j} \circ f^{-1}$
from $D \setminus \tilde K_{\tilde T_j}$ to $\mathbb H$, then
$(\tilde f_{\tilde T_j})^{-1}$ sends $\tilde\gamma(\tilde T_j)$ to $0$ and $b$ to
$\infty$, and $(\tilde T_{j+1})$ is the first time $\tilde\gamma[\tilde T_j,\infty)$
hits $\tilde{\mathbb H}_{\hat T_j} \setminus \tilde G_{j+1}$, where
$\tilde{\mathbb H}_{\tilde T_j} = {\mathbb H} \setminus \tilde K_{\hat T_j}$
and $\tilde G_{j+1} = \tilde f_{\tilde T_j}(C(0,\varepsilon))$.
Therefore, $\{ \tilde T_j \}$ is a sequence of stopping times like
those used in the definition of the spatial Markov property and, thanks
to the relation $\hat X_j = f^{-1}(\tilde X_j)$, the fact that $\tilde K_t$
satisfies the spatial Markov property implies that $\hat X_j$ is a Markov chain.
We also note that the fact that the hitting distribution of $\tilde\gamma(t)$
is determined by Cardy's formula implies the same for the hitting distribution
of $\hat\gamma(t)$, thanks to the conformal invariance of Cardy's formula.
We next use these properties %of the hull $\hat K_t$
to show that $\hat\gamma_{\varepsilon}$ is distributed like $\gamma_{\varepsilon}$.

To do so, we first note that $g_{T_j}$ and $\hat g_{\hat T_j}$ are random
and their distributions are functionals of those of the hulls $K_{T_j}$ and
$\hat K_{\hat T_j}$, since there is a one-to-one correspondence between hulls and
conformal maps (with the normalization we have chosen in~(\ref{conf-map1})--(\ref{conf-map2})).
Therefore, since $\hat K_{\hat T_1}$ is distributed like $K_{T_1}$ (see Lemma~\ref{hull}),
$g_{T_1}$ and $\hat g_{\hat T_1}$ have the same distribution, which also implies that
$\hat T_1$ is distributed like $T_1$ because, due to the parametrization by capacity of
$\gamma$ and $\hat\gamma$, $2 T_1$ is exactly the coefficient of the term $1/z$ in the
expansion at infinity of $g_{T_1}$, and $2 \hat T_1$ is exactly the coefficient of the
term $1/z$ in the expansion at infinity of $\hat g_{\hat T_1}$.
Moreover, it is also clear that $\hat\gamma(\hat T_1)$ is distributed like $\gamma(T_1)$,
because their distributions are both determined by Cardy's formula, and so
$\hat g_{\hat T_1}(\hat\gamma(\hat T_1))$ is distributed like $g_{T_1}(\gamma(T_1))$.
Notice that the law of the hull $\hat K_{\hat T_1}$ is conformally invariant because,
by Lemma~\ref{hull}, it coincides with the law of the $SLE_6$ hull $K_{T_1}$.

Using now the Markovian character of $\hat X_j$, which implies that,
conditioned on $\hat X_1=(\hat K_{\hat T_1},\hat\gamma(\hat T_1))$,
$\hat K_{\hat T_2} \setminus \hat K_{\hat T_1}$ and $\hat\gamma(\hat T_2)$
are determined by Cardy's formula in $\hat G_2$, from the fact that
$\hat K_{\hat T_1}$ is equidistributed with $K_{T_1}$ and therefore $\hat G_2$
is equidistributed with $G_2$, we obtain that the hull $\hat K_{\hat T_2}$ is
distributed like $K_{T_2}$ and
%$\hat K_{\hat T_2} = \hat K_{\hat T_1} \cup \hat g^{-1}_{\hat T_1}(K + \hat g_{T_1}(\hat\gamma(\hat T_1)))$,
%where $K$ is a hull (independent of $\hat K_{\hat T_1}$) distributed like $\hat K_{\hat T_1}$ or $K_{T_1}$,
%so that the distribution of $\hat K_{\hat T_2}$ is the same as the distribution of
%$K_{T_1} \cup g^{-1}_{T_1}(K'_{T_1} + g_{T_1}(\gamma(T_1)))$.
its ``tip" $\hat\gamma(\hat T_2)$ is distributed like the ``tip" $\gamma(T_2)$
of the hull $K_{T_2}$.
We can then conclude that the joint distribution of
$\{ \hat\gamma(\hat T_1),\hat\gamma(\hat T_2) \}$ is the same as that of
$\{ \gamma(T_1),\gamma(T_2) \}$.
It also follows immediately that $\hat g_{\hat T_2}$ is equidistributed
with $g_{T_2}$ and $\hat\tau_2$ is equidistributed with $\tau_2$ or indeed
with $\tau_1$.

By repeating this recursively, using the Markovian character of the hulls and tips,
we obtain that, for all $j$, $\{ \hat\gamma(\hat T_1),\ldots,\hat\gamma(\hat T_j) \}$
is equidistributed with $\{ \gamma(T_1),\ldots,\gamma(T_j) \}$.
This immediately implies that $\hat\gamma_{\varepsilon}$ is equidistributed with
$\gamma_{\varepsilon}$.

To conclude the proof, we just have to show that, as $\varepsilon \to 0$,
$\hat\gamma_{\varepsilon}$ converges to $\hat\gamma$ and $\gamma_{\varepsilon}$
to $\gamma$ in the uniform metric~(\ref{distance}) on continuous curves.
This, however, follows easily from the properties of the continuous curves
we are considering (see the discussion after~(\ref{hulls})), if we can
show that $\hat T_{j+1} - \hat T_j = \hat\tau_{j+1}$ and
$T_{j+1}-T_j=\tau_{j+1}$ go to $0$ as $\varepsilon \to 0$.
To see this, we recall that $\hat\tau_{j+1}$ and $\tau_{j+1}$ are distributed like
$\tau_1$ and use Lemma~2.1 of~\cite{lsw6}, which implies the (deterministic) bound
$\tau_1(\varepsilon) \leq \varepsilon^2/2$, which follows from the well-known bound
$a(t) \leq \varepsilon^2$ for the
half-plane capacity $a(t) = 2 \, t$ of~(\ref{conf-map2}).~\fbox{}

\begin{remark} \label{remark-difference}
The procedure for constructing the polygonal approximations of $\hat\gamma$
and $\gamma$ and the recursive strategy for proving that they have the same
distribution include significant modifications to the sketched argument for
convergence of the percolation exploration process to chordal $SLE_6$ proposed
by Smirnov in~\cite{smirnov} and \cite{smirnov-long}.
One modification is that we use ``conformal semi-balls" instead of balls
(see~\cite{smirnov,smirnov-long}) to define the sequences of stopping times
$\{ \hat T_j \}$ and $\{ T_j \}$.
Since the paths we are dealing with touch themselves (or almost do), if one
were to use ordinary balls, some of them would intersect multiple disjoint
pieces of the past hull, making it impossible to use Cardy's formula in the
``triangular setting" proposed by (Carleson and) Smirnov and used here.
The use of conformally mapped semi-balls ensures, thanks to the choice of the
conformal maps, that the domains used to define the stopping times intersect
a single piece of the past hull.
This is a natural choice (exploiting the conformal invariance) to obtain
a good polygonal approximation of the paths while still being able to use
Cardy's formula to determine hitting distributions.
\end{remark}

%Smirnov's proof of the convergence of crossing probabilities to Cardy's
%formula implies that the exit distribution of the hull of any subsequential
%limit of the percolation exploration path is determined by Cardy's formula.
%Next, we want to show that the hull of any subsequential scaling limit
%of the percolation exploration path also satisfies the spatial Markov property.
%In order to do that, we first need an additional result, Theorem~\ref{strong-cardy}
%below, about the convergence of crossing probabilities to Cardy's formula in
%Jordan domains.

We will next prove a version of Smirnov's result (Theorem~\ref{cardy-smirnov}
above) extended to cover the convergence of crossing probabilities to Cardy's
formula for the case of sequences of admissible domains
(see the definition of admissible in Section~\ref{explo}).
The statement of Theorem~\ref{strong-cardy} below is certainly
not optimal, but it is sufficient for our purposes.
We remark that a weaker statement restricted, for instance, only
to Jordan domains would not be sufficient -- see
Figure~\ref{cut-point-fig} and the discussion referring to it in
the proof of Theorem~\ref{spatial-markov} below.

\begin{theorem} \label{strong-cardy}
Consider a sequence $\{(D_k,a_k,c_k,b_k,d_k)\}$ of domains $D_k$ containing
the origin, admissible with respect to the points $a_k,c_k,d_k$
on $\partial D_k$, and with $b_k$
belonging to the interior of the counterclockwise arc $\overline{c_k d_k}$
of $\partial D_k$.
Assume that, as $k \to \infty$, $b_k \to b$ and there is convergence
in the metric~(\ref{distance}) of the counterclockwise arcs
$\overline{d_k a_k}$, $\overline{a_k c_k}$, $\overline{c_k d_k}$ to the
corresponding counterclockwise arcs
$\overline{d a}$, $\overline{a c}$, $\overline{c d}$ of $\partial D$,
%$(D_k,a_k,c_k,d_k)$ converges to
%$(D,a,c,d)$ and $b_k$ to $b$,
where $D$ is a domain containing the origin,
admissible with respect to $(a,c,d)$, and $b$ belongs to the interior of
$\overline{cd}$.
Then, for any sequence $\delta_k \downarrow 0$, the probability
$\Phi^{\delta_k}_k (\equiv \Phi^{\delta_k}_{D_k})$ of a blue crossing
inside $D_k$ from $\overline{a_k c_k}$
to $\overline{b_k d_k}$
converges, as $k \to \infty$, to Cardy's formula
$\Phi_D$ (see~(\ref{cardy-formula})) for a blue crossing inside $D$ from
$\overline{a c}$ to
$\overline{c d}$.
\end{theorem}

\noindent {\bf Proof.} When $D$ is Jordan, one can use essentially the
same arguments as in Lemma~\ref{equal} below (see also Theorem~1 of~\cite{cns1})
to construct for each $\varepsilon>0$, $\varepsilon$-approximations,
$(\tilde D, \tilde a, \tilde b, \tilde c, \tilde d)$ and $(\hat D, \hat a, \hat b, \hat c, \hat d)$,
to $(D,a,b,c,d)$ so that
$\Phi_{\tilde D(\varepsilon)}, \Phi_{\hat D(\varepsilon)}
\stackrel{\varepsilon \to 0}{\longrightarrow} \Phi_D$
while
\begin{equation} \label{bounds}
\Phi_{\tilde D(\varepsilon)} = \liminf_{k \to \infty} \Phi^{\delta_k}_{\tilde D(\varepsilon)}
\leq \liminf_{k \to \infty} \Phi^{\delta_k}_{D_k} \leq \limsup_{k \to \infty} \Phi^{\delta_k}_{D_k}
\leq \limsup_{k \to \infty} \Phi^{\delta_k}_{\hat D(\varepsilon)} = \Phi_{\hat D(\varepsilon)}.
\end{equation}

When $D$ is not Jordan but admissible, one can do a similar construction
on a Riemann surface with a cut starting from $a$ to separate the touching
arcs $\overline{da}$ and $\overline{ac}$ of $\partial D$.
This is similar to an argument in~\cite{sw} replacing an annulus in the plane
by its universal cover. \fbox{}

\begin{remark}
A construction for the non-Jordan case of approximating $\tilde D$ and $\hat D$
without the use of a cut surface may be found in the appendix of~\cite{cn1}.
It has been suggested to us by a referee and by V.~Beffara that existing proofs
of convergence to Cardy's formula for fixed domains (see, in particular, \cite{beffara2})
should also work in the context of Theorem~\ref{strong-cardy}.
\end{remark}

%We note that Theorem~\ref{strong-cardy}, combined with the continuity
%of Cardy's formula in the shape of the domain (for admissible domains)
%and positions of the four points on the boundary (Lemma~\ref{cont-cardy}),
%implies the convergence of crossing probabilities to Cardy's formula
%\emph{locally uniformly} in the shape of the domain with respect to the
%uniform metric on curves, and in the location of the four points on the
%boundary with respect to the Euclidean metric; i.e., for $(D,a,b,c,d)$
%an admissible domain with $a,b,c,d \in \partial D$ (with the notation
%used in Theorem~\ref{strong-cardy}), $\forall \varepsilon>0$,
%$\exists \alpha_0=\alpha_0(\varepsilon)$ and $\delta_0=\delta_0(\varepsilon)$
%such that for all admissible domains $(D',a',b',c',d')$ with
%$\max{(\text{d}(\partial D, \partial D'),|a-a'|,|b-b'|,|c-c'|,|d-d'|) \leq \alpha_0}$
%and $\delta \leq \delta_0$,
%$|\Phi_{D'}(a',c';b',d') - \Phi^{\delta}_{D'}(a',c';b',d')| \leq \varepsilon$,
%where $\Phi_{D'}(a',c';b',d')$ is Cardy's formula and
%$\Phi^{\delta}_{D'}(a',c';b',d')$ is the corresponding crossing probability.

\section{Boundary of the Hull and the Scaling Limit} \label{boundary}

We give here some important results which are needed
in the proofs of the main theorems.
We start with two lemmas from~\cite{cn1,cn2}, which are
consequences of~\cite{ab}, of standard bounds on the
probability of events corresponding to having a certain
number of disjoint monochromatic crossings of an annulus
(see Lemma~5 of~\cite{ksz}, Appendix~A of~\cite{lsw5},
and also~\cite{ada}).
Afterwards we give two related lemmas that are more suited
to this paper and whose proofs are modified revisions of those
of the first two lemmas.

\begin{lemma} \label{sub-conv}
Let $\gamma^{\delta}_{{\mathbb D},-i,i}$ be the percolation
exploration path on the edges of $\delta {\cal H}$ inside
(a $\delta$-approximation of) $\mathbb D$ between
(e-vertices close to) $-i$ and $i$.
For any fixed point $z \in {\mathbb D}$, chosen independently
of $\gamma^{\delta}_{{\mathbb D},-i,i}$, as $\delta \to 0$,
$\gamma^{\delta}_{{\mathbb D},-i,i}$ and the boundary
$\partial {\mathbb D}^{\delta}_{-i,i}(z)$ of the domain
${\mathbb D}^{\delta}_{-i,i}(z)$ that contains $z$ jointly
have limits in distribution along subsequences of $\delta$
with respect to the uniform metric~(\ref{distance}) on
continuous curves.
Moreover, any subsequence limit of
$\partial {\mathbb D}^{\delta}_{-i,i}(z)$ is almost surely
a simple loop~\cite{ada}.
\end{lemma}

\noindent {\bf Proof.}
The first part of the lemma is a direct consequence
of~\cite{ab}; it is enough to notice that the (random)
polygonal curves $\gamma^{\delta}_{{\mathbb D},-i,i}$
and $\partial {\mathbb D}^{\delta}_{-i,i}(z)$
satisfy the conditions in~\cite{ab} and thus have a
scaling limit in terms of continuous curves, at least
along subsequences of $\delta$.

To prove the second part, we use standard percolation
bounds (see Lemma~5 of~\cite{ksz} and Appendix A
of~\cite{lsw5}) to show that, in the limit $\delta \to 0$,
the loop $\partial {\mathbb D}^{\delta}_{-i,i}(z)$ does
not collapse on itself but remains a  simple loop~\cite{ada}.

Let us assume that this is not the case and that
the limit $\tilde\gamma$ of
$\partial {\mathbb D}^{\delta_k}_{-i,i}(z)$ along some
subsequence $\{ \delta_k \}_{k \in {\mathbb N}}$ touches
itself, i.e., $\tilde\gamma(t_0)=\tilde\gamma(t_1)$ for
$t_0 \neq t_1$ with positive probability.
If so, we can take $\varepsilon>\varepsilon'>0$ so small that the annulus
$B(\tilde\gamma(t_1),\varepsilon) \setminus B(\tilde\gamma(t_1),\varepsilon')$
is crossed at least four times by $\tilde\gamma$
(here $B(u,r)$ is the ball of radius $r$ centered at $u$).

Because of the choice of topology, the convergence in
distribution of $\partial {\mathbb D}^{\delta_k}_{-i,i}(z)$
to $\tilde\gamma$ implies that we can find coupled versions of
$\partial {\mathbb D}^{\delta_k}_{-i,i}(z)$ and $\tilde\gamma$
%on the Lebesgue probability space $([0,1], {\cal B}[0,1], \lambda)$
%(where ${\cal B}[0,1]$ is the Borel $\sigma$-algebra of the
%unit interval and $\lambda$ the Lebesgue measure) such that
on some $(\Omega',{\cal B}',{\mathbb P}')$ such that
$\text{d}(\partial {\mathbb D}^{\delta}_{-i,i}(z),\tilde\gamma) \to 0$,
for all $\omega' \in \Omega'$ as $k \to \infty$ (see, e.g.,
Corollary~1 of~\cite{billingsley1}).

Using this coupling, we can choose $k$ large enough (depending
on $\omega'$) so that $\partial {\mathbb D}^{\delta_k}_{-i,i}(z)$
stays in an $\varepsilon'/2$-neighborhood
${\cal N}(\tilde\gamma,\varepsilon'/2) \equiv \bigcup_{u \in \tilde\gamma} B(u,\varepsilon'/2)$
of $\tilde\gamma$.
This however would correspond to an event
${\cal A}_{\tilde\gamma(t_1)}(\varepsilon,\varepsilon')$
that (at least) four
paths of one color (corresponding to the four crossings
by $\partial {\mathbb D}^{\delta_k}_{-i,i}(z)$) and two
of the other color cross the annulus
$B(\tilde\gamma(t_1),\varepsilon-\varepsilon'/2) \setminus
B(\tilde\gamma(t_1),3 \, \varepsilon'/2)$.
As $\delta_k \to 0$, we can let $\varepsilon' \to 0$
(keeping $\varepsilon$ fixed),
in which case, we claim that the
probability of seeing the event just described somewhere inside
$\mathbb D$ goes to zero,
%CHANGE (02-21-06): REFERENCES REMOVED
%SINCE THE RELEVANT ONE APPEARS BELOW %\cite{ksz,lsw5},
leading to a contradiction.
%
%CHANGE (02-21-06): PARAGRAPH MOVED HERE FROM PROOF OF LEMMA 6.2.
%Let us call ${\cal A}_w(\varepsilon,\varepsilon')$ the
%event described above, where $\gamma(t_1)=w$;
This is because a standard
bound~\cite{ksz} on the probability of six disjoint crossings
(not all of the same color) of an annulus gives that the
probability of ${\cal A}_w(\varepsilon,\varepsilon')$ scales as
$(\frac{\varepsilon'}{\varepsilon})^{2+\alpha}$ with $\alpha>0$.
As $\delta \to 0$, we can let $\varepsilon' \to 0$
(keeping $\varepsilon$ fixed); then
the probability of ${\cal A}_w(\varepsilon,\varepsilon')$
goes to zero sufficiently rapidly with $\varepsilon'$
to conclude that the probability to see such an event
anywhere in $\mathbb D$ goes to zero. \fbox{} \\

The second lemma states that, for every subsequence limit, the
discrete boundaries converge to the boundaries of the domains
generated by the limiting continuous curve.
In order to insure this, we need to show that whenever the discrete
exploration path comes at distance of order $\delta$ from the boundary
of the exploration domain or from its past filling, producing a ``fjord"
(see~\cite{ada}) and causing touching in the limit $\delta \to 0$, with
high probability the discrete path already closes the fjord by touching
the boundary of the exploration domain or %its past filling,
by ``touching" itself (i.e., getting to distance $\delta$
of itself, just one hexagon away), so that no
discrepancy arises, as $\delta \to 0$, between the limit of the discrete
filling and the filling of the limiting continuous curve.
This issue will come up in the proof of Theorem~\ref{spatial-markov}
below, and is one of the main technical issues of this paper.

\begin{lemma} \label{boundaries}
Using the notation of Lemma~\ref{sub-conv}, let
$\gamma_{{\mathbb D},-i,i}$ be the limit in distribution
of $\gamma^{\delta}_{{\mathbb D},-i,i}$ as $\delta \to 0$
along some convergent subsequence $\{ \delta_k \}$ and
$\partial {\mathbb D}_{-i,i}(z)$ be the boundary of the domain
${\mathbb D}_{-i,i}(z)$ of ${\mathbb D} \setminus \gamma_{D,-i,i}[0,1]$
that contains $z$.
Then, as $k \to \infty$,
$(\gamma^{\delta_k}_{{\mathbb D},-i,i},\partial {\mathbb D}^{\delta_k}_{-i,i}(z))$
converges in distribution to
$(\gamma_{{\mathbb D},-i,i},\partial {\mathbb D}_{-i,i}(z))$.
\end{lemma}

\noindent {\bf Proof.}
Let $\{ \delta_k \}_{k \in {\mathbb N}}$ be a convergent
subsequence for $\gamma^{\delta}_{{\mathbb D},-i,i}$
and $\gamma \equiv \gamma_{{\mathbb D},-i,i}$
the limit in distribution of $\gamma^{\delta_k}_{{\mathbb D},-i,i}$
as $k \to \infty$.
%the trace of chordal $SLE_6$ inside $\mathbb D$ from $-i$ to $i$.
For simplicity of notation, we now drop the $k$ and write $\delta$
instead of $\delta_k$.
Because of the choice of topology, the convergence in distribution
of $\gamma^{\delta} \equiv \gamma^{\delta}_{{\mathbb D},-i,i}$ to
$\gamma$ implies that we can find coupled versions of $\gamma^{\delta}$
and $\gamma$ on some probability space $(\Omega',{\cal B}',{\mathbb P}')$
%the Lebesgue probability space $([0,1], {\cal B}[0,1], \lambda)$
such that $\text{d}(\gamma^{\delta}(\omega'),\gamma(\omega')) \to 0$,
for all $\omega'$ as $k \to \infty$ (see, for example, Corollary~1
of~\cite{billingsley1}).
Using this coupling, our first task will be to prove the following claim:
\begin{itemize}
\item[(C)] For two (deterministic) points $u,v \in {\mathbb D}$,
the probability that ${\mathbb D}_{-i,i}(u) = {\mathbb D}_{-i,i}(v)$ but
${\mathbb D}^{\delta}_{-i,i}(u) \neq {\mathbb D}^{\delta}_{-i,i}(v)$
or vice versa goes to zero as $\delta \to 0$.
\end{itemize}
%As a result, we will get that for any $z \in D$,
%$\partial D_{a,b}^{\delta}(z)$ converges in distribution
%to $\partial D_{a,b}(z)$ in the topology induced by~(\ref{distance}).
%This will tell us that $\partial D_{a,b}(z)$ is a Jordan
%curve, and therefore that $D_{a,b}(z)$ is a Jordan domain.

Let us consider first the case of $u,v$ such that
${\mathbb D}_{-i,i}(u) = {\mathbb D}_{-i,i}(v)$ but
${\mathbb D}^{\delta}_{-i,i}(u) \neq {\mathbb D}^{\delta}_{-i,i}(v)$.
Since ${\mathbb D}_{-i,i}(u)$ is an open subset of $\mathbb C$,
there exists a continuous curve $\gamma_{u,v}$ joining
$u$ and $v$ and a constant $\varepsilon>0$ such that
the $\varepsilon$-neighborhood
${\cal N}(\gamma_{u,v},\varepsilon)$ of the curve is contained
in ${\mathbb D}_{-i,i}(u)$, which implies that $\gamma$ does
not intersect ${\cal N}(\gamma_{u,v},\varepsilon)$.
%by~(\ref{coupling}), that the probability that
%$\gamma^{\delta}$ does not intersect
%${\cal N}_{\varepsilon/2} (\gamma_{u,v})$ is close
%to one for $\delta$ small enough (depending on
%$\varepsilon$).
Now, if $\gamma^{\delta}$ does not intersect
${\cal N}(\gamma_{u,v},\varepsilon/2)$, for
$\delta$ small enough, then there is a $\cal T$-path
$\pi$ of unexplored hexagons connecting the hexagon
that contains $u$ with the hexagon that contains $v$,
and we conclude that
${\mathbb D}^{\delta}_{-i,i}(u) = {\mathbb D}^{\delta}_{-i,i}(v)$.

This shows that the event that ${\mathbb D}_{-i,i}(u) = {\mathbb D}_{-i,i}(v)$
but ${\mathbb D}^{\delta}_{-i,i}(u) \neq {\mathbb D}^{\delta}_{-i,i}(v)$
implies the existence of a curve $\gamma_{u,v}$ whose
$\varepsilon$-neighborhood ${\cal N}(\gamma_{u,v},\varepsilon)$
is not intersected by $\gamma$ but whose $\varepsilon/2$-neighborhood
${\cal N}(\gamma_{u,v},\varepsilon/2)$ is intersected by $\gamma^{\delta}$.
This implies that $\forall u,v \in {\mathbb D}$, $\exists \varepsilon>0$ such
that ${\mathbb P}'({\mathbb D}_{-i,i}(u) = {\mathbb D}_{-i,i}(v) \text{ but }
{\mathbb D}^{\delta}_{-i,i}(u) \neq {\mathbb D}^{\delta}_{-i,i}(v)) \leq
{\mathbb P}'(\text{d}(\gamma^{\delta},\gamma) \geq \varepsilon/2)$.
But the right hand side goes to zero for every $\varepsilon>0$ as
$\delta \to 0$, which concludes the proof of one direction of the claim.
%Because of the coupling between $\gamma$ and $\gamma^{\delta}$, for
%every $\omega'$ there is a $\delta_0=\delta_0(\omega')$ such that
%this does not happen for all $\delta < \delta_0$.
%This implies that, as $\delta \to 0$, the probability of the event
%described tends to zero, which proves one direction of the claim.

To prove the other direction, we consider two points
$u,v \in {\mathbb D}$ such that ${\mathbb D}_{-i,i}(u)
\neq {\mathbb D}_{-i,i}(v)$ but ${\mathbb D}^{\delta}_{-i,i}(u)
= {\mathbb D}^{\delta}_{-i,i}(v)$.
Assume that $u$ is trapped before $v$ by $\gamma$ and
suppose for the moment that ${\mathbb D}_{-i,i}(u)$ is
a domain of type 3 or 4 (as defined at the end of
Section~\ref{lapa}); the case of a domain of type 1 or 2
is analogous and will be treated later.
Let $t_1$ be the first time $u$ is trapped by
$\gamma$ with $\gamma(t_0)=\gamma(t_1)$ the double
point of $\gamma$ where the domain ${\mathbb D}_{-i,i}(u)$
containing $u$ is ``sealed off.''
At time $t_1$, a new domain containing $u$ is
created and $v$ is disconnected from $u$.

Choose $\varepsilon>0$ small enough so that neither $u$ nor
$v$ is contained in the ball $B(\gamma(t_1),\varepsilon)$
of radius $\varepsilon$ centered at $\gamma(t_1)$, nor in
the $\varepsilon$-neighborhood
${\cal N}(\gamma[t_0,t_1],\varepsilon)$ of the portion
of $\gamma$ which surrounds $u$. %$D_{-i,i}(u)$.
Then it follows from the coupling that, for $\delta$
small enough, there are appropriate parameterizations
of $\gamma$ and $\gamma^{\delta}$ such that the portion
$\gamma^{\delta}[t_0,t_1]$ of $\gamma^{\delta}(t)$ is
inside ${\cal N}(\gamma[t_0,t_1],\varepsilon)$, and
$\gamma^{\delta}(t_0)$ and $\gamma^{\delta}(t_1)$ are
contained in $B(\gamma(t_1),\varepsilon)$.

For $u$ and $v$ to be contained in the same domain
in the discrete construction, there must be a
$\cal T$-path $\pi$ of unexplored hexagons
connecting the hexagon that contains $u$ to the
hexagon that contains $v$.
From what we said in the previous paragraph, any
such $\cal T$-path connecting $u$ and $v$
would have to go though a ``bottleneck'' in
$B(\gamma(t_1),\varepsilon)$.

Assume now, for concreteness but without loss of
generality, that ${\mathbb D}_{-i,i}(u)$ is a domain
of type~3, which means that $\gamma$ winds around
$u$ counterclockwise, and consider the hexagons
to the ``left" of $\gamma^{\delta}[t_0,t_1]$.
Those hexagons form a ``quasi-loop'' around $u$
since they wind around it (counterclockwise) and
the first and last hexagons are both contained in
$B(\gamma(t_1),\varepsilon)$.
The hexagons to the left of $\gamma^{\delta}[t_0,t_1]$
belong to the set $\Gamma_Y(\gamma^{\delta})$, which
can be seen as a (nonsimple) path by connecting the
centers of the hexagons in $\Gamma_Y(\gamma^{\delta})$
by straight segments.
Such a path shadows $\gamma^{\delta}$, with the difference
that it can have double (or even triple) points, since
the same hexagon can be visited more than once.
Consider $\Gamma_Y(\gamma^{\delta})$ as a path
$\hat\gamma^{\delta}$ with a given parametrization
$\hat\gamma^{\delta}(t)$, chosen so that
$\hat\gamma^{\delta}(t)$ is inside
$B(\gamma(t_1),\varepsilon)$ when $\gamma^{\delta}(t)$ is,
and it winds around $u$ together with $\gamma^{\delta}(t)$.

Now suppose that there were two times, $\hat t_0$ and
$\hat t_1$, such that $\hat\gamma^{\delta}(\hat t_1)
= \hat\gamma^{\delta}(\hat t_0) \in B(\gamma(t_1),\varepsilon)$
and $\hat\gamma^{\delta}[\hat t_0,\hat t_1]$ winds
around $u$.
This would imply that the ``quasi-loop'' of explored
yellow hexagons around $u$ is actually completed, and
that ${\mathbb D}^{\delta}_{-i,i}(v) \neq {\mathbb D}^{\delta}_{-i,i}(u)$.
Thus, for $u$ and $v$ to belong to the same discrete
domain, this cannot happen.

For any $0<\varepsilon'<\varepsilon$, if we take $\delta$
small enough, $\hat\gamma^{\delta}$ will be contained
inside ${\cal N}(\gamma,\varepsilon')$, due to the coupling.
Following the considerations above, the fact that $u$
and $v$ belong to the same domain in the discrete
construction but to different domains in the continuum
construction implies, for $\delta$ small enough, that
there are four disjoint yellow $\cal T$-paths
crossing the annulus $B(\gamma(t_1),\varepsilon)
\setminus B(\gamma(t_1),\varepsilon')$ (the paths
have to be disjoint because, as we said,
$\hat\gamma^{\delta}$ cannot, when coming back to
$B(\gamma(t_1),\varepsilon)$ after winding around
$u$, touch itself inside $B(\gamma(t_1),\varepsilon)$).
Since $B(\gamma(t_1),\varepsilon) \setminus
B(\gamma(t_1),\varepsilon')$ is also crossed
by at least two blue $\cal T$-paths from
$\Gamma_B(\gamma^{\delta})$, there is a total
of at least six $\cal T$-paths, not all of
the same color, crossing the annulus
$B(\gamma(t_1),\varepsilon) \setminus B(\gamma(t_1),\varepsilon')$.
As $\delta \to 0$, we can let $\varepsilon' \to 0$
(keeping $\varepsilon$ fixed) and conclude, like in the
proof of Lemma~\ref{sub-conv}, that the probability to
see such an event anywhere in $\mathbb D$ goes to zero.

%Let us call ${\cal A}_w(\varepsilon,\varepsilon')$ the
%event described above, where $\gamma(t_1)=w$; a standard
%bound~\cite{ksz} on the probability of six disjoint crossings
%(not all of the same color) of an annulus gives that the
%probability of ${\cal A}_w(\varepsilon,\varepsilon')$ scales as
%$(\frac{\varepsilon'}{\varepsilon})^{2+\alpha}$ with $\alpha>0$.
%As $\delta \to 0$, we can let $\varepsilon'$ go to
%zero (keeping $\varepsilon$ fixed); when we do this,
%the probability of ${\cal A}_w(\varepsilon,\varepsilon')$
%goes to zero sufficiently rapidly with $\varepsilon'$
%to conclude, like in the proof of Lemma~\ref{sub-conv},
%that the probability to see such an event anywhere in
%$\mathbb D$ goes to zero.

%As $\delta \to 0$, we can let $\varepsilon'$ go to
%zero and use a standard bound~\cite{ksz} to show that
%the probability of six disjoint crossings of
%$B_{\varepsilon}(\gamma(t_1))
%\setminus B_{\varepsilon'}(\gamma(t_1))$ not all
%of the same color goes to zero as $\varepsilon' \to 0$
%(for $\varepsilon$ fixed).

In the case in which $u$ belongs to a domain of type~1
or 2, let $\cal E$ be the excursion that traps $u$
and $\gamma(t_0) \in \partial {\mathbb D}$ be the point
on the boundary of $\mathbb D$ where $\cal E$ starts and
$\gamma(t_1) \in \partial {\mathbb D}$ the point where
it ends.
Choose $\varepsilon>0$ small enough so that neither $u$
nor $v$ is contained in the balls
$B(\gamma(t_0),\varepsilon)$ and $B(\gamma(t_1),\varepsilon)$
of radius $\varepsilon$ centered at $\gamma(t_0)$ and
$\gamma(t_1)$, nor in the $\varepsilon$-neighborhood
${\cal N}({\cal E},\varepsilon)$ of the excursion $\cal E$.
Because of the coupling, for $\delta$ small enough
(depending on $\varepsilon$), $\gamma^{\delta}$ shadows
$\gamma$ along $\cal E$, staying within
${\cal N}({\cal E},\varepsilon)$.
If this is the case, any $\cal T$-path of unexplored
hexagons connecting the hexagon that contains $u$ with
the hexagon that contains $v$ would have to go through
one of two ``bottlenecks,'' one contained in
$B(\gamma(t_0),\varepsilon)$ and the other in
$B(\gamma(t_1),\varepsilon)$.

Assume for concreteness (but without loss of generality)
that $u$ is in a domain of type~1, which means that $\gamma$
winds around $u$ counterclockwise.
If we parameterize $\gamma$ and $\gamma^{\delta}$ so that
$\gamma^{\delta}(t_0) \in B(\gamma(t_0),\varepsilon)$
and $\gamma^{\delta}(t_1) \in B(\gamma(t_1),\varepsilon)$,
$\gamma^{\delta}[t_0,t_1]$ forms a ``quasi-excursion''
around $u$ since it winds around it (counterclockwise)
and it starts inside $B_{\varepsilon}(\gamma(t_0))$ and
ends inside $B_{\varepsilon}(\gamma(t_1))$.
Notice that if $\gamma^{\delta}$ touched
$\partial {\mathbb D}^{\delta}$, inside both
$B_{\varepsilon}(\gamma(t_0))$ and
$B_{\varepsilon}(\gamma(t_1))$, this would imply
that the ``quasi-excursion'' is a real excursion and
that ${\mathbb D}^{\delta}_{-i,i}(v) \neq {\mathbb D}^{\delta}_{-i,i}(u)$.

For any $0<\varepsilon'<\varepsilon$, if we take $\delta$
small enough, $\gamma^{\delta}$ will be contained inside
${\cal N}(\gamma,\varepsilon')$, due to the coupling.
Therefore, the fact that
${\mathbb D}^{\delta}_{-i,i}(v) = {\mathbb D}^{\delta}_{-i,i}(u)$
implies, with probability going to one as $\delta \to 0$, that for
$\varepsilon>0$ fixed and any $0<\varepsilon'<\varepsilon$, $\gamma^{\delta}$
enters the ball $B(\gamma(t_i),\varepsilon')$ and does not touch
$\partial {\mathbb D}^{\delta}$ inside the larger ball
$B(\gamma(t_i),\varepsilon)$, for $i=0$ or $1$.
This is equivalent to having at least two yellow and one
blue $\cal T$-paths (contained in ${\mathbb D}^{\delta}$)
crossing the annulus
$B(\gamma(t_i),\varepsilon) \setminus B(\gamma(t_i),\varepsilon')$.
Let us call ${\cal B}_w(\varepsilon,\varepsilon')$ the
event described above, where $\gamma(t_i)=w$; a standard
bound~\cite{lsw5} (this bound can also be derived from the
one obtained in~\cite{ksz}) on the probability of disjoint
crossings (not all of the same color) of a semi-annulus in
the upper half-plane gives that the probability of
${\cal B}_w(\varepsilon,\varepsilon')$ scales as
$(\frac{\varepsilon'}{\varepsilon})^{1+\beta}$ with $\beta>0$.
(We can apply the bound to our case because the unit disc
is a convex subset of the half-plane $\{ x+iy:y>-1 \}$
and therefore the intersection of an annulus centered at
say $-i$ with the unit disc is a subset of the intersection
of the same annulus with the half-plane $\{ x+iy:y>-1 \}$.)
As $\delta \to 0$, we can let $\varepsilon' \to 0$
(keeping $\varepsilon$ fixed), concluding that the probability
that such an event occurs anywhere on the boundary of the disc
goes to zero.

We have shown that, for two fixed points $u,v \in {\mathbb D}$,
having ${\mathbb D}_{-i,i}(u) \neq {\mathbb D}_{-i,i}(v)$ but
${\mathbb D}^{\delta}_{-i,i}(u) = {\mathbb D}^{\delta}_{-i,i}(v)$
or vice versa implies the occurrence of an event whose
probability goes to zero as $\delta \to 0$, and the
proof of the claim is concluded.

The Hausdorff distance %$\text{d}_{\text{H}}(A,B)$
between two closed nonempty subsets of $\overline{\mathbb D}$ is
\begin{equation} \label{hausdorff-dist}
\text{d}_{\text{H}}(A,B) \equiv \inf \{ \ell \geq 0 : B \subset
\cup_{a \in A} B(a,\ell), \, A \subset
\cup_{b \in B} B(b,\ell) \}.
\end{equation}
With this metric, the collection of closed subsets of
$\overline{\mathbb D}$ is a compact space.
We will next prove that
$\partial {\mathbb D}^{\delta}_{-i,i}(z)$ converges
in distribution to $\partial {\mathbb D}_{-i,i}(z)$
as $\delta \to 0$, in the topology induced
by~(\ref{hausdorff-dist}).
(Notice that the coupling between $\gamma^{\delta}$
and $\gamma$ provides a coupling between
$\partial {\mathbb D}^{\delta}_{-i,i}(z)$ and
$\partial {\mathbb D}_{-i,i}(z)$, seen as
boundaries of domains produced by the two paths.)

%We will show that for any $\varepsilon,\varepsilon'>0$,
%there exists a $delta_0$ such that for all
%$\delta<\delta_0$,
%\begin{equation}
%\nu_{\delta_k}
%[\text{Dist}(\partial {\mathbb D}^{\delta}_{-i,i}(z),
%\partial {\mathbb D}_{-i,i}(z) > \varepsilon] < \varepsilon'.
%\end{equation}

We will now use Lemma~\ref{sub-conv} and take a further subsequence
$k_n$ of the $\delta$'s that for simplicity of notation we denote
by $\{ \delta_n \}_{n \in {\mathbb N}}$ such that, as $n \to \infty$,
$\{ \gamma^{\delta_n},\partial {\mathbb D}^{\delta_n}_{-i,i}(z) \}$
converge jointly in distribution to $\{ \gamma,\tilde\gamma \}$,
where $\tilde\gamma$ is a simple loop.
For any $\varepsilon>0$, since $\tilde\gamma$ is a compact set,
we can find a covering of $\tilde \gamma$ by a finite number
of balls of radius $\varepsilon/2$ centered at points on $\tilde\gamma$.
Each ball contains both points in the interior
$\text{int}(\tilde\gamma)$ of $\tilde\gamma$ and in the exterior
$\text{ext}(\tilde\gamma)$ of $\tilde \gamma$, and we can choose
(independently of $n$) one point from $\text{int}(\tilde\gamma)$
and one from $\text{ext}(\tilde\gamma)$ inside each ball.

Once again, the convergence in distribution of
$\partial {\mathbb D}^{\delta_n}_{-i,i}(z)$
to $\tilde\gamma$ implies the existence of a coupling
such that, for $n$ large enough, the selected points
that are in $\text{int}(\tilde\gamma)$ are contained
in ${\mathbb D}^{\delta_n}_{-i,i}(z)$, and those that
are in $\text{ext}(\tilde\gamma)$ are contained in the
complement of $\overline{{\mathbb D}^{\delta_n}_{-i,i}(z)}$.
But by claim (C), each one of the selected points that is contained
in ${\mathbb D}^{\delta_n}_{-i,i}(z)$ is also contained in
${\mathbb D}_{-i,i}(z)$ with probability going to $1$ as $n \to \infty$;
analogously, each one of the selected points contained in the complement
of $\overline{{\mathbb D}^{\delta_n}_{-i,i}(z)}$ is also contained in the
complement of $\overline{{\mathbb D}_{-i,i}(z)}$ with probability
going to $1$ as $n \to \infty$.
This implies that $\partial {\mathbb D}_{-i,i}(z)$ crosses each
one of the balls in the covering of $\tilde\gamma$, and therefore
$\tilde\gamma \subset \cup_ {u \in \partial {\mathbb D}_{-i,i}(z)} B(u,\varepsilon)$.
From this and the coupling between
$\partial {\mathbb D}^{\delta_n}_{-i,i}(z)$ and $\tilde\gamma$,
it follows immediately that, for $n$ large enough,
$\partial {\mathbb D}^{\delta_n}_{-i,i}(z) \subset
\cup_ {u \in \partial {\mathbb D}_{-i,i}(z)} B(u,\varepsilon)$
with probability close to one.

A similar argument (analogous to the previous one but simpler,
since it does not require the use of $\tilde\gamma$), with the roles
of ${\mathbb D}^{\delta_n}_{-i,i}(z)$ and ${\mathbb D}_{-i,i}(z)$
inverted, shows that $\partial {\mathbb D}_{-i,i}(z) \subset
\cup_ {u \in \partial {\mathbb D}^{\delta_n}_{-i,i}(z)} B(u,\varepsilon)$
with probability going to $1$ as $n \to \infty$.
Therefore, for all $\varepsilon>0$,
${\mathbb P}(\text{d}_{\text{H}}(\partial {\mathbb D}^{\delta_n}_{-i,i}(z),
\partial {\mathbb D}_{-i,i}(z)) > \varepsilon) \to 0$ as $n \to \infty$,
which implies convergence in distribution of
$\partial {\mathbb D}^{\delta_n}_{-i,i}(z)$ to
$\partial {\mathbb D}_{-i,i}(z)$, as $\delta_n \to 0$,
in the topology of~(\ref{hausdorff-dist}).
But Lemma~\ref{sub-conv} implies that
$\partial {\mathbb D}^{\delta_n}_{-i,i}(z)$
converges in distribution (using~(\ref{distance})) to a
simple loop; therefore $\partial {\mathbb D}_{-i,i}(z)$
must also be simple and we have convergence in the topology
of~(\ref{distance}).

It is also clear that the argument above is independent
of the subsequence $\{ \delta_n \}$, so the limit of
$\partial {\mathbb D}^{\delta}_{-i,i}(z)$ is unique and
coincides with $\partial {\mathbb D}_{-i,i}(z)$.
Hence, we have convergence in distribution of
$\partial {\mathbb D}^{\delta}_{-i,i}(z)$ to
$\partial {\mathbb D}_{-i,i}(z)$, as $\delta \to 0$,
in the topology of~(\ref{distance}), and indeed joint convergence
of $(\gamma^{\delta},\partial {\mathbb D}^{\delta}_{-i,i}(z))$
to $(\gamma,\partial {\mathbb D}_{-i,i}(z))$. \fbox{} \\

We next give two new lemmas which mostly follow from the previous
ones (or their proofs) but are more suitable for the purposes of
this paper.
Let $D$ be a Jordan domain, with $a$ and $b$ two distinct points on $\partial D$,
and consider Jordan sets $D^{\delta}$, for $\delta > 0$, from $\delta \cal H$
such that $(D^{\delta},a^{\delta},b^{\delta}) \to (D,a,b)$, as $\delta \to 0$,
where $a^{\delta},b^{\delta} \in \partial D^{\delta}$ are two distinct e-vertices
on $\partial D^{\delta}$.
This means that $(D^{\delta},a^{\delta},b^{\delta})$ is a $\delta$-approximation
of $(D,a,b)$.
Denote by $\gamma^{\delta}_{D,a,b}$ the percolation exploration path inside
$D^{\delta}$ from $a^{\delta}$ to $b^{\delta}$.

Let $f$ be a conformal map from the upper half-plane $\mathbb H$
to the Jordan domain $D$ and assume that $f^{-1}(a)=0$ and
$f^{-1}(b)=\infty$.
(Since $D$ is a Jordan domain, the map $f^{-1}$ has a continuous extension from $D$
to $D \cup \partial D$ -- see Theorem~\ref{cara-thm} of Appendix~\ref{rado} -- and,
by a slight abuse of notation, we do not distinguish between $f^{-1}$ and its extension;
the same applies to $f$.)
Denote by $C(0,\varepsilon) = \{ z:|z|<\varepsilon \} \cap {\mathbb H}$ the
semi-ball of radius $\varepsilon$ centered at the origin of the real line.
Let $G \equiv f(C(0,\varepsilon))$, $c' \equiv f(\varepsilon)$,
$d' \equiv f(-\varepsilon)$.
Also denote by $\partial^* G$ the following subset of the boundary of
$G$: $\partial^* G \equiv f(\{ z:|z|=\varepsilon \} \cap {\mathbb H})$.

Analogously, let $f^{\delta}$ be a conformal map from the upper half-plane $\mathbb H$
to the Jordan set $D^{\delta}$, assume that $(f^{\delta})^{-1}(a^{\delta})=0$ and
$(f^{\delta})^{-1}(b^{\delta})=\infty$, and define $G^{\delta} \equiv f^{\delta}(C(0,\varepsilon))$
and $\partial^* G^{\delta} \equiv f^{\delta}(\{ z:|z|=\varepsilon \} \cap {\mathbb H})$.
Note that since $\partial D^{\delta} \to \partial D$, by an application of
Corollary~\ref{cor2-unif-conv} of Appendix~\ref{rado}, we can and do choose $f^{\delta}$
so that it converges to $f$ uniformly in $\overline{\mathbb H}$.
(We remark that the full strength of Corollary~\ref{cor2-unif-conv} is not needed here
since we are dealing with Jordan domains).
With this choice, $\partial^* G^{\delta} \to \partial^* G$ in the metric~(\ref{distance}).

%Let $\gamma^{\delta}_{D,a,b}$ be the percolation exploration path
%on the edges of $\delta {\cal H}$ inside (a $\delta$-approximation
%of) $D$ between (e-vertices close to) $a$ and $b$,
Let $T^{\delta}$
be the first time that
$\gamma^{\delta}_{D,a,b}$
intersects $\partial^* G^{\delta}$,
%FEDERICO: HERE I FOUND AN INCONSISTENCY BETWEEN THE TEX FILE AND THE PDF FILE --
%THIS IS THE VERSION IN THE TEX FILE
and let $K^{\delta}_{T^{\delta}}$ be the ({\bf discrete})
{\bf filling} of $\gamma^{\delta}_{D,a,b}[0,T^{\delta}]$, i.e., the union
of the hexagons explored up to time $T^{\delta}$ and those unexplored
hexagons from which it is not possible to reach $b$ without crossing an
explored hexagon or $\partial D$ (in other words, this is the set of
hexagons that at time $T^{\delta}$ have been explored or are disconnected
from $b$ by the exploration path).
Notice that even though time variables appear explicitly in the next two
lemmas, the time parametrizations of the curves are irrelevant and do not
need to be specified.

\begin{lemma} \label{new-sub-conv}
With the above notation, as $\delta \to 0$,
$\gamma^{\delta}_{D,a,b}$ and the boundary of $K^{\delta}_{T^{\delta}}$
jointly have limits in distribution along subsequences of $\delta$ with
respect to the uniform metric~(\ref{distance}) on continuous curves.
Moreover, any subsequence limit of $K^{\delta}_{T^{\delta}}$ is almost
surely a hull that touches $\partial^* G$ at a single point.
\end{lemma}

\begin{figure}[!ht]
\begin{center}
\includegraphics[width=8cm]{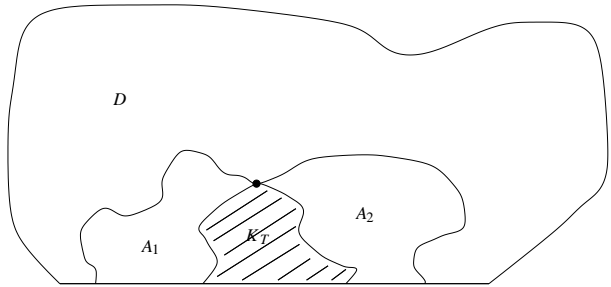}
\caption{Schematic figure representing $G \setminus K_T = A_1 \cup A_2$.}
\label{fig-sec5}
\end{center}
\end{figure}

\noindent {\bf Proof.} As in Lemma~\ref{sub-conv}, the first part
of this lemma is a direct consequence of~\cite{ab}.
The fact that the scaling limit of $K^{\delta_k}_{T^{\delta_k}}$ along
any convergent subsequence $\delta_k \downarrow 0$ touches
$\partial^* G$ at a single point, is a consequence of Lemma~\ref{double-crossing}
of Section~\ref{convergence} below.
(For any fixed $k$, the statement that $K^{\delta_k}_{T^{\delta_k}}$
touches $\partial^* G^{\delta_k}$ at a single point is a consequence of the
definition of the stopping time $T^{\delta_k}$, but a priori, this
could fail to be true in the limit $k \to \infty$.)
Therefore, if we remove that single point, the scaling limit of the
boundary of $K^{\delta_k}_{T^{\delta_k}}$ splits into a left and a right
part (corresponding to the scaling limit of the leftmost yellow and the
rightmost blue $\cal T$-paths of hexagons explored by
$\gamma^{\delta_k}_{D,a,b}$, respectively) each of which does not touch
$\partial^* G$ (like in Figure~\ref{fig-sec5}).

Moreover, Lemmas~\ref{double-crossing} and~\ref{close-encounters} below
(with $\tilde D=D$, $\hat D=G$ and in the limit where the target region
$J'=\partial^* G$ is fixed while $J \to \partial D \setminus \{a\}$)
imply that if $\gamma^{\delta_k}_{D,a,b}$ has a ``close encounter" with
$\partial D^{\delta_k}$, then the fjord produced by $\gamma^{\delta_k}_{D,a,b}$
is closed nearby with probability going to $1$.
%then it ``touches" $\partial D$ (or at least touches another segment
%of $\gamma^{\delta_k}_{D,a,b}$ that touches $\partial D$ nearby, which
%suffices for our purposes).
Analogously, the standard bound on the probability of six crossings
of an annulus~\cite{ksz}, used repeatedly before, implies that wherever
$\gamma^{\delta_k}_{D,a,b}$ has a ``close encounter" with itself, there
is ``touching" (see the proof of Lemma~\ref{sub-conv}).
These two observations assure that the complement of the scaling limit
of $K^{\delta_k}_{T^{\delta_k}}$ is almost surely connected, which means
that the scaling limit of $K^{\delta_k}_{T^{\delta_k}}$ has (almost surely)
the properties of a filling.
From the same bound on the probability of six crossings of an annulus,
we can also conclude that the scaling limits of the left and right
boundaries of $K^{\delta_k}_{T^{\delta_k}}$ are almost surely simple,
as in the proof of Lemma~\ref{sub-conv}.

It is also possible to conclude that the intersection of the scaling limit
of the left and right boundaries of $K^{\delta_k}_{T^{\delta_k}}$ with the
boundary of $D$ almost surely does not contain arcs of positive length.
In fact, if that were the case, one could find a subdomain $D'$ with three
points $z_1,z_2,z_3$ in counterclockwise order on $\partial D'$ such that
the probability that an exploration path started at $z_1$ and stopped
when it first hits the arc $\overline{z_2 z_3}$ of $\partial D'$ has a
positive probability, in the scaling limit, of hitting at $z_2$ or $z_3$,
contradicting Cardy's formula (which, by Theorem~\ref{strong-cardy}, holds
for all subsequential scaling limits).
Thus the scaling limit of $K^{\delta_k}_{T^{\delta_k}}$ almost
surely satisfies the condition in~(\ref{hulls}) and is therefore a hull. \fbox{}

\begin{lemma} \label{new-boundaries}
Using the notation of Lemma~\ref{new-sub-conv},
let $\gamma_{D,a,b}$ be the limit in distribution of
$\gamma^{\delta}_{D,a,b}$ as $\delta \to 0$ along some
convergent subsequence $\{ \delta_k \}$.
Denote by $T$ the first time that $\gamma_{D,a,b}$ exits
$G$ and by $K_T$ the filling of $\gamma_{D,a,b}[0,T]$.
Then, as $k \to \infty$,
$(\gamma^{\delta_k}_{D,a,b},K^{\delta_k}_{T^{\delta_k}})$
converges in distribution to $(\gamma_{D,a,b},K_T)$.
Moreover, $\gamma_{D,a,b}$ satisfies the properties (i)--(iii)
stated after~(\ref{hulls}) above, and its hull $K_T$ is
equidistributed with that of chordal $SLE_6$ at the
corresponding stopping time.
\end{lemma}

\noindent {\bf Proof.} Let $A^{\delta_k}_1$ and $A^{\delta_k}_2$
be the two domains of
$G^{\delta_k} \setminus K^{\delta_k}_{T^{\delta_k}}$,
and $A_1$ and $A_2$ the two domains of $G \setminus K_T$ (see
Figure~\ref{fig-sec5}).
Since hulls are characterized by their ``envelope" (see Lemma~\ref{hull}
and the discussion preceding it), the joint convergence in distribution
of $\{ \partial A^{\delta_k}_1,\partial A^{\delta_k}_2 \}$ to
$\{ \partial A_1,\partial A_2 \}$ would be enough to
conclude that $K^{\delta_k}_{T^{\delta_k}}$ converges to $K_T$ as $k \to \infty$,
and in fact that $(\gamma^{\delta_k}_{D,a,b},K^{\delta_k}_{T^{\delta_k}})$
converges in distribution to $(\gamma_{D,a,b},K_T)$.

In order to obtain the convergence of $\{ \partial A^{\delta_k}_1,\partial A^{\delta_k}_2 \}$,
we can use the convergence in distribution of $\gamma^{\delta_k}_{D,a,b}$
to $\gamma_{D,a,b}$ and apply almost the same arguments as used in
the proof of Lemma~\ref{boundaries}.
In fact, the domains $A^{\delta_k}_1,A^{\delta_k}_2$ and $A_1,A_2$
are of the same type as those treated in Lemma~\ref{boundaries}.
We just need to extend the definitions of the domains $D^{\delta}_{a,b}(z)$
and $D_{a,b}(z)$, as at the end of Section~\ref{lapa} and in Lemma~\ref{boundaries}
(where $(D,a,b)$ was taken to be (${\mathbb D},-i,i$)), to cover the case
in which the domain $D$ is replaced by a subset $G$ of $D$ and the target
point $b$ on the boundary of $D$ by an arc $\partial^* G$ of the boundary of $G$.
In our case, the subdomain $G$ and the arc $\partial^* G$ are defined just
before Lemma~\ref{new-sub-conv}.

The definitions are as before but with a deterministic target point
replaced by the random hitting point at the stopping time, i.e., we define
$G^{\delta}_{a,\partial^* G}(z) \equiv G^{\delta}_{a,\gamma^{\delta}_{D,a,b}(T^{\delta})}(z)$
and $G_{a,\partial^* G}(z) \equiv G_{a,\gamma_{D,a,b}(T)}(z)$.
$A^{\delta}_i$ (resp., $A_i$) for $i=1,2$ is a domain of type $G^{\delta}_{a,\partial^* G}(z_i)$
(resp., $G_{a,\partial^* G}(z_i)$) for some $z_i \in G$.
%
%We just need to modify the definition at the end of Section~\ref{lapa}
%of domains $D^{\delta}_{a,b}(z)$ used in Lemma~\ref{boundaries} to include
%$A^{\delta}_1$ and $A^{\delta}_2$ by considering the set
%$\Gamma(\gamma^{\delta}_{D,a,b}[0,T^{\delta}]) =
%\Gamma_Y(\gamma^{\delta}_{D,a,b}[0,T^{\delta}]) \cup
%\Gamma_B(\gamma^{\delta}_{D,a,b}[0,T^{\delta}])$.
%The set $G^{\delta} \setminus \Gamma(\gamma^{\delta}_{D,a,b}[0,T^{\delta}])$
%is the union of its connected components (in the lattice sense), which
%are simply connected.
%Given a point $z \in {\mathbb C}$ contained in
%$G^{\delta} \setminus \Gamma(\gamma^{\delta}_{D,a,b}[0,T^{\delta}])$,
%we will denote by $D^{\delta}_{a,G}(z)$ the domain corresponding to the
%unique element of $G^{\delta} \setminus \Gamma(\gamma^{\delta}_{D,a,b}[0,T^{\delta}])$
%that contains $z$ (notice that for a deterministic $z \in G$, $D^{\delta}_{a,G}(z)$
%is well defined with high probability for $\delta$ small, i.e., when
%$z \in G^{\delta}$ and $z \notin \Gamma(\gamma[0,T^{\delta}])$).
%Analogously, we will denote by $D_{a,G}(z)$ the domain corresponding to
%the unique element of $G \setminus \gamma_{D,a,b}[0,T]$ that contains $z$.
%
With these definitions, we need to prove the following claim.
\begin{itemize}
\item[($\text{C}^{\prime}$)] For two (deterministic) points $u,v \in G$, the probability
that $D_{a,\partial^* G}(u) = D_{a,\partial^* G}(v)$ but
$D^{\delta}_{a,\partial^* G}(u) \neq D^{\delta}_{a,\partial^* G}(v)$
or vice versa goes to zero as $\delta \to 0$.
\end{itemize}

The proof is the same as that of claim (C) in Lemma~\ref{boundaries},
except that here we cannot use the bound on the probability of three
crossings of an annulus centered at a boundary point because we may
not have a convex domain.
To replace that bound we use once again Lemmas~\ref{double-crossing}
and~\ref{close-encounters} below (as in the proof of Lemma~\ref{new-sub-conv}).

We have proved the convergence in distribution of
$(\gamma^{\delta_k}_{D,a,b},K^{\delta_k}_{T^{\delta_k}})$ to
$(\gamma_{D,a,b},K_T)$. As a consequence of that
and Smirnov's result on the convergence of
crossing probabilities (see Theorem~\ref{strong-cardy}), the hitting
distribution of $\gamma_{D,a,b}$ at the stopping time $T$ is determined
by Cardy's formula, which allows us to apply Lemma~\ref{hull} to conclude
that $K_T$ is equidistributed with the hull of chordal $SLE_6$ at the
corresponding stopping time.

It remains to prove (i)--(iii) stated after~(\ref{hulls}).
%To prove (i), note that the filling of a continuous curve is a hull as
%long as the curve does not ``stick" to the boundary of the domain.
%This follows from Cardy's formula as in the proof in Lemma~\ref{new-sub-conv}
%that any limit of $K^{\delta}_{T^\delta}$ is a hull.
%The proof of property (i) is now essentially the same as the proof in
%Lemma~\ref{new-sub-conv} that any limit of $K^{\delta}_{T^\delta}$ is a hull.
Property (i) is immediate in our case. Properties (ii) and (iii) are
consequences of Lemmas \ref{double-crossing}, \ref{close-encounters} and
six arm estimates. \fbox{}

\section{Convergence of the Exploration Path} \label{convergence}

Next we show that the filling of any subsequential scaling limit
of the percolation exploration process satisfies the spatial Markov property.
Let us start with some notation.
First, suppose that $D$ is Jordan, with $a$ and $b$ distinct
points on $\partial D$, and consider a sequence $D^{\delta}$
%$\{(D^{\delta},a^{\delta},b^{\delta})\}$
of Jordan sets  from $\delta \cal H$
such that $(D^{\delta},a^{\delta},b^{\delta}) \to (D,a,b)$, as
$\delta \to 0$, where $a^{\delta},b^{\delta} \in \partial D^{\delta}$ are two
distinct e-vertices on $\partial D^{\delta}$.
This means that $(D^{\delta},a^{\delta},b^{\delta})$ is a $\delta$-approximation
of $(D,a,b)$.
As before, denote by $\gamma^{\delta}_{D,a,b}$
%$\gamma^{\delta}_{D,a,b} \equiv
%\delta \gamma^1_{\frac{1}{\delta}D^{\delta},\frac{1}{\delta}a^{\delta},\frac{1}{\delta}b^{\delta}}$
the percolation exploration path inside $D^{\delta}$ from $a^{\delta}$ to
$b^{\delta}$.

%Our next task is to show that the filling of any subsequential scaling limit
%of the percolation exploration process satisfies the spatial Markov property.
%Let us start with some notation.
%First of all, suppose that $D$ is a Jordan domain, with $a$ and $b$ two distinct
%points on $\partial D$, and consider a sequence $\{(D_k,a_k,b_k)\}$ of Jordan
%domains such that $(D_k,a_k,b_k) \to (D,a,b)$, where $a_k,b_k \in \partial D_k$
%are two distinct points on $\partial D_k$.
%Denote by $\gamma^{\delta}_k \equiv \gamma^{\delta}_{D_k,a_k,b_k}$
%the percolation exploration path inside (a $\delta$-approximation
%$D^{\delta}_k$ of) $D_k$ from (some e-vertex close to) $a_k$ to
%(some e-vertex close to) $b_k$.

%CHANGE (02-05-06)? DELETE NEXT SENTENCE? DO WE EVER USE THIS COUPLING?
%Notice that we can couple the paths $\gamma^{\delta}_{D,a,b}$ simultaneously
%for all values of $\delta$ by using the same percolation configuration
%to generate all of them.
We can apply the results of~\cite{ab} to conclude that there exist
subsequences $\delta_k \downarrow 0$ such that the law of
%CHANGE (02-06-06) [SHOULD WE WRITE $\gamma_k$ INSTEAD OF $\gamma^{\delta_k}_k$
%HERE AND LATER?]
$\gamma^{\delta_k}_k \equiv \gamma^{\delta_k}_{D,a,b}$
(i.e., the percolation exploration path inside $D^{\delta_k}$ from $a^{\delta_k}$
to $b^{\delta_k}$) converges to some limiting law for a process $\tilde\gamma$
supported on (H\"older) continuous curves inside $D$ from $a$ to $b$.
%CHANGE (01-20-06): SENTENCE ADDED
%All the sequences $\delta_k \downarrow 0$ that appear in the rest of the
%paper are assumed to have this property.
The curves are defined up to (monotonic) reparametrizations; in the next
theorem and its proof, even where the time variable appears explicitly,
we do not specify a parametrization since it is irrelevant.
The filling $\tilde K_t$ of $\tilde\gamma[0,t]$, appearing in the next
theorem, is defined just above~(\ref{hulls}).

%Let $C' \subset D$ be a closed subset of $\overline D$ such that
%$a \notin C'$, $b \in C'$, and $D' \equiv D \setminus C'$ is a bounded Jordan
%domain whose boundary contains the counterclockwise arc $\overline {cd}$ that
%does not belong to $\partial D$ (except for its endpoints $c$ and $d$).
%Let $\tilde K_t$ denote the hull of $\tilde\gamma(t)$ up to time $t$ and
%$D_t = D \setminus \tilde K_t$,
%$T'=\inf \{ t:\tilde K_t \cap C' \neq \emptyset \}$ the first time that
%$\tilde\gamma(t)$ exits $D'$, and $\tilde\nu_{D',a,c,d}$ the distribution
%of the hull $\tilde K_{T'}$.
%
%Consider now a closed subset $C''$ of $D \setminus \tilde K_{T'}$ such that
%$\tilde\gamma(T') \notin C''$, $b \in C''$, and $D''  = D_{T'} \setminus C''$
%is a bounded Jordan domain whose boundary contains the counterclockwise arc $\overline{ef}$
%that does not belong to $\partial D_{T'}$ (except for its endpoints $e$ and $f$),
%and let $T''$ be the first time $\tilde\gamma(t)$ exits $D''$.

\begin{theorem} \label{spatial-markov}
For any subsequencial limit $\tilde\gamma$ of the percolation exploration
path $\gamma^{\delta}_{D,a,b}$ defined above, the filling $\tilde K_t$ of
$\tilde\gamma[0,t]$, as a process, satisfies the spatial Markov property.
\end{theorem}

\noindent {\bf Proof.} Let $\delta_k \downarrow 0$ be a subsequence such that
the law of $\gamma^{\delta_k}_k$ converges to some limiting law supported on
continuous curves $\tilde\gamma$ in $D$ from $a$ to $b$.
We will prove the spatial Markov property by showing that
$(\tilde K_{\tilde T_j}, \tilde\gamma(\tilde T_j))$ as defined in the proof
of Theorem~\ref{characterization} are \emph{jointly} distributed like the
corresponding $SLE_6$ hull variables, which do have the spatial Markov property.
Since $\gamma^{\delta_k}_k$ converges in distribution to $\tilde\gamma$,
we can find coupled versions of $\gamma^{\delta_k}_k$ and $\tilde\gamma$
on some probability space $(\Omega',{\cal B}',{\mathbb P}')$ such that
%on the Lebesgue probability space $(\Omega,{\cal B},\lambda)$ such that
$\gamma^{\delta_k}_k$ converges to $\tilde\gamma$ for all $\omega' \in \Omega'$;
in the rest of the proof we work with these new versions which, with a slight
abuse of notation, we denote with the same names as the original ones.

Let $\tilde f_0$ be a conformal transformation that maps $\mathbb H$ to $D$ such
that ${\tilde f_0}^{-1}(a)=0$ and ${\tilde f_0}^{-1}(b) = \infty$ and
let $\tilde T_1=\tilde T_1(\varepsilon)$ denote the first time $\tilde\gamma(t)$
hits $D \setminus \tilde G_1$, with $\tilde G_1 \equiv \tilde f_0(C(0,\varepsilon))$
and (as in the proof of Lemma~\ref{boundaries})
$C(0,\varepsilon) = \{ z: |z|<\varepsilon \} \cap {\mathbb H}$.
Define recursively $\tilde T_{j+1}$ as the first time $\tilde\gamma(t)$ hits
$D_{\tilde T_j} \setminus \tilde G_{j+1}$, with
$\tilde G_{j+1} \equiv \tilde f_{\tilde T_j}(C(0,\varepsilon))$ and
$D_{\tilde T_j} \equiv D \setminus \tilde K_{\tilde T_j}$,
where $\tilde f_{\tilde T_j}$ is a conformal map from $\mathbb H$ to
$D_{\tilde T_j}$ whose inverse maps $\tilde\gamma(\tilde T_j)$ to $0$
and $b$ to $\infty$, chosen as in the definition (before
Theorem~\ref{characterization}
above) of the spatial Markov property.
We also define $\tilde\tau_j \equiv \tilde T_{j+1} - \tilde T_j$, so that
$\tilde T_j = \tilde\tau_1+\ldots+\tilde\tau_j$, and the (discrete-time)
stochastic process $\tilde X_j \equiv (\tilde K_{\tilde T_j},\tilde\gamma(\tilde T_j))$.
%The conformal maps $\tilde f_{\tilde T_j}$ are only defined up to a scaling
%factor $\lambda_j>0$; we need to prove that, for any choice of the scaling
%factors for the maps $\tilde f_{\tilde T_j}$, $\tilde X_j$ is a Markov chain.

Analogous quantities can be defined for the trace of chordal $SLE_6$.
For clarity, they will be indicated here by the superscript $SLE_6$;
e.g., $f_{T_j}^{SLE_6}$, $K_{T_j}^{SLE_6}$, $G_j^{SLE_6}$ and $X_j^{SLE_6}$.
We choose $f_0^{SLE_6}=\tilde f_0$, so that $G_1^{SLE_6}=\tilde G_1$.
%$K_{T_1}^{SLE_6}=\tilde K_{\tilde T_1}$ and $X_1^{SLE_6}=\tilde X_1$

For each $k$, let $K^k_t$ denote the filling at time $t$ of $\gamma^{\delta_k}_k$
(see the definition of discrete filling just before Lemma~\ref{new-sub-conv}).
It follows from the Markovian character of the percolation exploration
process that, for all $k$, the filling $K^k_t$ satisfies a suitably adapted
(to the discrete setting) spatial Markov property.
(In fact, the percolation exploration path satisfies a stronger property -- roughly
speaking, that for \emph{all} times $t$ the future of the path given the filling
of the past is distributed as a percolation exploration path in the original domain
from which the filling up to time $t$ has been removed.)

Let now $f^k_0$ be a conformal transformation that maps $\mathbb H$ to
$D_k \equiv D^{\delta_k}$ such that $(f^k_0)^{-1}(a_k)=0$ and
$(f^k_0)^{-1}(b_k) = \infty$ and let $T^k_1=T^k_1(\varepsilon)$ denote the first
exit time of $\gamma^{\delta_k}_k(t)$ from $G^k_1 \equiv f_0^k(C(0,\varepsilon))$
defined as the first time that $\gamma^{\delta_k}_k$ intersects the image under
$f_0^k$ of the semi-circle $\{ z : |z| = \varepsilon \} \cap {\mathbb H}$.
Define recursively $T^k_{j+1}$ as the first exit time of $\gamma^{\delta_k}_k[T^k_j,\infty)$
from $G^k_{j+1} \equiv f^k_{T^k_j}(C(0,\varepsilon))$, where $f^k_{T^k_j}$
is a conformal map from $\mathbb H$ to $D_k \setminus K^k_{T^k_j}$ whose
inverse maps $\gamma^{\delta_k}_k(T^k_j)$ to $0$ and $b_k$ to $\infty$.
The maps $f^k_{T^k_j}$, for $j \geq 1$, are defined only up to a scaling factor.
We also define $\tau^k_{j+1} \equiv T^k_{j+1} - T^k_j$, so that
$T^k_j=\tau^k_1+\ldots+\tau^k_j$, and the (discrete-time) stochastic process
$X^k_j \equiv (K^k_{T^k_j},\gamma^{\delta_k}_k(T^k_j))$ for $j=1,2,\ldots$ .
The Markovian character of the percolation exploration process implies that,
for every $k$, $X^k_j$ is a Markov chain (in $j$).
We want to show recursively that, for any $j$, as $k \to \infty$,
$\{ X^k_1,\ldots,X^k_j \}$ converge jointly in distribution to
$\{ \tilde X_1,\ldots,\tilde X_j \}$.
By recursively applying Theorem~\ref{strong-cardy} and Lemma~\ref{hull}, we
will then be able to conclude that $\{ \tilde X_1,\ldots,\tilde X_j \}$ are
jointly equidistributed with the corresponding $SLE_6$ hull variables
(at the corresponding stopping times) $\{ X_1^{SLE_6},\ldots,X_j^{SLE_6} \}$.
Since the latter do satisfy the spatial Markov property, so will the
former, as desired.

The zeroth step consists in noticing that the convergence of $(D_k,a_k,b_k)$
to $(D,a,b)$ as $k \to \infty$ allows us to use Corollary~\ref{cor2-unif-conv}
%(Rad\`o's theorem -- see Theorem~\ref{rado-thm} of Appendix~\ref{rado} -- is
%sufficient in this case, but not in the later steps)
to select a sequence of conformal maps $f^k_0$ that
%scaling factors so that $f^k_0$ converges
converge to $f_0^{SLE_6}=\tilde f_0$ uniformly in $\overline{\mathbb H}$
as $k \to \infty$, which implies that the boundary $\partial G^k_1$
of $G^k_1=f^k_0(C(0,\varepsilon))$ converges to the boundary $\partial\tilde G_1$
of $\tilde G_1 = \tilde f_0(C(0,\varepsilon))$ in the uniform metric on continuous
curves.

Starting from there, the first step of our recursion argument is organized
as follows, where all limits and equalities are in distribution:
\begin{itemize}
\item[(i)] $K^k_{T^k_1} \to \tilde K_{\tilde T_1} = K_{T_1}^{SLE_6}$ by Lemma~\ref{new-boundaries}.
%%and Remark~\ref{remark-generalization},
%%``number of arms" percolation bounds~\cite{ksz} and Lemma~\ref{double-crossing} below,
%but also $K^k_{T^k_1} \to K_{T_1}^{SLE_6}$ by Lemma~\ref{hull}
%(with the use of Theorem~\ref{strong-cardy} and also Lemmas~\ref{new-sub-conv}
%and~\ref{new-boundaries}). CHANGE %--- see Remark~\ref{hull-convergence}).
\item[(ii)] by i), $D_k \setminus K^k_{T^k_1} \to D \setminus \tilde K_{\tilde T_1}
= D \setminus K_{T_1}^{SLE_6}$.
%and $D_k \setminus K^k_{T^k_1} \to D \setminus K_{T_1}^{SLE_6}$, so
%that $D \setminus K_{T_1}^{SLE_6} = D \setminus \tilde K_{\tilde T_1}$.
\item[(iii)] by (ii), $f_{T_1}^{SLE_6}=\tilde f_{\tilde T_1}$,
and by Corollary~\ref{cor2-unif-conv} of Appendix~\ref{rado} we can select
a sequence $f^k_{T^k_1} \to \tilde f_{\tilde T_1} = f_{T_1}^{SLE_6}$.
\item[(iv)] by (iii), $G^k_2 \to \tilde G_2 = G_2^{SLE_6}$.
%and $G^k_2 \to G_2^{SLE_6}$, so that $G_2^{SLE_6} = \tilde G_2$.
\end{itemize}

\begin{figure}[!ht]
\begin{center}
\includegraphics[width=8cm]{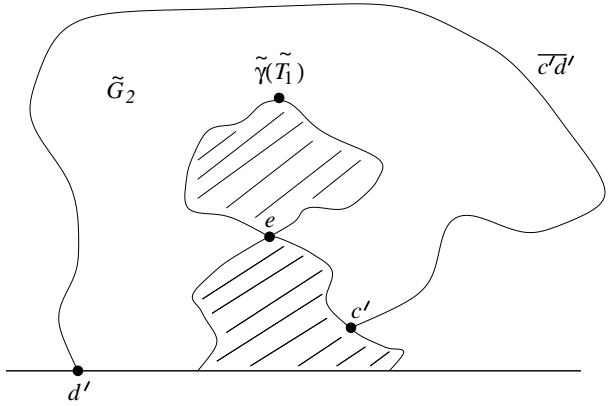}
\caption{Schematic figure representing a hull (shaded) with a cut-point $e$,
resulting in a non-Jordan, but admissible, $\tilde G_2$.}
\label{cut-point-fig}
\end{center}
\end{figure}

We remark that Lemmas~\ref{new-sub-conv} and~\ref{new-boundaries}
imply
%since Cardy's formula applies to all subsequence limits,
%the argument in the proof of Lemma~\ref{new-sub-conv} implies
that the filling
$\tilde K_{\tilde T_1}$ is a hull, and its ``envelope" is therefore composed
of two simple curves.
It follows that $D \setminus \tilde K_{\tilde T_1}$ and $\tilde G_2$ are
admissible, since the part of the boundary of either
$D \setminus \tilde K_{\tilde T_1}$ or $\tilde G_2$ that belongs to the
boundary of $\tilde K_{\tilde T_1}$ can be split up, by removing the single
point $\tilde\gamma(\tilde T_1)$, into two simple curves, while the
remaining part of the boundary of either $D \setminus \tilde K_{\tilde T_1}$
or $\tilde G_2$ is a Jordan arc whose interior does not touch the hull
$\tilde K_{\tilde T_1}$.
This allows us to use Theorem~\ref{strong-cardy} (and therefore Lemma~\ref{hull}),
Corollary~\ref{cor2-unif-conv} and Lemmas~\ref{double-crossing}
and~\ref{close-encounters}.
(Note that $D \setminus \tilde K_{\tilde T_1}$ and $\tilde G_2$
need not be Jordan because $\tilde K_{\tilde T_1}$ has cut-points
with positive probability -- see Fig.~\ref{cut-point-fig}.)

At this point, we are in the same situation as at the zeroth step,
but with $G^k_1$, $\tilde G_1$ and $G_1^{SLE_6}$ replaced by $G^k_2$,
$\tilde G_2$ and $G_2^{SLE_6}$ respectively, and we can proceed by
induction, as follows.
(As explained above, the theorem will then follow from the fact that
the $SLE_6$ hull variables do possess the spatial Markov property.)

%Before continuing with the proof, we note that we can no longer apply
%Lemma~\ref{new-boundaries} (and Remark~\ref{remark-generalization})
%because the domains in question are no longer Jordan domains.
%Instead, we will follow the proof of Lemma~\ref{new-boundaries} but
%together with Lemmas~\ref{double-crossing}-\ref{touching}, we will
%use a standard bound (see~\cite{ksz}) on the probability of six
%crossing of an annulus.
%In what follows we will show that the ``six-arm" bound~\cite{ksz}
%and Lemma~\ref{double-crossing} below imply the convergence of
%$K^k_{T^k_j}$ to $\tilde K_{\tilde T_j}$ for general $j>1$, and that
%the conditions to apply Theorem~\ref{strong-cardy}, Lemma~\ref{double-crossing}
%and Corollary~\ref{cor2-unif-conv} are always satisfied.
%(This last point boils down to showing that the domains
%$D \setminus \tilde K_{\tilde T_j}$ and $\tilde G_j$ are
%admissible for all $j$.)
%As explained in Remark~\ref{remark-generalization}, a suitable version
%of Lemma~\ref{new-boundaries} implies the joint convergence in
%distribution of $(K^k_{T^k_1},\gamma^{\delta_k}_k(T^k_1))$ to
%$(\tilde K_{\tilde T_1},\tilde\gamma(\tilde T_1))$, which concludes
%the first step of the argument.

The next step consists in proving that
$((K^k_{T^k_1},\gamma^{\delta_k}_k(T^k_1)),(K^k_{T^k_2},\gamma^{\delta_k}_k(T^k_2)))$
converges in distribution to
$((\tilde K_{\tilde T_1},\tilde\gamma(\tilde T_1)),(\tilde K_{\tilde T_2},\tilde\gamma(\tilde T_2)))$.
Since we have already proved the convergence of $(K^k_{T^k_1},\gamma^{\delta_k}_k(T^k_1))$ to
$(\tilde K_{\tilde T_1},\tilde\gamma(\tilde T_1))$, we claim that all we really need to prove
is the convergence of $(K^k_{T^k_2} \setminus K^k_{T^k_1},\gamma^{\delta_k}_k(T^k_2))$
to $(\tilde K_{\tilde T_2} \setminus \tilde K_{\tilde T_1},\tilde\gamma(\tilde T_2))$.
To do this, notice that $K^k_{T^k_2} \setminus K^k_{T^k_1}$ is distributed like
the filling of a percolation exploration path inside $D_k \setminus K^k_{T^k_1}$.
Besides, the convergence in distribution of $(K^k_{T^k_1},\gamma^{\delta_k}_k(T^k_1))$
to $(\tilde K_{\tilde T_1},\tilde\gamma(\tilde T_1))$ implies that we can find versions
of $(\gamma^{\delta_k}_k,K^k_{T^k_1})$ and $(\tilde\gamma,\tilde K_{\tilde T_1})$
on some probability space $(\Omega',{\cal B}',{\mathbb P}')$ such that
$\gamma^{\delta_k}_k(\omega')$ converges to $\tilde\gamma(\omega')$ and
$(K^k_{T^k_1},\gamma^{\delta_k}_k(T^k_1))$ converges to
$(\tilde K_{\tilde T_1},\tilde\gamma(\tilde T_1))$ for all $\omega' \in \Omega'$.
These two observations imply that, if we work with the coupled versions of
$(\gamma^{\delta_k}_k,K^k_{T^k_1})$ and $(\tilde\gamma,\tilde K_{\tilde T_1})$,
we are in the same situation as before, but with $D_k$ (resp., $D$) replaced by
$D_k \setminus K^k_{T^k_1}$ (resp., $D \setminus \tilde K_{\tilde T_1}$) and $a_k$
(resp., $a$) by $\gamma^{\delta_k}_k(T^k_1)$ (resp., $\tilde\gamma(\tilde T_1)$).

\begin{figure}[!ht]
\begin{center}
\includegraphics[width=8cm]{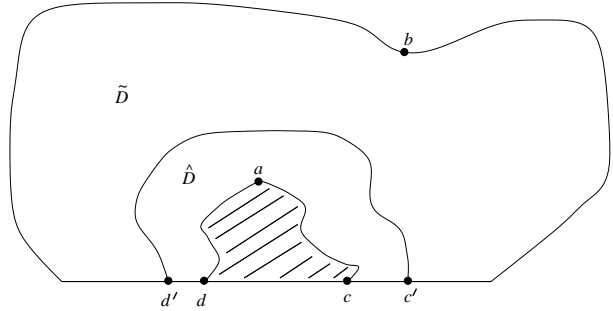}
\caption{Schematic figure representing the situation described in
Lemma~\ref{double-crossing}. Note that the shaded
region is not part of the domains $\tilde D$ and $\hat D$ and that
$\tilde D$ contains $\hat D$.
%An exploration path started at $a$ and aimed at $b$ is defined inside
%(a $\delta$-approximation of) $\tilde D$, and stopped when it first exits $\hat D$.
In the application to Theorem~\ref{spatial-markov}, the shaded region
represents the hull of the past (which may have cut-points, as in
Figure~\ref{cut-point-fig}) and the counterclockwise arc $\overline{c'd'}$
of $\partial \hat D$ represents the conformal image of a semicircle.
Notice that $c$ and $c'$ (resp., $d$ and $d'$) coincide if the right
(resp., left) endpoint of the conformal image of the semicircle lies
on the boundary of the shaded region (see, e.g., Figure~\ref{cut-point-fig}).
In the application to Lemmas~\ref{new-sub-conv}-\ref{new-boundaries},
$D=\tilde D$, $G=\hat D$ and $\partial^*G$ corresponds to the
counterclockwise arc $\overline{c'd'}$ of $\partial \hat D$.}
\label{fig-lemma6-1}
\end{center}
\end{figure}

Then, the conclusion that
$(K^k_{T^k_2} \setminus K^k_{T^k_1},\gamma^{\delta_k}_k(T^k_2))$
%$((K^k_{T^k_1},\gamma^{\delta_k}_k(T^k_1)),(K^k_{T^k_2},\gamma^{\delta_k}_k(T^k_2)))$
converges in distribution to
$(\tilde K_{\tilde T_2} \setminus \tilde K_{\tilde T_1},\tilde\gamma(\tilde T_2))$
%$((\tilde K_{\tilde T_1},\tilde\gamma(\tilde T_1)),(\tilde K_{\tilde T_2},\tilde\gamma(\tilde T_2)))$
follows, as before, by arguments like those used for Lemma~\ref{new-boundaries}
%and Remark~\ref{remark-generalization}
-- i.e., by using a standard bound on the probability of six disjoint
monochromatic crossings~\cite{ksz} and Lemmas~\ref{double-crossing}
and~\ref{close-encounters}, as we now explain.
$G^k_2$ (resp., $\tilde G_2$) is a domain admissible with respect to
$(\gamma^{\delta_k}_k(T^k_1),c'_k,d'_k)$ (resp., $(\tilde\gamma(\tilde T_1),c',d')$),
where $c'_k$ and $d'_k$ (resp., $c'$ and $d'$) are the unique points where the image
of $\overline{\partial {\mathbb D} \cap {\mathbb H}}$ under $f^k_{T^1_k}$ (resp.,
$\tilde f_{\tilde T_1}$) meets either the envelope of $K^k_{T^k_1}$ (resp.,
$\tilde K_{\tilde T_1}$) or $\partial D_k$ (resp., $\partial D$) -- see Figure~\ref{cut-point-fig}.
The envelope of the filling $K^k_{T^k_1}$ (resp., $\tilde K_{\tilde T_1}$) is
part of the boundary of the domain explored by $\gamma^{\delta_k}_k[T^k_1,T^k_2]$
(resp., $\tilde\gamma[\tilde T_1,\tilde T_2]$).
Close encounters of $\gamma^{\delta_k}_k[T^k_1,T^k_2]$ (resp.,
$\tilde\gamma[\tilde T_1,\tilde T_2]$) with this part of the boundary can be dealt
with using again a standard bound~\cite{ksz} on the probability of six disjoint
monochromatic crossings, as explained in the proofs of Lemmas~\ref{sub-conv}-\ref{boundaries}
(see also~\cite{ada}).
If part of the boundary of $G^k_2$ (resp., $\tilde G_2$) coincides with part of
$\partial D_k$ (resp., $\partial D$), we can use Lemmas~\ref{double-crossing}
and~\ref{close-encounters} as in Lemmas~\ref{new-sub-conv} and~\ref{new-boundaries},
to obtain the same conclusions.
Notice that Lemmas~\ref{double-crossing} and~\ref{close-encounters} are adapted to
the situation we encounter here, with the boundary of the exploration domain
divided in three parts (corresponding here to the envelope of the past filling,
part of the boundary of the original domain, and the semi-circle conformally
mapped from the upper half-plane -- see Figures~\ref{cut-point-fig} and~\ref{fig-lemma6-1}).

%The only slight difference with the previous case is that now the domain explored by %$\gamma^{\delta_k}_k[T^k_1,T^k_2]$ (resp., $\tilde\gamma[\tilde T_1,\tilde T_2]$) is
%admissible, but not necessarily Jordan, and the envelope of the filling $K^k_{T^k_1}$
%(resp., $\tilde K_{\tilde T_1}$) is part of the boundary of the domain explored by
%$\gamma^{\delta_k}_k[T^k_1,T^k_2]$ (resp., $\tilde\gamma[\tilde T_1,\tilde T_2]$),
%However, this is not a problem since close encounters with the past filling
%can be dealt with using again a standard bound~\cite{ksz} on the probability
%of six disjoint monochromatic crossings, as explained in the proofs of
%Lemmas~\ref{sub-conv}-\ref{boundaries} (see also~\cite{ada}).

%In order to get claim (C'), in places where the exploration path comes close
%to the boundary of the past filling we can use the bound on the probability of
%six crossings of an annulus in the plane (like, for instance, in the proofs
%of Lemmas~\ref{sub-conv} and~\ref{boundaries}), while in places where it
%comes close to the remaining portion of the boundary (i.e., $\partial D$
%or the Jordan arc $\overline{cd}$ in Figure~\ref{cut-point-fig}) we can
%use Lemmas~\ref{double-crossing}-\ref{touching}.

We can now iterate
the above arguments $j$ times, for any $j>1$.
It is in fact easy to see by induction that the domains
$D \setminus \tilde K_{\tilde T_j}$ and $\tilde G_j$
that appear in the successive steps are admissible for all $j$.
Therefore we can keep using
Theorem~\ref{strong-cardy} (and Lemma~\ref{hull}),
Corollary~\ref{cor2-unif-conv} and Lemmas~\ref{double-crossing}
and~\ref{close-encounters}.
If we keep track at each step of the previous ones, this provides
the \emph{joint} convergence of all the curves and fillings involved
at each step, and concludes the proof of Theorem~\ref{spatial-markov}. \fbox{}

\begin{remark} \label{explanation}
The key technical problem in proving Theorem~\ref{spatial-markov} is
showing that for the exploration path $\gamma^{\delta_k}_k$, one can
interchange the limit ${\delta}_k \to 0$ with the process of filling.
This requires showing two things about the exploration
path: (1) the return of a
(macroscopic) segment of the path close to an earlier segment
(and away from $\partial D_k$) without nearby (microscopic)
touching does not occur (probably), and (2) the close approach of a
(macroscopic) segment to $\partial D_k$ without nearby (microscopic) touching
either of $\partial D_k$ itself or else of another segment that touches
$\partial D_k$ does not occur (probably). In the related proof of
Lemma~\ref{boundaries} where $D$ was just the unit disk, these were
controlled by known estimates on probabilities of six-arm events in
the full plane for (1) and of three-arm events in the half-plane for (2).
When $D$ is not necessarily convex, as in Theorem~\ref{spatial-markov},
the three-arm event argument for (2) appears to break down. Our replacement
is the use of Lemmas~\ref{double-crossing}-\ref{mushroom}.
Basically, these control (2) by
%known estimates on six-arm events in the
%full plane, but applied to points on $\partial D$, combined with
a novel
argument about ``mushroom events'' on $\partial D$ (see Lemma~\ref{mushroom}),
which is based on continuity of Cardy's formula with respect to changes
in $\partial D$.
\end{remark}

%In the next lemma, used in the proof of Theorem~\ref{spatial-markov},
%in order to keep the notation simple, when defining the exploration path,
%we will not distinguish between the domains $\tilde D_k$ and their
%$\delta$-approximations $\tilde D^{\delta}_k$.
%This abuse of notation is innocuous since $\tilde D_k^{\delta} \to \tilde D_k$
%as $\delta \to 0$ (i.e., $\partial \tilde D^{\delta}_k \to \partial \tilde D_k$
%in the uniform metric on continuous curves).
The situation described in the next lemma is depicted in Figure~\ref{fig-lemma6-1}
and corresponds to those in the proof of Theorem~\ref{spatial-markov} and in
Lemmas~\ref{new-sub-conv}-\ref{new-boundaries}.
In the next lemmas and their proofs, when we write that a percolation exploration
path in $\delta{\cal H}$ touches itself we mean that it gets to distance $\delta$
of itself, just one hexagon away.

\begin{lemma} \label{double-crossing}
Let $\{ (\hat D_k,a_k,c_k,c'_k,d'_k,d_k) \}$ be a sequence of Jordan
domains with five points (not necessarily all distinct) on their boundaries
in counterclockwise order.
Assume that $\hat D_k \subset \tilde D_k$, where $\tilde D_k$ is a Jordan
set from $\delta_k \cal H$, that $a_k,c_k,c'_k,d'_k,d_k \in \partial \tilde D_k$,
that the counterclockwise arcs $\overline{d'_k c'_k}$ of $\partial \hat D_k$ and
$\partial \tilde D_k$ coincide, and that $a_k$ is an e-vertex of $\partial \tilde D_k$.
Consider a second e-vertex $b_k \in \partial \tilde D_k, b_k \notin \partial \hat D_k$,
and denote by $\gamma^{\delta_k}_k$ the percolation exploration path in $\tilde D_k$
started at $a_k$, aimed at $b_k$, and stopped when it first hits the counterclockwise arc
$J'_k=\overline{c'_k d'_k} \subset J_k=\overline{c_k d_k}$ of $\partial \hat D_k$.
Assume that, as $k \to \infty$, $\delta_k \downarrow 0$ and
$(\tilde D_k,a_k,c_k,d_k) \to (\tilde D,a,c,d)$,
$(\hat D_k,a_k,c_k,d_k) \to (\hat D,a,c,d)$, where $\tilde D$ and $\hat D$ are domains
admissible with respect to $(a,c,d)$.
Assume also that $J'_k$ converges in the metric~(\ref{distance}) to the counterclockwise
arc $J' \equiv \overline{c'd'}$ of $\partial \hat D$, a subset of the counterclockwise arc
$J \equiv \overline{cd}$ of $\partial \hat D$, and that
$b_k \to b \in \partial \tilde D$, $b \notin \partial \hat D$.

Let ${\cal E}_k(J_k;\varepsilon,\varepsilon') =
\{ \bigcup_{v \in J_k \setminus J'_k} {\cal B}_k(v;\varepsilon,\varepsilon') \}
\cup \{ \bigcup_{v \in J'_k} {\cal A}_k(v;\varepsilon,\varepsilon') \}$,
where ${\cal A}_k(v;\varepsilon,\varepsilon')$ is the event that
$\gamma^{\delta_k}_k$ contains a segment that stays within $B(v,\varepsilon)$ and
has a double crossing of the annulus $B(v,\varepsilon) \setminus B(v,\varepsilon')$
without that segment touching $\partial \hat D_k$, and
${\cal B}_k(v;\varepsilon,\varepsilon')$ is the event that
$\gamma^{\delta_k}_k$ enters $B(v,\varepsilon')$, but is stopped outside
$B(v,\varepsilon)$ and does not touch $\partial \hat D_k \cap B(v,\varepsilon)$.
Then, for any $\varepsilon>0$, %$0 < \varepsilon < \min \{ |a-c|, |a-d| \}$,
\begin{equation} \label{double}
\lim_{\varepsilon' \to 0} \, \limsup_{k \to \infty}
\, {\mathbb P}({\cal E}_k(J_k;\varepsilon,\varepsilon')) = 0.
\end{equation}
\end{lemma}

% Before giving the proof of Lemma~\ref{double-crossing}, we present a new
Our next result is a lemma in which the part of Lemma~\ref{double-crossing}
concerning $J \setminus J'$ can be strengthened to conclude that if $\tilde\gamma$
touches $J \setminus J'$ somewhere, then for $k$ large enough, $\gamma^{\delta_k}_k$
touches either $J_k \setminus J'_k$ or its past hull ``nearby."
Lemmas~\ref{double-crossing} and~\ref{close-encounters} are used in the proof
of Theorem~\ref{spatial-markov} to show that in the limit there is no discrepancy
between the hull generated by $\tilde\gamma$ and the limit as $k \to \infty$ of the
hull generated by $\gamma^{\delta_k}_k$.

\begin{lemma} \label{close-encounters}
With the notation of Lemma~\ref{double-crossing}, let
${\cal C}_k(v;\varepsilon,\varepsilon')$ be the event that $\gamma^{\delta_k}_k$
contains a segment that stays within $B(v,\varepsilon)$ and has a double crossing of
the annulus $B(v,\varepsilon) \setminus B(v,\varepsilon')$ without that segment touching
either $\partial \hat D_k$ or any other segment of $\gamma^{\delta_k}_k$ that stays within
$B(v,\varepsilon)$ and touches $\partial \hat D_k$.
Then, for any $\varepsilon>0$,
\begin{equation} \label{close}
\lim_{\varepsilon' \to 0} \, \limsup_{k \to \infty} \,
{\mathbb P}(\bigcup_{v \in J_k \setminus J'_k} {\cal C}_k(v;\varepsilon,\varepsilon')) = 0.
\end{equation}
\end{lemma}

The proofs of Lemmas~\ref{double-crossing} and~\ref{close-encounters}
are partly based on relating the failure of~(\ref{double}) or~(\ref{close})
to the occurrence with strictly positive probability of certain continuum
limit ``mushroom" events (see Lemma~\ref{mushroom}) that we will show must
have zero probability because otherwise there would be a contradiction
to Lemma~\ref{equal}, which itself is a consequence of the continuity
of Cardy's formula with respect to the domain boundary
(see Lemma~\ref{cont-cardy} of Appendix~\ref{rado}).
In both of the next two lemmas, we denote by $\mu$ any subsequence
limit as $\delta=\delta_k \to 0$ of the probability measures for the
collection of all colored (blue and yellow) $\cal T$-paths on \emph{all}
of ${\mathbb R}^2$, in the Aizenman-Burchard sense (see Remark~\ref{ab}).
We recall that in our notation, $D$ represents an open domain and
$\overline{z_1 z_2}$, $\overline{z_3 z_4}$ represent closed segments
of its boundary.
In Lemma~\ref{equal} below, we restrict attention to a Jordan domain $D$
since that case suffices for the use of Lemma~\ref{equal} in the proofs
of Lemmas~\ref{double-crossing} and~\ref{close-encounters}.

\begin{lemma} \label{equal}
For $(D,z_1,z_2,z_3,z_4)$, with $D$ a Jordan domain, consider the following
crossing events, ${\cal C}^*_i = {\cal C}^*_i(D,z_1,z_2,z_3,z_4)$, where $*$
denotes either blue or yellow, a $*$ path denotes a segment of a $*$ curve,
and $i=1,2,3$:
\begin{equation} \nonumber
{\cal C}^*_1 = \{ \exists \text{ a $*$ path in the closure } \overline D
\text{ from } \overline{z_1 z_2} \text{ to } \overline{z_3 z_4} \},
\end{equation}
\begin{equation} \nonumber
{\cal C}^*_2 = \{ \exists \text{ a $*$ path in } D \text{ from the interior of }
\overline{z_1 z_2} \text{ to the interior of } \overline{z_3 z_4} \},
\end{equation}
\begin{equation} \nonumber
{\cal C}^*_3 = \{ \exists \text{ a $*$ path starting and ending \emph{outside} }
\overline D \text{ whose restriction to } D \text{ is as in } {\cal C}^*_2 \}.
\end{equation}
Then $\mu({\cal C}^*_1)=\mu({\cal C}^*_2)=\mu({\cal C}^*_3)=\Phi_D(z_1,z_2;z_3,z_4)$.
\end{lemma}

\noindent {\bf Proof.}
%The proof is similar to that of Theorem~\ref{strong-cardy},
%but easier because $D$ is here a Jordan domain.
%Indeed, it is enough to construct a new Jordan domain $\tilde D(\varepsilon)$
%(with appropriately selected points $\tilde z_1(\varepsilon),\tilde z_2(\varepsilon),
%\tilde z_3(\varepsilon),\tilde z_4(\varepsilon)$  on the boundary and corresponding
%events $\tilde{\cal C}^*_i$) such that the occurrence of $\tilde{\cal C}^*_1$
%in $\tilde D(\varepsilon)$ implies the occurrence of ${\cal C}^*_3$ in $D$ and
%with $(\tilde D(\varepsilon),\tilde z_1(\varepsilon),\tilde z_2(\varepsilon),
%\tilde z_3(\varepsilon),\tilde z_4(\varepsilon)) \to (D,z_1,z_2,z_3,z_4)$ as
%$\varepsilon \to 0$.
%The continuity of Cardy's formula (Lemma~\ref{cont-cardy} in Appendix~\ref{rado})
%does the rest. \fbox{}
Recall that %because of the choice of topology, the
convergence to $\mu$
implies a coupling of the lattice and continuum processes
on some $(\Omega',{\cal B}',{\mathbb P}')$ such that the distance between the set
of $\cal T$-paths and the set of continuum paths tends to zero as $\delta_k \to 0$
for all $\omega' \in \Omega'$ (see, e.g., Corollary~1 of~\cite{billingsley1}).

We will construct for each small $\varepsilon>0$,
two domains with boundary points, denoted by
$(\tilde D, \tilde z_1, \tilde z_2, \tilde z_3, \tilde z_4)$
and $(\hat D, \hat z_1, \hat z_2, \hat z_3, \hat z_4)$, approximating $(D,z_1,z_2,z_3,z_4)$
in such a way that
$\tilde\Phi_{\varepsilon} \equiv \Phi_{\tilde D}(\tilde z_1, \tilde z_2; \tilde z_3, \tilde z_4)
\stackrel{\varepsilon \to 0}{\longrightarrow} \Phi \equiv \Phi_D(z_1,z_2;z_3,z_4)$ and the same
for $\hat\Phi_{\varepsilon}$, and with the property that
$\tilde\Phi_{\varepsilon} \leq \mu({\cal C}^*_i) \leq \hat\Phi_{\varepsilon}$
for $i=1,2,3$.
This will yield the desired result.
The construction of the approximating domains uses fairly straightforward conformal
mapping arguments.
We provide details for $\tilde D$; the construction of $\hat D$ is analogous.

\begin{figure}[!ht]
\begin{center}
\includegraphics[width=8cm]{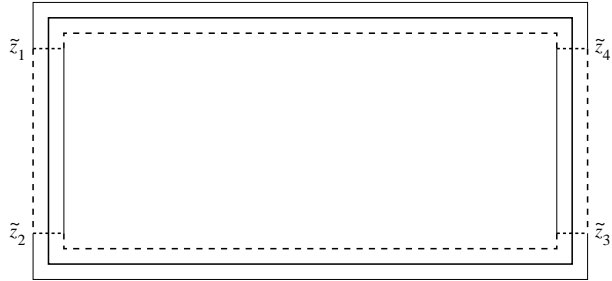}
\caption{Here, the middle rectangle is the domain $D$ while the
boundaries of the two rectangles outside and inside $D$ are
$\overline E$ and $\underline E$ respectively.
$\tilde D$ is the domain with a dashed boundary.}
\label{fig1-lemma7-3}
\end{center}
\end{figure}

To construct $\tilde D$ we will need continuous simple loops, $\underline E(D,\varepsilon)$
and $\overline E(D,\varepsilon)$, that are inner and outer approximations to $\partial D$
in the sense that $\underline E$ is surrounded by $\partial D$ which is surrounded by
$\overline E$ with
\begin{equation}
\text{d}(\partial D, \underline E) \leq \varepsilon, \,\,\,
\text{d}(\partial D, \overline E) \leq \varepsilon.
\end{equation}
We will also need four simple curves $\{ \partial_1,\partial_2,\partial_3,\partial_4 \}$
in the interior of the (topological) annulus between $\underline E$ and $\overline E$
and connecting their endpoints
$\{ (\underline z_1,\overline z_1),(\underline z_2,\overline z_2),(\underline z_3,\overline z_3),
(\underline z_4,\overline z_4) \}$
on $\underline E$ and $\overline E$ with each touching $\partial D$ at exactly one point
which is either in the interior of the counterclockwise segment $\overline{z_1 z_2}$
(for $\partial_1$ and $\partial_2$) or else the counterclockwise segment $\overline{z_3 z_4}$
(for $\partial_3$ and $\partial_4$).
Furthermore each of these connecting curves is close to its corresponding point $z_1,z_2,z_3$,
or $z_4$; i.e., $\text{d}(\partial_1,z_1) \leq \varepsilon$, etc.
In the special case where $D$ is a rectangle, the construction of
$\overline E, \underline E, \tilde D$ (and $\hat D$) is easily
done --- see Figure~\ref{fig1-lemma7-3}.

Returning to a general Jordan domain $D$, we will take
$\tilde z_1 = \overline z_1$, $\tilde z_2 = \overline z_2$,
$\tilde z_3 = \overline z_3$, and $\tilde z_4 = \overline z_4$ with $\partial \tilde D$
the concatenation of: $\partial_1$ from $\underline z_1$ to $\overline z_1$, the portion
of $\overline E$ from $\overline z_1$ to $\overline z_2$ counterclockwise, $\partial_2$
from $\overline z_2$ to $\underline z_2$, the portion of $\underline E$ from $\underline z_2$
to $\underline z_3$ counterclockwise, $\partial_3$ from $\underline z_3$ to $\overline z_3$,
the portion of $\overline E$ from $\overline z_3$ to $\overline z_4$ counterclockwise,
$\partial_4$ from $\overline z_4$ to $\underline z_4$, and the portion of $\underline E$
from $\underline z_4$ to $\underline z_1$ counterclockwise.
It is important that (for fixed $\varepsilon$ and $\tilde D$) there is a strictly
positive minimal distance between $\underline E \cup \overline E$ and $\partial D$,
and between $\partial\tilde D$ and the union of the two counterclockwise segments
$\overline{z_2 z_3}$ and $\overline{z_4 z_1}$ of $\partial D$
(see Figure~\ref{fig1-lemma7-3}).
These features will guarantee that for fixed $\varepsilon$, once $k$ is large enough,
a $\delta_k$-lattice crossing within $\tilde D$ that corresponds to the crossing
event whose (limiting) probability is
$\Phi_{\tilde D}(\tilde z_1,\tilde z_2;\tilde z_3,\tilde z_4)$
must have a subsegment that satisfies the conditions defining ${\cal C}^*_i$
(in $D$).

\begin{figure}[!ht]
\begin{center}
\includegraphics[width=5cm]{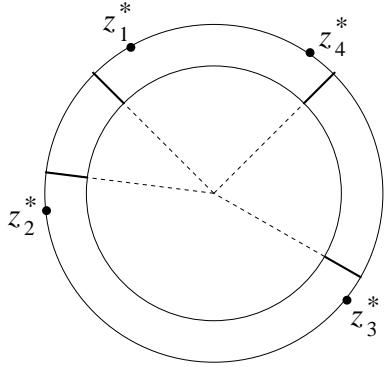}
\caption{The figure shows $\partial^*_1, \partial^*_2, \partial^*_3, \partial^*_4$
represented as heavy segments between the unit circle and the circle of radius
$1-\varepsilon'$ near $z_1^*,z_2^*,z_3^*,z_4^*$.}
\label{fig2-lemma7-3}
\end{center}
\end{figure}

We construct the parts of $\partial\tilde D$ that are inside and outside
$\partial D$ separately and then paste them together (with some care to
make sure that they ``match up").
Let $\phi$ be the conformal map from $\mathbb D$ onto $D$ with $\phi(0)=0$ and $\phi'(0)>0$,
and consider the image $\phi(\partial{\mathbb D}_{1-\varepsilon'})$ of the circle
$\partial{\mathbb D}_{1-\varepsilon'} = \{ z : |z| = 1 - \varepsilon' \}$ under $\phi$
and the inverse images $z_1^*,z_2^*,z_3^*,z_4^*$ under $\phi^{-1}$ of $z_1,z_2,z_3,z_4$.
Let $\partial^*_1(\varepsilon',\varepsilon_1)$ be the straight line between
$e^{-i \varepsilon_1} z_1$ on the unit circle $\partial {\mathbb D}$, and
$(1-\varepsilon') e^{-i \varepsilon_1} z_1$ on the circle
$\partial {\mathbb D}_{1-\varepsilon'}$, and define $\partial^*_2,\partial^*_3$,
and $\partial^*_4$ similarly,  but using \emph{clockwise} rotations by
$e^{+i \varepsilon_2}$ and $e^{+i \varepsilon_4}$ for $z_2$ and $z_4$ (see
Figure~\ref{fig2-lemma7-3}).
$\phi(\partial {\mathbb D}_{1-\varepsilon'})$ is a candidate for $\underline E(D,\varepsilon)$
and $\phi(\partial^*_{\sharp}(\varepsilon',\varepsilon_{\sharp}))$ is a candidate for half of
$\partial_{\sharp}$ (where ${\sharp}=1$ or $2$ or $3$ or $4$), so we must choose
$\varepsilon'$ and the $\varepsilon_{\sharp}$'s small enough so that
$\text{d}(\partial D, \phi({\mathbb D}_{1-\varepsilon'})) \leq \varepsilon$,
$\text{d}(\phi(\partial^*_1(\varepsilon',\varepsilon_1)),z_1) \leq \varepsilon$, etc.
We do an analogous construction using a conformal mapping from the exterior
of $\mathbb D$ onto the exterior of $D$, to obtain candidates for
$\overline E(D,\varepsilon)$ and for the exterior half of the $\partial_{\sharp}$'s.
Finally we use the freedom to choose the exterior values for $\varepsilon'$
and the $\varepsilon_{\sharp}$'s differently from the interior ones to make sure
that the interior and exterior halves of the $\varepsilon_{\sharp}$'s match up.

It should be clear that for a given approximation $\partial\tilde D$ of $\partial D$
constructed as described above there is a strictly positive $\tilde\varepsilon$ such
that the distance between $\partial D$ and the portions of $\partial\tilde D$ that
belong to $\underline E$ and $\overline E$ is not smaller than $\tilde\varepsilon$,
and the distance between the union $\overline{z_2 z_3} \cup \overline{z_4 z_1}$ of
two counterclockwise segments of $\partial D$ and
$\partial_1 \cup \partial_2 \cup \partial_3 \cup \partial_4$ is also not smaller
than $\tilde\varepsilon$.
This implies (see the definition of blue crossing just before Theorem~\ref{cardy-smirnov})
that for $k$ large enough, any (lattice) blue path crossing inside $\tilde D$ from the
counterclockwise segment $\overline{\overline z_1 \overline z_2}$ of $\partial \tilde D$
to the counterclockwise segment $\overline{\overline z_3 \overline z_4}$ of
$\partial \tilde D$ (and thanks to the coupling, any limiting continuum counterpart)
must (with high probability) have a subpath that satisfies the conditions of ${\cal C}^*_i$.
Thus, for $i=1,2,3$,
\begin{equation} \label{bound}
\mu({\cal C}^*_i) \geq
\lim_{k \to \infty} \tilde\Phi^{\delta_k}_{\tilde D(\varepsilon)} = \tilde\Phi_{\varepsilon}
\end{equation}
as desired (the equality uses Theorem~\ref{cardy-smirnov} for the Jordan domain
$\tilde D(\varepsilon)$).

We now note that as $\varepsilon \to 0$,
$(\tilde D, \tilde z_1, \tilde z_2, \tilde z_3, \tilde z_4) \to (D,z_1,z_2,z_3,z_4)$.
This allows us to use the continuity of Cardy's formula (Lemma~\ref{cont-cardy} in
Appendix~\ref{rado}) to obtain
\begin{equation} \label{cardy2}
\lim_{\varepsilon \to 0} \tilde\Phi_{\varepsilon} = \Phi.
\end{equation}
From this and~(\ref{bound}) it follows that for $i=1,2,3$,
\begin{equation} \label{upper}
\mu({\cal C}^*_i) \geq \Phi.
\end{equation}

The remaining part of the proof involves defining a domain $\hat D$
analogous to $\tilde D$ (i.e., such that
$\lim_{\varepsilon \to 0} \lim_{k \to \infty}
\hat \Phi^{\delta_k}_{\hat D(\varepsilon)} = \Phi$)
but with the property that any blue $\delta_k$-lattice path crossing inside
$D$ from the counterclockwise segment $\overline{z_1 z_2}$ to the counterclockwise
segment $\overline{z_3 z_4}$ of $\partial D$ must have a subpath in $\hat D$ that
crosses between the counterclockwise segment $\overline{\overline z_1 \overline z_2}$
and the counterclockwise segment $\overline{\overline z_3 \overline z_4}$ of
$\partial\hat D$.
(The details of the construction of $\hat D$ are analogous to those of
$\tilde D$; we leave them to the reader.)
Then, for $i=1,2,3$,
\begin{equation} \label{lower}
\mu({\cal C}^*_i) \leq \Phi,
\end{equation}
which, combined with~(\ref{upper}), implies $\mu({\cal C}^*_i) = \Phi$
and concludes the proof. \fbox{} \\

\begin{figure}[!ht]
\begin{center}
\includegraphics[width=7cm]{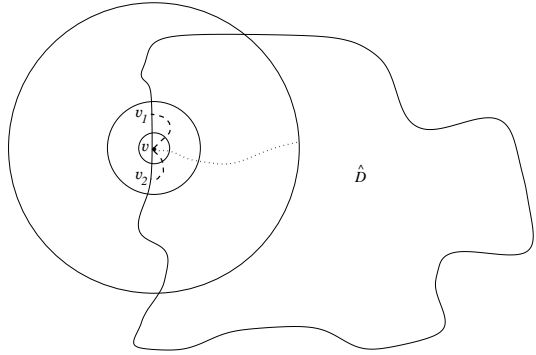}
\caption{A yellow ``mushroom" event. The dashed path is blue
and the dotted path is yellow.
The three circles centered at $v$ in the figure have radii
$\varepsilon/8$, $\varepsilon/3$, and $\varepsilon$ respectively.}
\label{fig1-lemma7-4}
\end{center}
\end{figure}

\begin{lemma} \label{mushroom}
For $(\hat D,a,c,d)$ as in Lemma~\ref{double-crossing},
$v \in J \equiv \overline{cd}$, and $\varepsilon>0$, we define
$U^{yellow}(\hat D,\varepsilon,v)$, the yellow ``mushroom" event (at $v$),
to be the event that there is a yellow path in the closure of $\hat D$
from $v$ to $\partial B(v,\varepsilon)$ and a blue path in the closure
of $\hat D$, between some pair of distinct points $v_1,v_2$ in
$\partial \hat D \cap \{ B(v,\varepsilon/3) \setminus B(v,\varepsilon/8) \}$,
that passes through $v$ and such that this blue path is \emph{between}
$\partial \hat D$ and the yellow path (see Figure~\ref{fig1-lemma7-4}).
We similarly define $U^{blue}(\hat D,\varepsilon,v)$ with the colors
interchanged and
$U^*(\hat D,\varepsilon,J) = \cup_{v \in J} U^*(\hat D,\varepsilon,v)$
where $*$ denotes blue or yellow.
Then for any deterministic domain $\hat D$ and any
$0 < \varepsilon < \min \{ |a-c|, |a-d| \}$, $\mu(U^*(\hat D,\varepsilon,J))=0$.
\end{lemma}

%STATEMENT OF LEMMA CHANGED
%\begin{lemma} \label{mushroom}
%Let $D$ be a bounded simply connected domain whose boundary $\partial D$
%is the finite union of Jordan arcs.
%For $\varepsilon>0$ smaller than the diameter of $D$, denote by
%$J(D,\varepsilon)$ the union of $v \in \partial D$ such that
%$\partial D \cap B(v,\varepsilon)$ is a Jordan arc.
%Let ....
%\end{lemma}

\begin{figure}[!ht]
\begin{center}
\includegraphics[width=7cm]{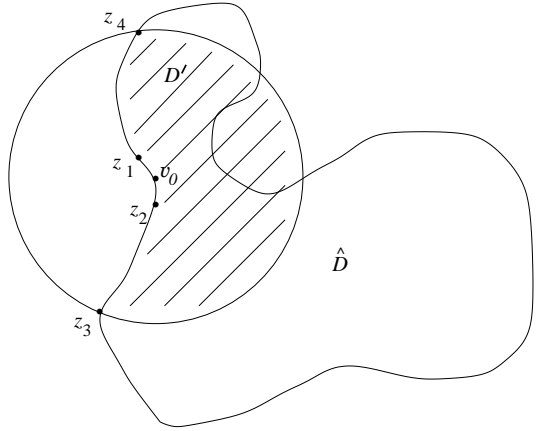}
\caption{Construction of the domain $D'$ (shaded)
used in the proof of Lemma~\ref{mushroom}.}
\label{fig2-lemma7-4}
\end{center}
\end{figure}

\noindent {\bf Proof.}
If $\mu(U^*(\hat D,\varepsilon,J))>0$ for some
$\varepsilon>0$, then there is some segment $\overline{z_1 z_2} \subset J$
of $\partial \hat D$ of diameter not larger than $\varepsilon/10$ such that
\begin{equation}
\mu(\cup_{v \in \overline{z_1 z_2}} U^*(\hat D,\varepsilon,v)) > 0.
\end{equation}
Choose any point $v_0 \in \overline{z_1 z_2}$ and consider the new domain $D'$
whose boundary consists of the correctly chosen (as we explain below) segment
of the circle $\partial B(v_0,\varepsilon/2)$ between the two points $z_3,z_4$
where $\partial \hat D$ first hits $\partial B(v_0,\varepsilon/2)$ on
either side of $v_0$, together with the segment from $z_4$ to $z_3$ of
$\partial \hat D$ (see Figure~\ref{fig2-lemma7-4}).
The correct circle segment between $z_3$ and $z_4$ is the (counter)clockwise
one if $v_0$
%is between
%CHANGE? I'M NOT SURE WHY ORDER OF $z_3$ and $z_4$ WAS CHANGED NEXT
%SINCE IT DOESN'T SEEM TO CHANGE THE MEANING AND IT LOOKS
%FUNNY TO HAVE TWO DIFF ORDERS UNLESS THERE'S A REASON --
%FEDERICO: THE ORDER MAKES A DIFFERENCE BECAUSE IT IDENTIFIES ONE OR THE OTHER
%OF TWO ARCS OF THE BOUNDARY; I THINK THAT THE ORDER SHOULD BE...
comes after $z_4$ and before $z_3$ along $\partial \hat D$ when
$\partial \hat D$ is oriented (counter)clockwise.
It is also not hard to see that since $\varepsilon < \min \{ |a-c|, |a-d| \}$,
$D'$ is a Jordan domain, so that Lemma~\ref{equal} can be applied.
In the new domain $D'$, $\overline{z_1 z_2}$ is the same curve segment as
it was in the old domain $\hat D$, but $\overline{z_3 z_4}$ is now a segment
of the circle $\partial B(v_0,\varepsilon/2)$.
It should be clear that
\begin{equation}
\cup_{v \in \overline{z_1 z_2}} U^*(\hat D,\varepsilon,v) \subset
{\cal C}^*_1(D',z_1,z_2,z_3,z_4) \setminus {\cal C}^*_3(D',z_1,z_2,z_3,z_4)
\end{equation}
which yields a contradiction of Lemma~\ref{equal} if
$\mu(U^*(\hat D,\varepsilon,J))>0$. \fbox{} \\

\noindent {\bf Proof of Lemma~\ref{double-crossing}.}
Let us first consider the simpler case of ${\cal A}_k(v;\varepsilon,\varepsilon')$
in which $v \in J'_k$.
For $\varepsilon>0$, let
$J'_k(\varepsilon) \equiv J'_k \setminus \{ B(c'_k,\varepsilon) \cup B(d'_k,\varepsilon) \}$ and
$\hat D_k(\varepsilon) \equiv \hat D_k \setminus \{ B(c'_k,\varepsilon) \cup B(d'_k,\varepsilon) \}$.
With this notation, we have
\begin{equation} \label{split}
{\mathbb P}(\bigcup_{v \in J'_k} {\cal A}_k(v;\varepsilon,\varepsilon')) \leq
{\mathbb P}(\bigcup_{v \in J'_k(\varepsilon)} {\cal A}_k(v;\varepsilon,\varepsilon'))
+ F(\varepsilon; \hat D_k, c'_k, d'_k),
\end{equation}
where $F(\varepsilon; \hat D_k, c'_k, d'_k)$ is the probability that $\gamma^{\delta_k}_k$
enters $B(c'_k,\varepsilon)$ or $B(d'_k,\varepsilon)$ before touching $J'_k$.
$F(\varepsilon; \hat D_k, c'_k, d'_k)$ can be expressed as the sum of two crossing
probabilities in $\hat D_k(\varepsilon)$: (1) the probability of a blue crossing
from the (portion of the) counterclockwise arc $\overline{ac}$ of
$\partial\hat D_k(\varepsilon)$ to the ``first exposed" arc of $B(d',\varepsilon)$
contained in $\hat D_k(\varepsilon)$ and (2) the analogous probability of a yellow
crossing from the (portion of the) counterclockwise arc $\overline{da}$ of
$\partial\hat D_k(\varepsilon)$ to the ``first exposed" arc of $B(c',\varepsilon)$
contained in $\hat D_k(\varepsilon)$.
Since crossing probabilities converge to Cardy's formula, we easily conclude that
\begin{equation} \label{f-to-0}
\lim_{\varepsilon \to 0} \, \limsup_{k \to \infty} F((\varepsilon; \hat D_k, c'_k, d'_k)) = 0.
\end{equation}
Noting that the probability in the left hand side of~(\ref{split}) is nonincreasing
in $\varepsilon$, we see that in order to obtain
\begin{equation} \label{a-events}
\lim_{\varepsilon' \to 0} \, \limsup_{k \to \infty} \,
{\mathbb P}(\bigcup_{v \in J'_k} {\cal A}_k(v;\varepsilon,\varepsilon')) = 0,
\end{equation}
it is enough to show that the limit as $\varepsilon \to 0$ of the left hand side
of~(\ref{a-events}) is zero.
Therefore, thanks to~(\ref{split}) and~(\ref{f-to-0}),
it suffices to prove that
for every $\varepsilon>0$
\begin{equation}
\lim_{\varepsilon' \to 0} \, \limsup_{k \to \infty} \,
{\mathbb P}(\bigcup_{v \in J'_k(\varepsilon)} {\cal A}_k(v;\varepsilon,\varepsilon')) = 0.
\end{equation}

To do so, we follow the exploration process until time $T$, when it
first touches $\partial B(v,\varepsilon')$ for some
$v \in J'_k(\varepsilon)$,
and consider the annulus $B(v,\varepsilon) \setminus B(v,\varepsilon')$.
Let $\pi_Y$ be the leftmost yellow $\cal T$-path and $\pi_B$
the rightmost blue $\cal T$-path in $\Gamma(\gamma^{\delta_k}_k)$
at time $T$ that cross $B(v,\varepsilon) \setminus B(v,\varepsilon')$.
$\pi_Y$ and $\pi_B$ split the annulus
$B(v,\varepsilon) \setminus B(v,\varepsilon')$ into three sectors
that, for simplicity, we will call the central sector, containing
the crossing segment of the exploration path, the yellow (left) sector,
with $\pi_Y$ as part of its boundary, and the blue (right) sector,
the remaining one, with $\pi_B$ as part of its boundary.

We then look for a yellow ``lateral" crossing within the yellow sector
from $\pi_Y$ to $\partial \hat D_k$ and a blue lateral crossing within
the blue sector from $\pi_B$ to $\partial \hat D_k$.
Notice that the yellow sector may contain ``excursions"
of the exploration path coming off $\partial B(v,\varepsilon)$,
producing nested yellow and blue excursions off
$\partial B(v,\varepsilon)$, and the same for the blue sector.
But for topological reasons, those excursions are such that
for every group of nested excursions, the outermost one is
always yellow in the yellow sector and blue in the blue sector.
Therefore, by standard percolation theory arguments, the conditional
probability (conditioned on $\Gamma(\gamma^{\delta_k}_k)$ at time $T$)
to find a yellow lateral crossing of the yellow sector from $\pi_Y$
to $\partial \hat D_k$ is bounded below by the probability to find a
yellow circuit in an annulus with inner radius $\varepsilon'$
and outer radius $\varepsilon$.
An analogous statement holds for the conditional probability
(conditioned on $\Gamma(\gamma^{\delta_k}_k)$ at time $T$ and
on the entire percolation configuration in the yellow sector)
to find a blue lateral crossing of the blue sector from $\pi_B$
to $\partial \hat D_k$.
Thus for any fixed $\varepsilon>0$, by the Russo-Seymour-Welsh
lemma~\cite{russo,sewe}, the conditional probability to find both
a yellow lateral crossing within the yellow sector from $\pi_Y$ to
$\partial \hat D_k$ and a blue lateral crossing within the blue
sector from $\pi_B$ to $\partial \hat D_k$ goes to one as
$\varepsilon'\to 0$.

But if such yellow and blue crossings are present, the exploration
path is forced to touch $J'_k$ before exiting $B(v,\varepsilon)$,
and if that happens, the exploration process is stopped, so that
it will never exit $B(v,\varepsilon)$ and the union over
$v \in J'_k(\varepsilon)$
of ${\cal A}_k(v;\varepsilon,\varepsilon')$ cannot occur.
This concludes the proof of this case.

Let us now consider the remaining case in which $v \notin J'_k$.
The basic idea of the proof is then that by straightforward weak
convergence and related coupling arguments, the failure of~(\ref{double})
would imply that \emph{some} subsequence limit $\mu$ would satisfy
$\mu(U^{yellow}(\hat D,\varepsilon,J) \cup U^{blue}(\hat D,\varepsilon,J))>0$,
which would contradict Lemma~\ref{mushroom}.
This is essentially because the close approach of an exploration
path on the $\delta_k$-lattice to $J_k \setminus J'_k$ without
quickly touching nearby yields one two-sided colored $\cal T$-path
(the ``perimeter" of the portion of the hull of the exploration
path seen from a boundary point of close approach) and a one-sided
$\cal T$-path of the other color belonging to the percolation
cluster not seen from the boundary point (i.e., shielded by the
two-sided path).
Both the two-sided path and the one-sided one are subsets of
$\Gamma(\gamma^{\delta_k}_k)$.

\begin{figure}[!ht]
\begin{center}
\includegraphics[width=8cm]{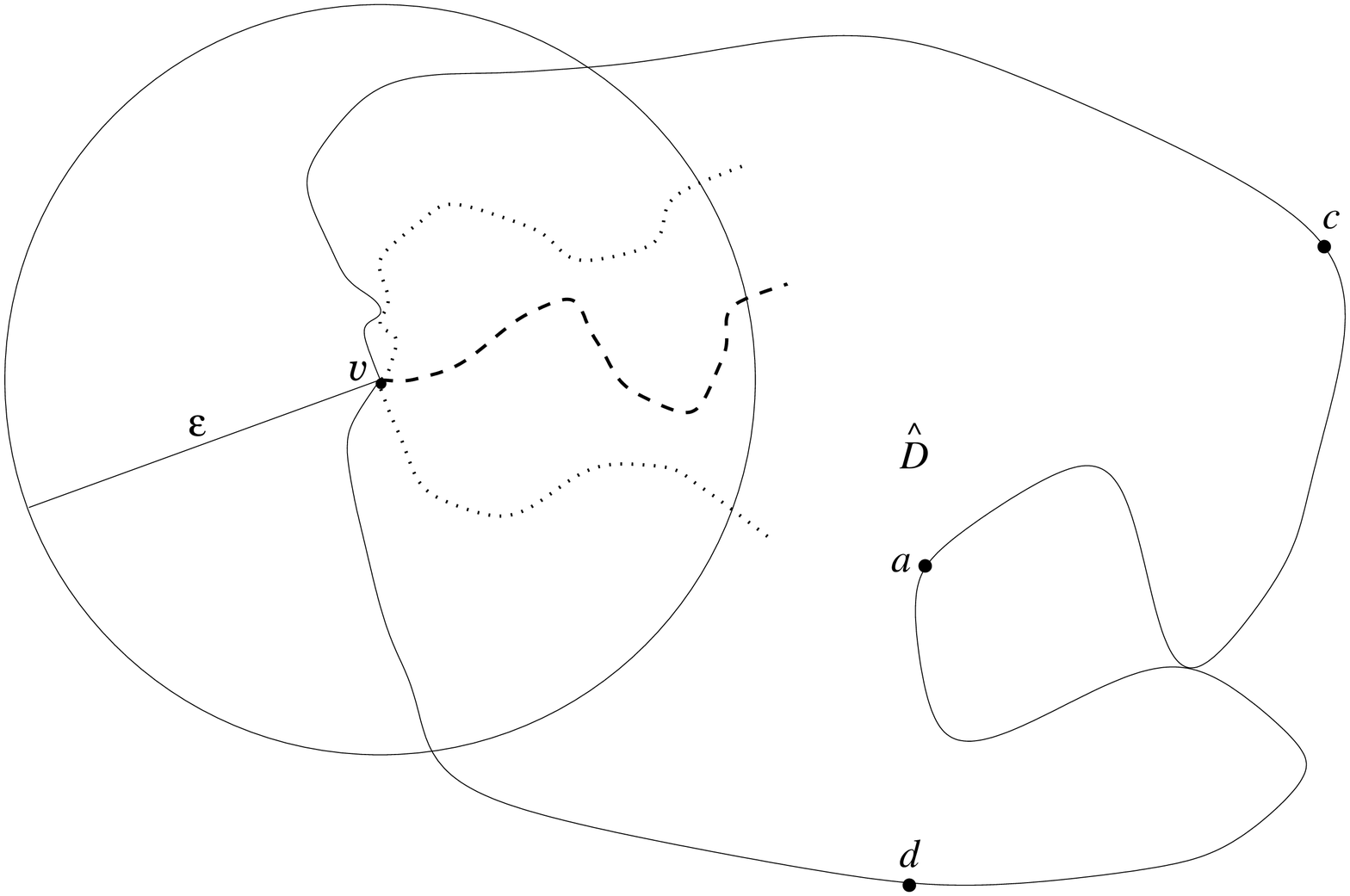}
\caption{The event consisting of a yellow double crossing
and a blue crossing used as a first step for obtaining a
blue mushroom event in the proof of Lemma~\ref{double-crossing}.
The dashed crossing is blue and the dotted crossing is yellow.}
\label{fig1-lemma7-2}
\end{center}
\end{figure}

We first note that since the probability in~(\ref{double}) is nonincreasing
in $\varepsilon$, we may assume that $\varepsilon < \min \{|a-c|,|a-d| \}$,
as requested by Lemma~\ref{mushroom}.
Assume by contradiction that~(\ref{double}) is false, so that close
encounters without touching happen with bounded away from zero probability.
Consider for concreteness an exploration path $\gamma^{\delta_k}_k$
that has a close approach to a point $v$ in the counterclockwise
arc $\overline{d'_k d_k}$.
The exploration path may have multiple close approaches to $v$ with
differing colors of the perimeter as seen from $v$, but for topological
reasons, the last time the exploration path comes close to $v$, it
must do so in such a way as to produce a yellow $\cal T$-path $\pi_Y$
(seen from $v$) that crosses $B(v,\varepsilon) \setminus B(v,\varepsilon')$
twice, and a blue path $\pi_B$ that crosses it once
(see Figure~\ref{fig1-lemma7-2}).
This is so because the exploration process that produced
$\gamma^{\delta_k}_k$ ended somewhere on $J'_k$ (and outside
$B(v,\varepsilon)$), which is to the right of (i.e., clockwise to) $v$.

%The presence of $\pi_Y$ implies that there are a yellow leftmost
%$\cal T$-path $\pi_L$ and a yellow rightmost $\cal T$-path $\pi_R$
%(looking at $v$ from inside $\hat D_k$) crossing the annulus
The portion of $\pi_Y$ inside $B(v,\varepsilon) \setminus B(v,\varepsilon')$
contains at least two yellow $\cal T$-paths crossing the annulus.
We denote by $\pi_L$ the leftmost (looking at $v$ from inside $\hat D_k$)
such path and denote by $\pi_R$ the rightmost one.
The paths $\pi_L$ and $\pi_R$ split
$B(v,\varepsilon) \setminus B(v,\varepsilon')$ into three sectors,
that we will call the central sector, containing $\pi_B$, the left
sector, with $\pi_L$ as part of its boundary, and the right sector,
with $\pi_R$ as part of its boundary.
Again for topological reasons, all other monochromatic crossings
of the annulus
associated with $\gamma^{\delta_k}_k$ are contained in the central sector,
including at least one blue path $\pi_B$.
As in the previous case, the left and right sectors can contain
nested monochromatic excursions off $\partial B(v,\varepsilon)$
(and in this case also excursions off $\partial B(v,\varepsilon')$),
but this time for every group of excursions, the outermost one is
yellow in both sectors.

Now consider the annulus $B(v,\varepsilon/3) \setminus B(v,\varepsilon/3)$.
We look for a yellow lateral crossing within the left sector
from $\pi_L$ to $\partial \hat D_k$ and a yellow lateral crossing
within the right sector from $\pi_R$ to $\partial \hat D_k$.
Since the outermost excursions in both sectors are yellow,
the conditional probability to find a yellow lateral crossing
within the left sector from $\pi_L$ to $\partial \hat D_k$ is bounded
below by the probability to find a yellow circuit in an annulus with
inner radius $\varepsilon/8$ and outer radius $\varepsilon/3$,
and an analogous statement holds for the conditional probability
to find a yellow lateral crossing within the right sector from
$\pi_R$ to $\partial \hat D_k$.
Thus for any fixed $\varepsilon>0$, by an application of the
Russo-Seymour-Welsh lemma~\cite{russo,sewe}, the conditional probability
to find both yellow lateral crossings
remains bounded away from zero as $k \to \infty$
and $\varepsilon' \to 0$.
But the presence of such yellow crossings would produce a (blue)
mushroom event, leading to a contradiction with Lemma~\ref{mushroom}. \fbox{} \\

\noindent {\bf Proof of Lemma~\ref{close-encounters}.} First of all, notice that
\begin{equation}
\bigcup_{v \in J_k \setminus J'_k}{\cal B}_k(v;\varepsilon,\varepsilon')
\subset \bigcup_{v \in J_k \setminus J'_k}{\cal C}_k(v;\varepsilon,\varepsilon')
\subset \bigcup_{v \in J_k \setminus J'_k}{\cal A}_k(v;\varepsilon,\varepsilon').
\end{equation}
By lemma~\ref{double-crossing}, we only need consider events in
$\bigcup_{v \in J_k \setminus J'_k} \{ {\cal C}_k(v;\varepsilon,\varepsilon')
\setminus {\cal B}_k(v;\varepsilon,\varepsilon') \}$, since
\begin{equation}
\lim_{\varepsilon' \to 0} \, \limsup_{k \to \infty} \,
{\mathbb P}(\bigcup_{v \in J_k \setminus J'_k} {\cal B}_k(v;\varepsilon,\varepsilon')) = 0.
\end{equation}
These are events such that $\gamma^{\delta_k}_k$ touches $\partial \hat D_k$
inside $B(v,\varepsilon)$, but it also contains a segment that stays within
$B(v,\varepsilon)$, does not touch $\partial \hat D_k$ (or any segment that touches
$\partial \hat D_k$ inside $B(v,\varepsilon)$), and has a double crossing of the
annulus $B(v,\varepsilon) \setminus B(v,\varepsilon')$.

Let us assume by contradiction that, for some fixed $\hat\varepsilon>0$,
\begin{equation} \label{not-close}
\limsup_{\varepsilon' \to 0} \, \limsup_{k \to \infty} \,
{\mathbb P}(\bigcup_{v \in J_k \setminus J'_k} {\cal C}_k(v;\hat\varepsilon,\varepsilon'))>0.
\end{equation}
This implies (by using a coupling argument) that we can find a subsequence $\{n\}$ of $\{k\}$
and sequences $\delta_n \to 0$ and $\varepsilon_n \to 0$ such that (with strictly positive
probability), as $n \to \infty$, $\gamma^{\delta_n}_n$ converges to a curve $\tilde\gamma$
that, for some $\bar v \in \partial \hat D$ (in fact, in $\overline{J \setminus J'}$), contains
a segment that stays in $B(\bar v,\hat\varepsilon)$, touches $\partial \hat D$ at $\bar v$, and
makes a double crossing of %$B(\bar v,\hat\varepsilon) \setminus B(\bar v,\varepsilon_n)$.
$B(\bar v,\hat\varepsilon) \setminus \{ \bar v \}$.
Since
the events that we are considering are in
$\bigcup_{v \in J_n \setminus J'_n} \{ {\cal A}_n(v;\hat\varepsilon,\varepsilon_n)
\setminus {\cal B}_n(v;\hat\varepsilon,\varepsilon_n) \}$, before the limit
$n \to \infty$ is taken, $\gamma^{\delta_n}_n$ has a segment
that stays in $B(\bar v,\hat\varepsilon)$, does not touch
$\partial \hat D_n$ and makes a double crossing of the annulus
$B(\bar v_n,\hat\varepsilon) \setminus B(\bar v_n,\varepsilon_n)$ for some
$\bar v_n \in \partial \hat D_n$ converging to $\bar v$.
Moreover, $\gamma^{\delta_n}_n$ must
touch $\partial \hat D_n$ inside $B(\bar v_n,\hat\varepsilon)$.

Consider first whether there exist $0<\tilde\varepsilon<\hat\varepsilon$ and
$n_0<\infty$ such that the point $v_n$, closest to $\bar v_n$, where $\gamma^{\delta_n}_n$
touches $\partial \hat D_n$ inside $B(\bar v_n,\hat\varepsilon)$ is outside
$B(\bar v,\tilde\varepsilon)$ for all $n \geq n_0$.
If this occurred with strictly positive probability, it would imply that for all
$\varepsilon<\tilde\varepsilon$,
\begin{equation}
\limsup_{n \to \infty} \,
{\mathbb P}(\bigcup_{v \in J_n \setminus J'_n} {\cal B}_n(v;\varepsilon,\varepsilon_n))>0,
\end{equation}
which would contradict~(\ref{double}).

\begin{figure}[!ht]
\begin{center}
\includegraphics[width=7cm]{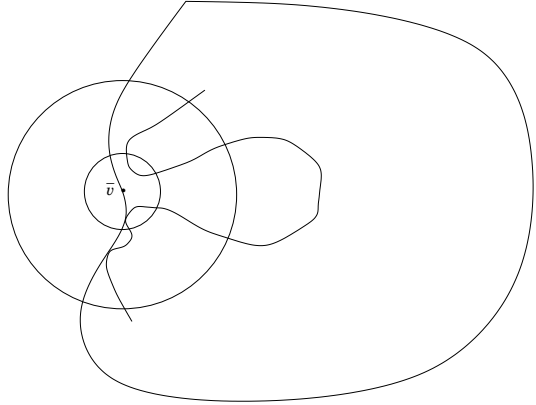}
\caption{One way in which a path inside a domain can have two double
crossings of an annulus centered at $\bar v$, one which does not
touch the boundary of the domain and one which does touch the boundary.}
\label{near-boundary1}
\end{center}
\end{figure}

\begin{figure}[!ht]
\begin{center}
\includegraphics[width=7cm]{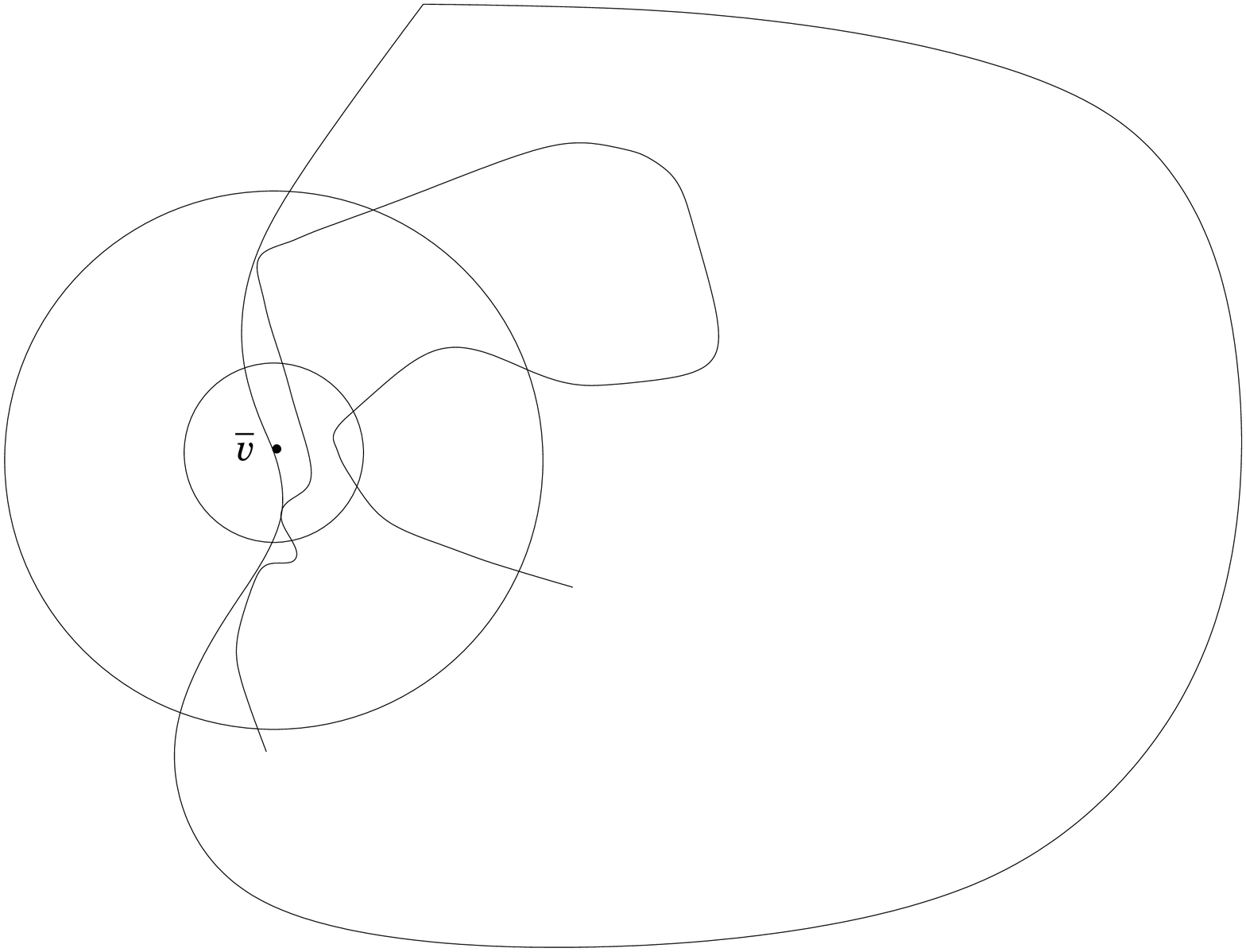}
\caption{The other (topologically distinct) way in which a path inside
a domain can have two double crossings of an annulus centered at $\bar v$,
one which does not touch the boundary of the domain and one which does
touch the boundary.}
\label{near-boundary2}
\end{center}
\end{figure}

\begin{figure}[!ht]
\begin{center}
\includegraphics[width=7cm]{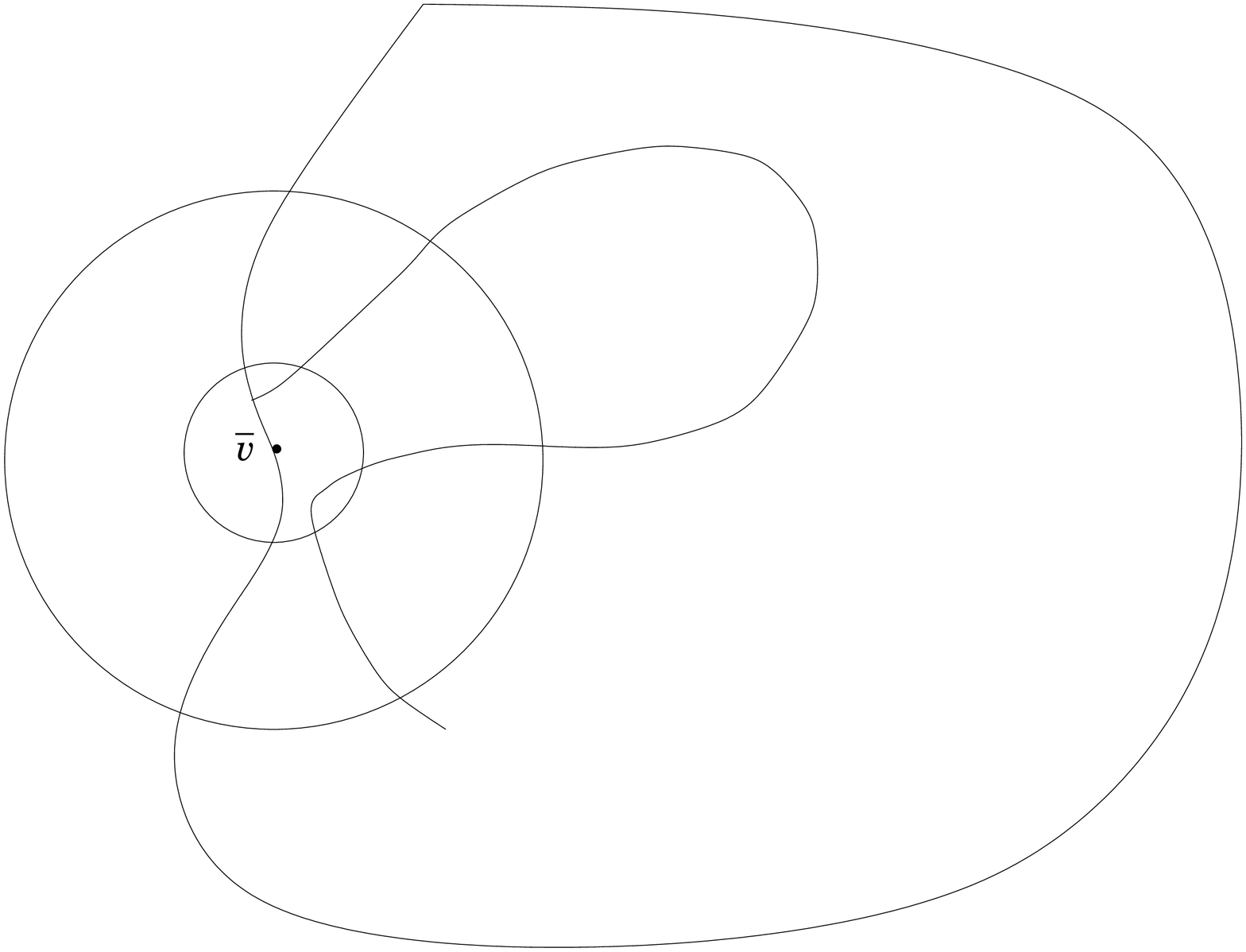}
\caption{A path in a domain making a double crossing of an annulus
centered at $\bar v$, without touching the domain boundary, and
a single crossing that ends on the boundary.}
\label{near-boundary3}
\end{center}
\end{figure}

This leaves the case in which the point $v_n$, closest to $\bar v_n$, where $\gamma^{\delta_n}_n$
touches $\partial \hat D_n$ inside $B(\bar v_n,\hat\varepsilon)$ converges to $\bar v$ (with strictly
positive probability).
Our assumption (\ref{not-close}) implies that, as $n \to \infty$, the segment of $\gamma^{\delta_n}_n$
that stays in $B(\bar v_n,\hat\varepsilon)$ and gets to within distance $\varepsilon_n$ of
$\partial \hat D_n$ without touching it, also does not touch the segment of $\gamma^{\delta_n}_n$
contained in $B(\bar v_n,\hat\varepsilon)$ that touches $\partial \hat D_n$ at $v_n$.
In that case, one can choose a subsequence $\{m\}$ of $\{n\}$ such that, for all $m$ large enough,
$|\bar v_m - v_m| \leq \varepsilon_m$ and there are five disjoint monochromatic crossings, not
all of the same color, of the annulus
$B(\bar v_m,\hat\varepsilon/2) \setminus B(\bar v_m,\varepsilon_m)$.
Three crossings of alternating colors are associated with the segment of $\gamma^{\delta_m}_m$ contained in $B(\bar v_m,\hat\varepsilon/2)$ which
does not touch $\partial\hat D_m$ and makes a double crossing of the annulus
$B(\bar v_m,\hat\varepsilon/2) \setminus B(\bar v_m,\varepsilon_m)$, while two
more crossings are associated with another segment which does touch $\partial\hat D_m$ and
makes either a single or a double crossing of the annulus (see Figures~\ref{near-boundary1},
\ref{near-boundary2} and~\ref{near-boundary3}).
This assures that $\gamma^{\delta_m}_m$ can have yet another segment that crosses the
annulus
$B(\bar v_m,\hat\varepsilon/2) \setminus B(\bar v_m,\varepsilon_m)$ (which would add at
least one more disjoint monochromatic crossing to the existing five) for at most finitely many $m$'s, since otherwise the standard six-arm bound proved in~\cite{ksz} would be
violated.
Therefore, we can choose $0<\bar\varepsilon<\hat\varepsilon/2$ such that for all $m$
large enough, $\gamma^{\delta_m}_m$
has at most two segments that stay in $B(\bar v_m,\bar\varepsilon)$ and cross
$B(\bar v_m,\bar\varepsilon) \setminus B(\bar v_m,\varepsilon_m)$, one of which does
not touch $\partial \hat D_m$ and makes a double crossing, while the other touches
$\partial\hat D_m$.

Since $\{m\}$ is a subsequence of $\{k\}$, it follows from~(\ref{not-close}) that
\begin{equation} \label{still-not-close}
%\limsup_{\varepsilon' \to 0} \, \limsup_{\substack{m \to \infty \\ \delta \to 0}} \,
%{\mathbb P}(\bigcup_{v \in J_m \setminus J'_m} {\cal C}^{\delta}_m(v;\bar\varepsilon,\varepsilon'))>0.
\limsup_{m \to \infty} \,
{\mathbb P}(\bigcup_{v \in J_m \setminus J'_m} {\cal C}_m(v;\bar\varepsilon,\varepsilon_m))>0.
\end{equation}
When the event ${\cal C}_m(v;\bar\varepsilon,\varepsilon_m)$ happens,
it follows from the above observations that $\gamma^{\delta_m}_m$ cannot touch
$\partial\hat D_m$ inside $B(\bar v_m,\bar\varepsilon)$ both before and after
the segment that makes a double crossing of
$B(\bar v_m,\bar\varepsilon) \setminus B(\bar v_m,\varepsilon_m)$ without touching $\partial \hat D_m$.
We will assume that, following the exploration path $\gamma^{\delta_m}_m$ from $a_m$ to
the target region $J'_m$, with positive probability $v_m$ is encountered after the segment
that makes a double crossing of $B(\bar v_m,\bar\varepsilon) \setminus B(\bar v_m,\varepsilon_m)$
without touching $\partial \hat D_m$.
To see that this assumption can be made without loss of generality, we note that we may consider
the case where the target area $J'_m$ is a single point $b_m$ and then consider the ``time-reversed"
exploration path in $\hat D_m$ that starts at $b_m$ and ends at $a_m$.

Now let $\gamma^{\delta_m}_m(t)$ be a parametrization of the exploration path $\gamma^{\delta_m}_m$
from $a_m$ to $J'_m$ and consider the stopping time $T$ (where the dependence on $m$ has been suppressed)
defined as the first time such that there is a point $v$ on $J_m$ with the property that
$\gamma^{\delta_m}_m$ up to time $T$ has never touched $\partial \hat D_m \cap B(v,\bar\varepsilon)$ and has
a segment that makes a double crossing of $B(v,\bar\varepsilon) \setminus B(v,\varepsilon_m)$.
It follows from~(\ref{still-not-close}) and the observation and assumption made right after
it that such a $T$ occurs with probability bounded away from zero as $m \to \infty$.

If this is the case, we can use the partial percolation configuration produced
inside $B(v,\bar\varepsilon) \setminus B(v,\varepsilon_m)$ by the exploration process
stopped at the stopping time $T$ to construct a mushroom event with positive probability,
as in the proof of Lemma~\ref{double-crossing}.
Notice, in fact, that the partial percolation configuration produced inside
$B(v,\bar\varepsilon) \setminus B(v,\varepsilon_m)$ by the exploration process stopped
at time $T$ has exactly the same properties as the partial percolation configuration of
the proof of Lemma~\ref{double-crossing} produced by the exploration process considered
there inside the annulus $B(v,\varepsilon) \setminus B(v,\varepsilon')$ -- see the second
part of the proof of Lemma~\ref{double-crossing}.
It follows that~(\ref{still-not-close}) implies positive probability of a mushroom event,
contradicting Lemma~\ref{mushroom}. \fbox{} \\

We are finally ready to prove the main result.

\begin{theorem} \label{conv-to-sle-1}
Let $(D,a,b)$ be a Jordan domain with two distinct selected points
on its boundary $\partial D$.
Then, for
%any sequence $\{(D^{\delta},a^{\delta},b^{\delta})\}$ of Jordan
Jordan sets $D^{\delta}$ from $\delta \cal H$
with two distinct selected e-vertices
$a^{\delta},b^{\delta}$ on their boundaries $\partial D^{\delta}$,
such that $(D^{\delta},a^{\delta},b^{\delta}) \to (D,a,b)$ as $\delta \to 0$,
the percolation exploration path
$\gamma^{\delta}_{D,a,b}$ inside $D^{\delta}$ from $a^{\delta}$ to $b^{\delta}$
converges in distribution to the trace $\gamma_{D,a,b}$ of chordal $SLE_6$
inside $D$ from $a$ to $b$, as $\delta \to 0$.
\end{theorem}

%\begin{theorem} \label{conv-to-sle-1}
%Consider a sequence CHANGE (01-20-06) $\delta_k \downarrow 0$ and a sequence
%$\{(D_k,a_k,b_k)\}$ of Jordan CHANGE (01-20-06) sets with two
%distinct selected CHANGE (01-20-06) e-vertices $a_k,b_k$ on their boundaries $\partial D_k$.
%Assume that $(D_k,a_k,b_k) \to (D,a,b)$, where $D$ is a Jordan domain
%with two distinct selected points on its boundary $\partial D$.
%Denote by CHANGE (01-20-06) $\gamma_k \equiv \gamma^{\delta_k}_{D_k,a_k,b_k}$
%the percolation exploration path inside CHANGE (01-20-06) $D_k$ from CHANGE (01-20-06)
%$a_k$ to CHANGE (01-20-06) $b_k$.
%Then, CHANGE (01-20-06) as $k \to \infty$, $\gamma_k$ converges in distribution to
%the trace $\gamma_{D,a,b}$ of chordal $SLE_6$ inside $D$ from $a$ to $b$.
%\end{theorem}

\noindent {\bf Proof.}
%It follows from~\cite{ab} that CHANGE (01-20-06) $\gamma_k$
%converges in distribution along subsequence limits $k_n$.
It follows from~\cite{ab} that $\gamma^{\delta}_{D,a,b}$ converges
in distribution along subsequence limits $\delta_k \downarrow 0$.
Since we have proved that the filling of any such subsequence
limit $\tilde\gamma$ satisfies the spatial Markov property
(Theorem~\ref{spatial-markov}) and the
hitting distribution of
$\tilde\gamma$ is determined by Cardy's formula (Theorem~\ref{strong-cardy}),
we can deduce from Theorem~\ref{characterization} that the limit
is unique and that the law of $\gamma^{\delta}_{D,a,b}$ converges,
as $\delta \to 0$, to the law of the trace $\gamma_{D,a,b}$ of
chordal $SLE_6$ inside $D$ from $a$ to $b$.~\fbox{}

\appendix
\refstepcounter{section}
\section*{Appendix \thesection: Sequences of Conformal Maps} \label{rado}

In this appendix, we give some results about sequences of conformal maps.
%that are used in various places throughout the paper.
A standard reference with more details is~\cite{pommerenke}.
%
%\begin{definition} \emph{(see Section~1.4 of~\cite{pommerenke})} \label{kernel-conv}
%Let $w_0 \in {\mathbb C}$ be given and let $\{ G_n \}$ be
%domains with $w_0 \in G_n \subset {\mathbb C}$.
%We say that $G_n$ converges to $G$ as $n \to \infty$ with respect to $w_0$
%(in the sense of {\bf kernel convergence}) if
%\begin{enumerate}
%\item either $G = \{ w_0 \}$, or else $G$ is a domain $\neq {\mathbb C}$
%with $w_0 \in G$ such that some neighborhood of every $w \in G$ lies
%in $G_n$ for large $n$; and
%\item for $w \in \partial G$ there exist $w_n \in \partial G_n$ such
%that $w_n \to w$ as $n \to \infty$.
%\end{enumerate}
%\end{definition}
%
%It is clear from the definition that every subsequence limit also converges
%to $G$ and it is also easy to see that the limit is uniquely determined.
%With this definition we can now state Carath\'eodory's kernel theorem~\cite{caratheodory}.
%
%\begin{theorem} \emph{(see Theorem~1.8 of~\cite{pommerenke})} \label{kernel-thm}
%Let $\phi_n$ map ${\mathbb D}$ conformally onto $G_n$ with $\phi_n(0)=w_0$
%and $\phi'_n(0)>0$.
%If $G = \{ w_0 \}$, let $\phi(z) \equiv w_0$; otherwise let $\phi$ map
%$\mathbb D$ conformally onto $G$ with $\phi(0)=w_0$ and $\phi'(0)>0$.
%Then, as $n \to \infty$, $\phi_n \to \phi$ locally uniformly in ${\mathbb D}$
%if and only if
%$G_n$ converges to $G$ with respect to $w_0$.
%\end{theorem}
%
There is a theorem attributed to both Courant~\cite{courant}
(see Theorem IX.14 of~\cite{tsuji}) and to Rad\'o~\cite{rado} (see Theorem~2.11
of~\cite{pommerenke}) that provides conditions under which conformal maps
$\phi_n$ from $\mathbb D$ onto Jordan domains $G_n$ converge uniformly on
all of $\overline{\mathbb D}$.
One of the purposes of this appendix is to provide in Corollary~\ref{cor-unif-conv}
an extension of the Courant-Rad\'o theorem to admissible domains (which are
not necessarily Jordan --- see the definition in Section~\ref{explo}).
Although this suffices here, it appears that there is a wider extension
(S.~Rohde, private communication) in which the domains $G_n$ have boundaries
given by continuous loops, without requiring admissibility.

To proceed, we need the next definition, in which a {\bf continuum} denotes
a compact connected set with more than one point.
\begin{definition} \emph{(Sect.~2.2 of~\cite{pommerenke})} \label{locally-connected}
The closed set $A \subset {\mathbb C}$ is called {\bf locally connected}
if for every $\varepsilon>0$ there is $\delta>0$ such that, for any two
points $a,b \in A$ with $|a-b|<\delta$, we can find a continuum $B$ with
diameter smaller than $\varepsilon$ and with $a,b \in B \subset A$.
\end{definition}
We remark that every continuous curve (with more than one point) is a
locally connected continuum (the converse is also true: every locally
connected continuum is a curve).
The concept of local connectedness gives a topological answer to the
problem of continuous extension of a conformal map to the domain boundary,
as follows.
\begin{theorem} \emph{(Sect.~2.1 of~\cite{pommerenke})}
\label{cont-thm}
Let $\phi$ map the unit disk $\mathbb D$ conformally onto $G \subset {\mathbb C} \cup \{\infty\}$.
Then $\phi$ has a continuous extension to $\overline{\mathbb D}$
if and only if $\partial G$ is locally connected.
\end{theorem}

When $\phi$ has a continuous extension to $\overline{\mathbb D}$,
we do not distinguish between $\phi$ and its extension.
This is the case for the conformal maps considered in this paper.
The problem of whether this extension is injective on $\overline{\mathbb D}$
has also a topological answer, as follows.
\begin{theorem} \emph{(Sect.~2.1 of~\cite{pommerenke})}
\label{cara-thm}
In the notation of Theorem~\ref{cont-thm}, the function $\phi$ has a continuous
and injective extension if and only if $\partial G$ is a Jordan curve.
\end{theorem}

%\begin{theorem} \emph{(see~\cite{pommerenke})} \label{continuous-boundary}
%Let $f$ map $\mathbb D$ conformally onto the bounded domain $G$.
%Then the following four conditions are equivalent:
%\begin{enumerate}
%\item $f$ has a continuous extension to $\overline{\mathbb D}$;
%\item $\partial G$ is a continuous curve;
%\item $\partial G$ is locally connected;
%\item ${\mathbb C} \setminus G$ is locally connected.
%\end{enumerate}
%\end{theorem}

When considering sequences of domains whose boundaries are locally connected
the following definition is useful.
\begin{definition} \emph{(Sect.~2.2 of~\cite{pommerenke})}
\label{unif-locally-connected}
The closed sets $A_n \subset {\mathbb C}$ are {\bf uniformly locally connected}
if, for every $\varepsilon>0$, there exists $\delta>0$ independent of $n$ such
that any two points $a_n,b_n \in A_n$ with $|a_n-b_n|<\delta$ can be joined by
continua $B_n \subset A_n$ of diameter smaller than $\varepsilon$.
\end{definition}

The convergence of domains used in this paper (i.e., $G_n \to G$ if
$\partial G_n \to \partial G$ in the uniform metric~(\ref{distance})
on continuous curves) allows us to use Carath\'eodory's kernel
convergence theorem (Theorem~1.8 of~\cite{pommerenke}).
However, we need uniform convergence in $\overline{\mathbb D}$.
This is guaranteed by the Courant-Rad\'o theorem in the case of Jordan
domains; in the non-Jordan case, sufficient conditions to have uniform
convergence are stated in the next theorem.

\begin{theorem} \emph{(Corollary~2.4 of~\cite{pommerenke})}
\label{unif-conv}
Let $\{ G_n \}$ be a sequence of bounded domains such that,
for some $0<r<R<\infty$, $B(0,r) \subset G_n \subset B(0,R)$
for all $n$ and such that $\{ {\mathbb C} \setminus G_n \}$
is uniformly locally connected.
Let $\phi_n$ map $\mathbb D$ conformally onto $G_n$ with
$\phi_n(0)=0$.
If $\phi_n(z) \to \phi(z)$ as $n \to \infty$ for each $z \in {\mathbb D}$,
then the convergence is uniform in $\overline{\mathbb D}$.
\end{theorem}

To use Theorem~\ref{unif-conv} we need the following lemma.
The definition of admissible and the related notion of
convergence are given respectively in Sect.~\ref{explo} and in
Theorem~\ref{strong-cardy}.

\begin{lemma} \label{lemma-unif-loc-conn}
Let $\{ (G_n,a_n,c_n,d_n) \}$ be a sequence of domains admissible
with respect to $(a_n,c_n,d_n)$ and assume that, as $n \to \infty$,
$(G_n,a_n,c_n,d_n) \to (G,a,c,d)$, where $G$ is a domain admissible
with respect to $(a,c,d)$.
Then the sequence of closed sets $\{ {\mathbb C} \setminus G_n \}$
is uniformly locally connected.
\end{lemma}

\noindent{\bf Proof.} %In order to prove the lemma, we claim that it
If the conclusion of the theorem is not valid, then for some $\varepsilon>0$,
there are indices $k$ (actually $n_k$, but we abuse notation a bit) and points
$u_k, v_k \in {\mathbb C} \setminus G_k$ with $|u_k-v_k| \to 0$ that cannot
be joined by a continuum of diameter $\leq \varepsilon$ in ${\mathbb C} \setminus G_k$.
We assume this and search for a contradiction.
By compactness, we may also assume that $u_k \to u, v_k \to v$, with $u=v$.
There is an easy contradiction (using a small disc around $u=v$ as the connecting
continuum) unless $u=v$ is on $\partial G$, and so we also assume that.
Further, by considering points on $\partial G_k$ near to $u_k, v_k$, we
can also assume that $u_k, v_k \in \partial G_k$.

Splitting $\partial G$ into three Jordan arcs,
$J_1=\overline{da}, J_2=\overline{ac}, J_3=\overline{cd}$, and $\partial G_k$
into the corresponding $J_{1,k}, J_{2,k}, J_{3,k}$, we note that there is an
easy contradiction (using arcs along $\partial G_k$ as the continua) if $u_k$
and $v_k$ both belong to $J_{1,k} \cup J_{3,k}$ or both belong to
$J_{2,k} \cup J_{3,k}$ for all $k$ large enough, since the concatenation of
$J_{1,k}$ with $J_{3,k}$ or of $J_{2,k}$ with $J_{3,k}$ is a Jordan arc.
%The above reasoning applies except when $u(=v)$ is on both $J_1$ and $J_2$.
%and (along a possible further subsequence) $u_k \in J_{1,k}$, $v_k \in J_{2,k}$
%(or vice versa) for all large $k$.
The above reasoning does not apply if $u=v$ is on both $J_1$ and $J_2$.
But when $u=v=a$, one can paste together small Jordan arcs on $J_{1,k}$
and $J_{2,k}$ to get a suitable continuum leading again to a contradiction.
The sole remaining case is when, for all $k$ large enough, $u_k$
belongs to the interior of $J_{1,k}$ and $v_k$ belongs to the interior of
$J_{2,k}$.

(Notice that we are ignoring the ``degenerate" case in which $c=d$
coincides with the ``last" [from $a$] double-point on $\partial G$,
and $J_3$ is a simple loop.
In that case $u_k$ and $v_k$ could converge to $u=v=c=d \in J_1 \cap J_2$
and $u_k$ or $v_k$ could still belong to $J_{3,k}$ for arbitrarily large $k$'s.
However, in that case one can find two distinct points on $J_3$, $c'$
and $d'$, such that $D$ is admissible with respect to $(a,c',d')$, and
points $c'_k$ and $d'_k$ on $J_{3,k}$ converging to $c'$ and $d'$
respectively, and define accordingly new Jordan arcs, $J'_1, J'_2, J'_3$
and $J'_{1,k}, J'_{2,k}, J'_{3,k}$, so that $u_k \in J'_{1,k}$ and
$v_k \in J'_{2,k}$ for $k$ large enough.
We assume that this has been done if necessary, and for simplicity
of notation drop the primes.)

In this case let $[u_k v_k]$ denote the closed straight line segment
in the plane between $u_k$ and $v_k$.
Imagine that $[u_k v_k]$ is oriented from $u_k$ to $v_k$ and let $v'_k$
be the first point of $J_{2,k}$ intersected by $[u_k v_k]$ and $u'_k$
be the previous intersection of $[u_k v_k]$ with $\partial G_k$.
Clearly, $u'_k \notin J_{2,k}$.
For $k$ large enough, $u'_k$ cannot belong to $J_{3,k}$ either, or
otherwise in the limit $k \to \infty$, $J_3$ would touch the interior
of $J_1$ and $J_2$.
We deduce that for all $k$ large enough, $u'_k \in J_{1,k}$.
Since $J_{1,k}$ and $J_{2,k}$ are continuous curves and therefore
locally connected, $u_k$ and $u'_k$ belong to a continuum $B_{1,k}$
contained in $J_{1,k}$ whose diameter goes to zero as $k \to \infty$,
and the same for $v_k$ and $v'_k$ (with $B_{1,k}$ and $J_{1,k}$
replaced by $B_{2,k}$ and $J_{2,k}$).

Since the interior of $[u'_k v'_k]$ does not intersect any
portion of $\partial G_k$, it is either contained in $G_k$
or in its complement ${\mathbb C} \setminus G_k$.
If $[u'_k v'_k] \subset {\mathbb C} \setminus G_k$, we have a
contradiction since the union of $[u'_k v'_k]$ with $B_{1,k}$
and $B_{2,k}$ is contained in ${\mathbb C} \setminus G_k$ and
is a continuum containing $u_k$ and $v_k$ whose diameter goes
to zero as $k \to \infty$.
%If the interior of $[u'_k v'_k]$ is contained in $G_k$, by the
%convergence of $G_k$ to $G$ and of $u'_k$ and $v'_k$ to $u=v$,
%we get that the $\liminf$ of the distance between a point $u$
%in the interior of $J_1$ and a point $v$ in the interior of $J_2$
%along a curve contained in $G$ is zero, contradicting
%Lemma~\ref{admissible} below.

If the interior of $[u'_k v'_k]$ is contained in $G_k$, let us
consider a conformal map $\phi_k$ from $\mathbb D$ onto $G_k$.
Since $\partial G_k$ is locally connected, the conformal map
$\phi_k$ extends continuously to the boundary of the unit disc.
Let $u'_k = \phi_k(u^*_k)$, $v'_k = \phi_k(v^*_k)$, $a_k = \phi_k(a^*_k)$,
$c_k = \phi_k(c^*_k)$ and $d_k = \phi_k(d^*_k)$.
The points $c^*_k, d^*_k, u^*_k, a^*_k, v^*_k$ are in counterclockwise
order on $\partial {\mathbb D}$, so that any curve in $\mathbb D$
from $a^*_k$ to the counterclockwise arc $\overline{c^*_k d^*_k}$
must cross the curve from $u^*_k$ to $v^*_k$ whose image under
$\phi_k$ is $[u'_k v'_k]$.
%is the line segment from $u'_k$ to $v'_k$.
This implies that any curve in $G_k$ going from $a_k$ to the
counterclockwise arc $\overline{c_k d_k}$ of $\partial G_k$ must
cross the (interior of the) line segment $[u'_k v'_k]$.
Then, in the limit $k \to \infty$, any curve in $G$ from $a$ to
the counterclockwise arc $\overline{cd}$ must contain the limit
point $u=\lim_{k \to \infty}u'_k=\lim_{k \to \infty}v'_k=v$.
On the other hand, except for its starting and ending point, any
such curve is completely contained in $G$, which implies that either
$u=v=a$ or else that (in the limit $k \to \infty$) the counterclockwise
arc $\overline{cd}$ is the single point at $u=v=c=d$.
We have already dealt with the former case.
In the latter case, one can paste together small Jordan arcs
from $u'_k$ to $d_k$, from $d_k$ to $c_k$, and from $c_k$ to $v'_k$,
and take the union with $B_{1,k}$ and $B_{2,k}$ (defined above) to
get a suitable continuum in ${\mathbb C} \setminus G_k$
containing $u_k$ and $v_k$, leading to a contradiction.
This concludes the proof. \fbox{} \\

%Let the {\bf internal distance} between two distinct points (or more
%precisely, two prime ends) $u,v$ on the boundary $\partial G$ of the
%domain $G$ be the $\inf$ over all paths in $G$ from $u$ to $v$ of the
%path arclength.
%Then the following holds.
%\begin{lemma} \label{admissible}
%Given a domain $(G,a,c,d)$ admissible with respect to $(a,c,d)$,
%let $J_1 = \overline{da}$, $J_2 =\overline{ac}$ and $J_3 = \overline{cd}$.
%Then, the internal distance between any two points (or prime ends) in the
%interiors of $J_1$ and $J_2$ respectively is strictly positive.
%\end{lemma}
%
%\noindent {\bf Proof.} Take $u$ and $v$ in the interior of
%$J_1$ and $J_2$ respectively.
%If $u \neq v$ (as points in the plane), the conclusion is obvious,
%so we just have to consider the case in which $u=v$.
%Assume, by contradiction, that the conclusion is false in that
%case, and consider a conformal map $\phi$ that maps $G$ onto
%the unit disc $\mathbb D$.
%Since $\partial G$ is locally connected, $\phi$ can be extended
%continuously to the boundary of $G$.
%If $u=v$, then $u$ and $v$ belong to two different prime ends and
%are mapped by $\phi$ to two different points $\phi(u) \neq \phi(v)$
%on $\partial {\mathbb D}$.
%If the internal distance between $u$ and $v$ were zero, there would
%be a sequence of points $w_n$ in $G$ (along a sequence of paths with
%arclength tending to zero) such that $|w_n - u|, |w_n - v| \to 0$.
%But then $\phi(w_n)$ would be converging simultaneously to $\phi(u)$
%and $\phi(v)$, which gives a contradiction. \fbox{} \\

Theorem~\ref{unif-conv}, together with Carath\'eodory's kernel
convergence theorem (Theorem~1.8 of~\cite{pommerenke}) and
Lemma~\ref{lemma-unif-loc-conn}, implies the following result.

\begin{corollary} \label{cor-unif-conv}
%Let $\{ G_n \}$ be a sequence of bounded simply connected domains
%(containing the origin) whose boundaries $\{ \partial G_n \}$ are continuous
%curves converging, as $n \to \infty$, to the continuous curve $\partial G$
%(boundary of the bounded simply connected domain $G$ containing the origin)
%in the uniform metric~(\ref{distance}) on continuous curves.
With the notation and assumptions of Lemma~\ref{lemma-unif-loc-conn}
(and also assuming that $G_n$ and $G$ contain the origin), let $\phi_n$
map ${\mathbb D}$ conformally onto $G_n$ with $\phi_n(0)=0$ and $\phi'_n(0)>0$,
and $\phi$ map $\mathbb D$ conformally onto $G$ with $\phi(0)=0$ and $\phi'(0)>0$.
Then, as $n \to \infty$, $\phi_n \to \phi$ uniformly in $\overline{\mathbb D}$.
\end{corollary}

\noindent {\bf Proof.} As already remarked, the convergence of $\partial G_n$
to $\partial G$ in the uniform metric~(\ref{distance}) on continuous curves
(which is part of the definition of $(G_n,a_n,c_n,d_n) \to (G,a,c,d)$)
easily implies that the conditions in Carath\'eodory's kernel theorem
(Theorem~8.1 of~\cite{pommerenke}) are satisfied and therefore that $\phi_n$
converges to $\phi$ locally uniformly in $\mathbb D$, as $n \to \infty$.
By an application of Lemma~\ref{lemma-unif-loc-conn}, the sequence
$\{ {\mathbb C} \setminus D_n \}$ is uniformly locally connected, so that
we can apply Theorem~\ref{unif-conv} to conclude that, as $n \to \infty$,
$\phi_n$ converges to $\phi$ uniformly in $\overline{\mathbb D}$. \fbox{} \\

The next result is a corollary of the previous one and is used in the
proof of Theorem~\ref{spatial-markov}.

\begin{corollary} \label{cor2-unif-conv}
Let $\{ (G_n,a_n,c_n,d_n) \}$ be a sequence of domains admissible
with respect to $(a_n,c_n,d_n)$ with $b_n$ in the interior of $\overline{c_n d_n}$,
and assume that, as $n \to \infty$, $(G_n,a_n,c_n,d_n) \to (G,a,c,d)$
and $b_n \to b$, where $G$ is a domain admissible with respect to $(a,c,d)$
and $b$ is in the interior of $\overline{cd}$.
Let $f$ be a conformal map from $\mathbb H$ to $G$ such that $f^{-1}(a)=0$
and $f^{-1}(b)=\infty$.
Then, there exists a sequence $\{f_n\}$ of conformal maps from $\mathbb H$
to $G_n$ with $f_n^{-1}(a_n)=0$ and $f_n^{-1}(b_n)=\infty$ and such that
$f_n$ converges to $f$ uniformly
in $\overline{\mathbb H}$.
\end{corollary}

\noindent {\bf Proof.}
The conformal transformation $f(\cdot)$ can be written as $\phi \circ \psi(\lambda \cdot)$,
where $\lambda$ is a positive constant, $\phi$ is the unique conformal
transformation that maps $\mathbb D$ onto $G$ with $\phi(0)=0$ and $\phi'(0)>0$
(we are assuming for simplicity that $G$ contains the origin; if that is not
the case, one can use a translated domain that does contain the origin),
and
\begin{equation}
\psi(z) = e^{i \theta_0} \left( \frac{(z+1)-z_0}{(z+1)-\overline{z_0}} \right)
\end{equation}
maps $\overline{\mathbb H}$ onto $\overline{\mathbb D}$.
$\theta_0$ is chosen so that $e^{i \theta_0} = \phi^{-1}(b)$ and $z_0$
so that $|1-z_0|=1$, $\text{Im}(z_0)>0$ and
$\phi^{-1}(b) (\frac{1-z_0}{1-\overline{z_0}}) = \phi^{-1}(a)$,
which implies that $f^{-1}$ indeed maps $a$ to $0$ and $b$ to $\infty$.

We now take $f_n(\cdot)$ of the form $\phi_n \circ \psi_n(\lambda\cdot)$,
where $\phi_n$ is the unique conformal transformation that maps $\mathbb D$
onto $G_n$ with $\phi_n(0)=0$ and $\phi_n'(0)>0$
(the assumption that $G$ contains the origin implies that $G_n$ also contains
the origin, for large $n$),
%(we are again assuming for simplicity that $G$ -- and hence $G_n$ for large $n$
%-- contains the origin; if that is not the case, one can use a translated domain
%that does contain the origin),
and
\begin{equation}
\psi_n(z) = e^{i \theta_n} \left( \frac{(z+1)-z_n}{(z+1)-\overline{z_n}} \right)
\end{equation}
maps $\overline{\mathbb H}$ onto $\overline{\mathbb D}$.
Note that the same $\lambda$ is used as appeared in the expression
$\phi \circ \psi(\lambda\cdot)$ for $f(\cdot)$.
$\theta_n$ is chosen so that $e^{i \theta_n} = \phi_n^{-1}(b_n)$ and $z_n$
is chosen so that $|1-z_n|=1$, $\text{Im}(z_n)>0$ and
$\phi_n^{-1}(b_n) (\frac{1-z_n}{1-\overline{z_n}}) = \phi_n^{-1}(a_n)$,
which implies that $f_n^{-1}$ indeed maps $a_n$ to $0$ and $b_n$ to $\infty$.

Corollary~\ref{cor-unif-conv} implies that, as $n \to \infty$, $\phi_n$
converges to $\phi$ uniformly in $\overline{\mathbb D}$.
This, together with the convergence of $a_n$ to $a$ and $b_n$ to $b$, implies
that $\phi_n^{-1}(a_n)$ converges to $\phi^{-1}(a)$ and $\phi_n^{-1}(b_n)$ to
$\phi^{-1}(b)$.
Therefore, we also have the convergence of $\psi_n$ to $\psi$ uniformly
in $\overline{\mathbb H}$, which implies that $f_n$ converges to $f$ uniformly
in $\overline{\mathbb H}$.
\fbox{} \\

We conclude this appendix with a simple lemma, used in the proof of
Theorem~\ref{strong-cardy}, about the continuity of Cardy's formula
with respect to the shape of the domain and the positions of the four
points on the boundary.

\begin{lemma} \label{cont-cardy}
For $\{ (D_n,a_n,c_n,b_n,d_n) \}$ and $(D,a,c,b,d)$ as in
Theorem~\ref{strong-cardy}, let $\Phi_n$ denote Cardy's formula
(see~(\ref{cardy-formula})) for a crossing inside $D_n$ from the
counterclockwise segment $\overline{a_n c_n}$ of $\partial D_n$ to
the counterclockwise segment $\overline{b_n d_n}$ of $\partial D_n$
and $\Phi$ the corresponding Cardy's formula for the limiting domain $D$.
Then, as $n \to \infty$, $\Phi_n \to \Phi$.
\end{lemma}

\noindent {\bf Proof.} Let $\phi_n$ be the conformal map that takes
$\mathbb D$ onto $D_n$ with $\phi_n(0)=0$ and $\phi'_n(0)>0$, and let
$\phi$ denote the conformal map from $\mathbb D$ onto $D$ with $\phi(0)=0$
and $\phi'(0)>0$; let $z_1 = \phi^{-1}(a)$, $z_2 = \phi^{-1}(c)$, $z_3 = \phi^{-1}(b)$,
$z_4 = \phi^{-1}(d)$, $z_1^n = \phi_n^{-1}(a_n)$, $z_2^n = \phi_n^{-1}(c_n)$,
$z_3^n = \phi_n^{-1}(b_n)$, and $z_4^n = \phi_n^{-1}(d_n)$.
We can apply Corollary~\ref{cor-unif-conv} to conclude that, as $n \to \infty$,
$\phi_n$ converges to $\phi$ uniformly in $\overline{\mathbb D}$.
This, in turn, implies that, as $n \to \infty$, $z_1^n \to z_1$,
$z_2^n \to z_2$, $z_3^n \to z_3$ and $z_4^n \to z_4$.

Cardy's formula for a crossing inside $D_n$ from the counterclockwise
segment $\overline{a_n c_n}$ of $\partial D_n$ to the counterclockwise
segment $\overline{b_n d_n}$ of $\partial D_n$ is given by
\begin{equation} \label{cardy1}
\Phi_n =
\frac{\Gamma(2/3)}{\Gamma(4/3) \Gamma(1/3)} \eta^{1/3}_n
{}_2F_1(1/3,2/3;4/3;\eta_n),
\end{equation}
where
\begin{equation}
\eta_n =
\frac{(z_1^n-z_2^n)(z_3^n-z_4^n)}{(z_1^n-z_3^n)(z_2^n-z_4^n)}.
\end{equation}
Because of the continuity of $\eta_n$ in $z_1^n$, $z_2^n$, $z_3^n$,
$z_4^n$, and the continuity of Cardy's formula~(\ref{cardy1}) in $\eta_n$,
the convergence of $z_1^n \to z_1$, $z_2^n \to z_2$, $z_3^n \to z_3$ and
$z_4^n \to z_4$ immediately implies the convergence of $\Phi_n$ to $\Phi$. \fbox{}

\bigskip
\bigskip

\noindent {\bf Acknowledgements.} We are grateful to Vincent
Beffara, Greg Lawler, Oded Schramm, Yuri Suhov and Wendelin Werner
for various interesting and useful conversations and to Stas
Smirnov for communications about work in progress.
We note that a discussion with Oded Schramm,
%at the November 2004 Northeast Probability Seminar at the CUNY Graduate Center,
about dependence of exploration paths with respect to small
changes of domain boundaries, pointed us in a direction that
eventually led to Lemmas~\ref{double-crossing}-\ref{mushroom}.
We are especially grateful to Vincent Beffara for pointing out
a gap in a preliminary version of the proof of the main result.
We thank Michael Aizenman, Oded Schramm, Vladas Sidoravicius,
and Lai-Sang Young for comments about presentation, Alain-Sol
Sznitman for his interest and encouragement, Steffen Rohde for
extensive discussions about extensions of Corollary~\ref{cor-unif-conv},
and an anonymous referee for useful comments.
F.~C. thanks Wendelin Werner for an invitation to Universit\'e
Paris-Sud, and acknowledges the kind hospitality of the
Courant Institute where part of this work was completed.
C.~M.~N. acknowledges the kind hospitality of the Vrije
Universiteit Amsterdam.

\end{document}